\documentclass[11pt,reqno]{amsart} 
\usepackage{graphicx}
\usepackage{amssymb,amsthm,listings,enumerate,geometry,setspace,float,csquotes,mathtools,url}
\usepackage{etoolbox}
\usepackage[labelformat=parens,labelfont={}]{subcaption}
\usepackage{hyperref}
\usepackage{pdfpages,xspace,xparse}
\usepackage{import,pst-plot}
\usepackage[numbers,sort]{natbib}
\usepackage{bm}
\usepackage[symbol,bottom]{footmisc}

\DeclareMathOperator{\sech}{sech}

\DeclareRobustCommand{\bms}[1]{\bm{#1}}
\newcommand{\MATLAB}{\textsc{Matlab}\xspace}

\def\Xint#1{\mathchoice
   {\XXint\displaystyle\textstyle{#1}}%
   {\XXint\textstyle\scriptstyle{#1}}%
   {\XXint\scriptstyle\scriptscriptstyle{#1}}%
   {\XXint\scriptscriptstyle\scriptscriptstyle{#1}}%
   \!\int}
\def\XXint#1#2#3{{\setbox0=\hbox{$#1{#2#3}{\int}$}
     \vcenter{\hbox{$#2#3$}}\kern-.5\wd0}}
\def\dashint{\Xint-}

\theoremstyle{definition}
\newtheorem{remark}{Remark}[section]
\numberwithin{equation}{section}
\numberwithin{figure}{section}

\makeatletter
\def\@settitle{\begin{center}%
  \baselineskip14\p@\relax
  \@title
  \end{center}%
}

\title[Num. soln. of SD lin. problems on the finite interval using UTM]{\Large{The numerical solution of semidiscrete linear evolution problems on the finite interval using the Unified Transform Method}\vspace{-15pt}}

\author[J. Cisneros \& B. Deconinck]{}

\begin{document}
\maketitle
\begin{center}
	Jorge Cisneros\footnote[1]{Corresponding author: \texttt{jorgec5@uw.edu}} \& Bernard Deconinck \\[5pt]
	Department of Applied Mathematics \\
	University of Washington\\
	Seattle, WA 98195-2420\\[5pt]
	\today
\end{center}

\begin{abstract}
We study a semidiscrete analogue of the Unified Transform Method introduced by A. S. Fokas, to solve initial-boundary-value problems for linear evolution partial differential equations with constant coefficients on the finite interval $x \in (0,L)$. The semidiscrete method is applied to various spatial discretizations of several first and second-order linear equations, producing the exact solution for the semidiscrete problem, given appropriate initial and boundary data. From these solutions, we derive alternative series representations that are better suited for numerical computations. In addition, we show how the Unified Transform Method treats derivative boundary conditions and ghost points introduced by the choice of discretization stencil and we propose the notion of ``natural'' discretizations. We consider the continuum limit of the semidiscrete solutions and compare with standard finite-difference schemes.\\[7pt]
\fontsize{8}{9.6}\selectfont{Keywords: continuum limit, finite difference, finite interval, ghost points, semidiscrete linear problem, Unified Transform Method}
\end{abstract}



\section{Introduction}

	We are interested in the numerical solution of the $M^{\text{th}}$-order quasilinear partial differential equation (PDE)
	\begin{equation}
		q_t = c \, q_{Mx} + F\left(q,q_x,\ldots,q_{(M-1)x} \right),\quad c \in \mathbb{C} \, \text{\textbackslash} \, \{ 0 \},
		\label{ibvp_eq}
	\end{equation} 
	on the finite interval $x \in (0,L)$. We assume initial and boundary conditions have been provided to make the initial-boundary-value problem (IBVP) well posed, and that the provided initial and boundary functions are sufficiently smooth and compatible.
	
\begin{figure}[b]
	\begin{center}
	\def\svgwidth{3.75in}
		\includegraphics[width=0.5\linewidth]{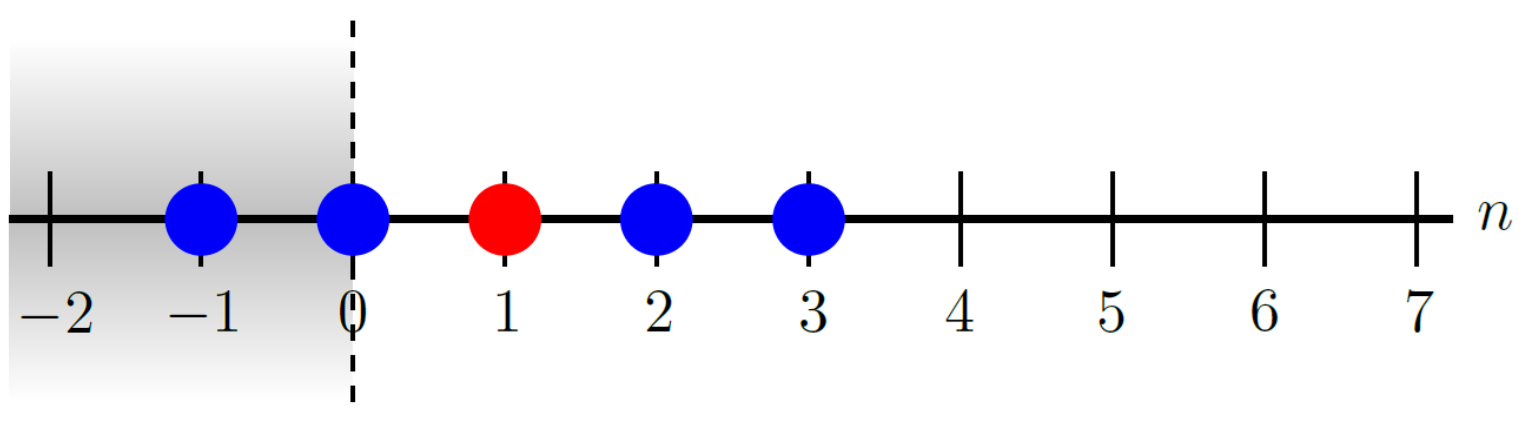}
		\caption{A stencil that requires information at a ghost point near $x = 0$.}
		\label{ghost_point}
	\end{center}
\end{figure}

	The application of finite-difference schemes on discretized spatial and temporal grids is a standard and intuitive approach to numerically solve finite-interval IBVPs at points $x_n \equiv n \Delta x$ and $t_j \equiv j \Delta t$. However, numerical methods with high-order spatial stencils often need data at grid points outside of the $x$ domain, referred to as \textit{ghost points} (see Figure \ref{ghost_point}). 
These points arise solely from the choice of spatial stencil, independent from the given IBVP. Note that with periodic boundary conditions, ghost points and their complications are nonexistent. The current heuristic techniques to determine data at ghost points can destabilize schemes shown to be stable in the full-line or periodic problem \cite{randy,dale} and are not easily generalizable to higher-order PDEs and their derivative boundary conditions \cite{gustafsson, cheema_thesis, coco, weno, black_hole, fd_finance, baronas1, baronas2}. Incorporating all boundary conditions correctly while addressing ghost points is a non-trivial issue \cite{iserles,strikwerda,trefethen,fd_pde_book} that we propose to undertake with an operator splitting approach. 
	
	Higher-order derivatives usually require higher-order stencils that embed more ghost points than lower-order stencils. We consider the class of quasilinear IBVPs \eqref{ibvp_eq} that are discretized in space $x_n = n h$ but are continuous in time $t$, \textit{i.e.}, semidiscretized, where the most nonlocal stencil is applied to the linear term $c \, q_{Mx}$. When solving the nonlinear problem \eqref{ibvp_eq}, the lower-order problem $q_t = F\left(q,q_x,\ldots,q_{(M-1)x} \right)$ is handled with established split-step techniques, while the linear problem $q_t = c \, q_{Mx}$ and its ghost points require special attention. Similar to the ideas presented for half-line problems \cite{SDUTM_HL}, this paper treats $q_t = c \, q_{Mx}$ on a finite interval using the Unified Transform Method on the semidiscrete $(n,t)$-plane. The split-step methods for the nonlinear IBVPs \eqref{ibvp_eq} will be explored in a future paper.


\section{The Continuous Unified Transform Method}\label{cont_UTM_steps}

	The Unified Transform Method (UTM) or Method of Fokas provides a systematic procedure to solve linear constant-coefficient evolution PDEs on the half-line and finite interval, although it originated with integrable nonlinear PDEs \cite{fokas_paper,fokas_collab,fokas_book}. For problems $q_t = c \, q_{Mx}$ with nonhomogeneous boundary conditions, the UTM gives an explicit analytical solution in terms of integrals along paths in the complex plane of a spectral parameter $k \in \mathbb{C}$ that can be numerically evaluated through contour parameterizations \cite{flyer,xin_paper,papa_heat}. Regardless of the type of boundary conditions, \textit{e.g.}, Dirichlet, Neumann, or Robin, the UTM is more widely applicable than classical methods to solve evolution IBVPs \cite{kernel_bcs}.
	
	The UTM for either half-line or finite-interval problems uses the following, algorithmic steps \cite{bernard_fokas,SDUTM_HL}: 
\begin{enumerate}
	\item rewrite the PDE in divergence form with a spectral parameter $k$ to obtain the local relation and the dispersion relation $\tilde{W}(k)$,\label{utm_1}
	\item integrate over the $(x,T)$-domain to obtain the global relation, \label{utm_2}
	\item invert the global relation to obtain a representation of the solution that depends on known and unknown boundary data, \label{utm_3}
	\item determine symmetries $\tilde{\nu}_j(k)$ of $\tilde{W}(k)$, \label{utm_4}
	\item determine where in $\mathbb{C}$ the global relations are valid when evaluated at $\tilde{\nu}_j(k)$, \label{utm_5}
	\item if necessary, appropriately deform integral paths involving boundary terms, \label{utm_6}
	\item solve for unknown boundary data using the global relations evaluated at $\tilde{\nu}_j(k)$, and \label{utm_7}
	\item check that integral terms involving $\hat{q}(\tilde{\nu}_j,T)$ vanish, resulting in a solution representation. \label{utm_8}
\end{enumerate}
	Although the calculations are more tedious for higher-order problems and boundary conditions, the UTM ultimately reduces the challenge of solving an IBVP to solving a set of algebraic equations involving the dispersion relation and its symmetries.
	
	For the finite-difference approach to nonlinear IBVPs, we apply the UTM to linear semidiscrete IBVPs. A method-of-lines formulation allows the UTM to address ghost points directly by providing an analytical solution to the semidiscrete problem.


\section{Semidiscrete UTM: Notation and Definitions}

	The UTM has received a lot of attention for continuous IBVPs, but not nearly as much for semidiscrete ones. Biondini and Hwang \cite{gino_nls}, Biondini and Wang \cite{gino_main}, and Moon and Hwang \cite{sd_finite} study semidiscrete problems in the context of the UTM, but from the perspective of a purely semidiscrete problem on closed contours with discretized Lax pairs. In \cite{SDUTM_HL}, we study semidiscrete problems on the half line, motivated by spatial discretizations of PDEs. In the present paper, we do the same for discretizations of PDEs on the finite interval.
	
	Our motivation is to improve the numerical solution of nonlinear IBVPs. To reach our goal, we apply the SD-UTM to the linear problem, addressing ghost points in finite-difference methods directly and analytically. As we will show, the SD-UTM formulas for the semidiscrete solution $q_n(T)$ are easier to work with than those obtained from the continuous UTM, but further approximations are needed in order to efficiently implement them into a split-step method (see Section \ref{small_time_sec}).
	
	In the following sections, we present the SD-UTM through examples of several linear semidiscretized IBVPs on the finite interval. We introduce an explicit mesh parameter $h = L/(N+1) \ll 1$ for $n = 0, 1, \ldots, N, N+1$ with $N$ interior-domain grid points, such that $n = 0$ and $n = N+1$ correspond to $x = 0$ and $x = L$, respectively. For each section, the first few concrete examples are followed by higher-order discretizations where ghost points play a significant role. We follow a similar procedure to Steps \eqref{utm_1}-\eqref{utm_8}, adapting the methodology as we move through the examples. We use the shift operator $\Delta Q_n = Q_{n+1} - Q_n$, which effectively replaces the spatial derivative with a forward difference. For IBVPs with Dirichlet boundary conditions on both ends of the interval, the Fourier transform pair can be written as
	\begin{subequations}
		\begin{equation}
			\hat{q}(k,t) = h \sum_{n = 1}^{N} e^{-ik nh} q_n(t) , \quad\quad k \in \mathbb{C},
			\label{fourier_cont_SD}
		\end{equation}
		\begin{equation}
			q_n(t) = \frac{1}{2 \pi} \int_{-\pi/h}^{\pi/h} e^{iknh} \hat{q}(k,t) \, dk, \quad\quad k \in \mathbb{C}.
			\label{inv_fourier_cont_SD}
		\end{equation}
		\label{fourier_transforms_SD}
	\end{subequations}
	If a Dirichlet condition is not given at $x = 0$, then we start \eqref{fourier_cont_SD} at $n = 0$. Similarly, if there is no Dirichlet data at $x = L$, the sum ends at $N +1$. We define the time transforms of spatial \textit{nodes} at and near the $n = 0$ boundary, including ghost points: 
	\begin{equation}
		f_j(W,T) = \int_{0}^T e^{Wt} q_j(t) \, dt, \quad j = \ldots -1, 0, 1, \ldots, \quad k \in \mathbb{C}, 
		\label{F}
	\end{equation}
	and, at and near the $n = N+1$ boundary, including ghost points: 
	\begin{equation}
		g_j(W,T) = \int_{0}^T e^{Wt} q_{N+1+j}(t) \, dt, \quad j = \ldots -1, 0, 1, \ldots, \quad k \in \mathbb{C},
		\label{G}
	\end{equation}
	for an arbitrary finite $T > 0$ and semidiscrete dispersion relation $W(k)$.
	
	To compare with standard numerical approaches, finite-difference solutions are obtained using a method-of-lines approach with the same spatial discretization as the semidiscrete UTM approach, so that we set up the full discretization of the PDE into a large, but sparse, system of ordinary differential equations (ODEs). For example, the centered-discretized heat equation with Dirichlet boundary conditions is formulated as
	\begin{align*}
		q_t(x,t) = q_{xx}(x,t) \quad
		\Rightarrow \quad \dot{q}_n(t) = \frac{q_{n+1}(t) - 2 q_{n}(t) + q_{n-1}(t)}{h^2} \quad
		\Rightarrow \quad \dot{Q}(t) = A Q(t) + g(t) + p(t), 
	\end{align*}
	where $Q(t) \in \mathbb{R}^{N}$ is a column vector composed of all $q_n(t)$, $A \in \mathbb{R}^{N \times N}$ is a sparse tridiagonal matrix, and both $g(t) \in \mathbb{R}^{N}$ and $p(t) \in \mathbb{R}^{N}$ are sparse column vectors that include boundary conditions from the left and right, respectively. The system of ODEs above is discretized in time and solved via the Forward Euler (FE), Runge-Kutta Fourth Order (RK4), Backward Euler (BE), and Trapezoidal (TR) methods. In summary, we compare the exact solution of an IBVP to these four classical numerical solutions and the SD-UTM explicit solution, which exactly solves the spatially discretized problem, \textit{i.e.}, requires no time-stepping.
	
	For the following examples, we compute the SD-UTM solutions within the interval $x \in (0,1)$ if Dirichlet boundary conditions are specified or $x \in [0,1]$ if Neumann boundary conditions are given. The solutions are implemented in \MATLAB using built-in functions, such as the vectorized \texttt{integral()}. To reduce computation time, we analytically evaluate sums, like those defining the forward discrete Fourier transform, and integrals when possible. In addition, all IBVPs have initial and boundary conditions compatible at $(x,t) = (0,0)$ and $(x,t) = (L,0)$. 


\section{Advection Equations}\label{advec_eqns_sec}

We begin by discussing advection equations in some detail, to show how the UTM is applied to semidiscrete problems on the finite interval. This allows us to introduce notation and to illustrate the kinds of numerical experiments we use here and for the more complicated examples of the following sections. Time $T$ is used to denote a fixed final time of interest in the time transforms \eqref{F} -- \eqref{G}, while $t < T$ refers to the generic time variable.


\subsection{Forward Discretization of $\bms{q_t = c\, q_{x}}$}\label{advec_forward_finiteinterval}
\begin{figure}[tb]
		\begin{center}
		\def\svgwidth{2.75in}
			\includegraphics[width=0.45\linewidth]{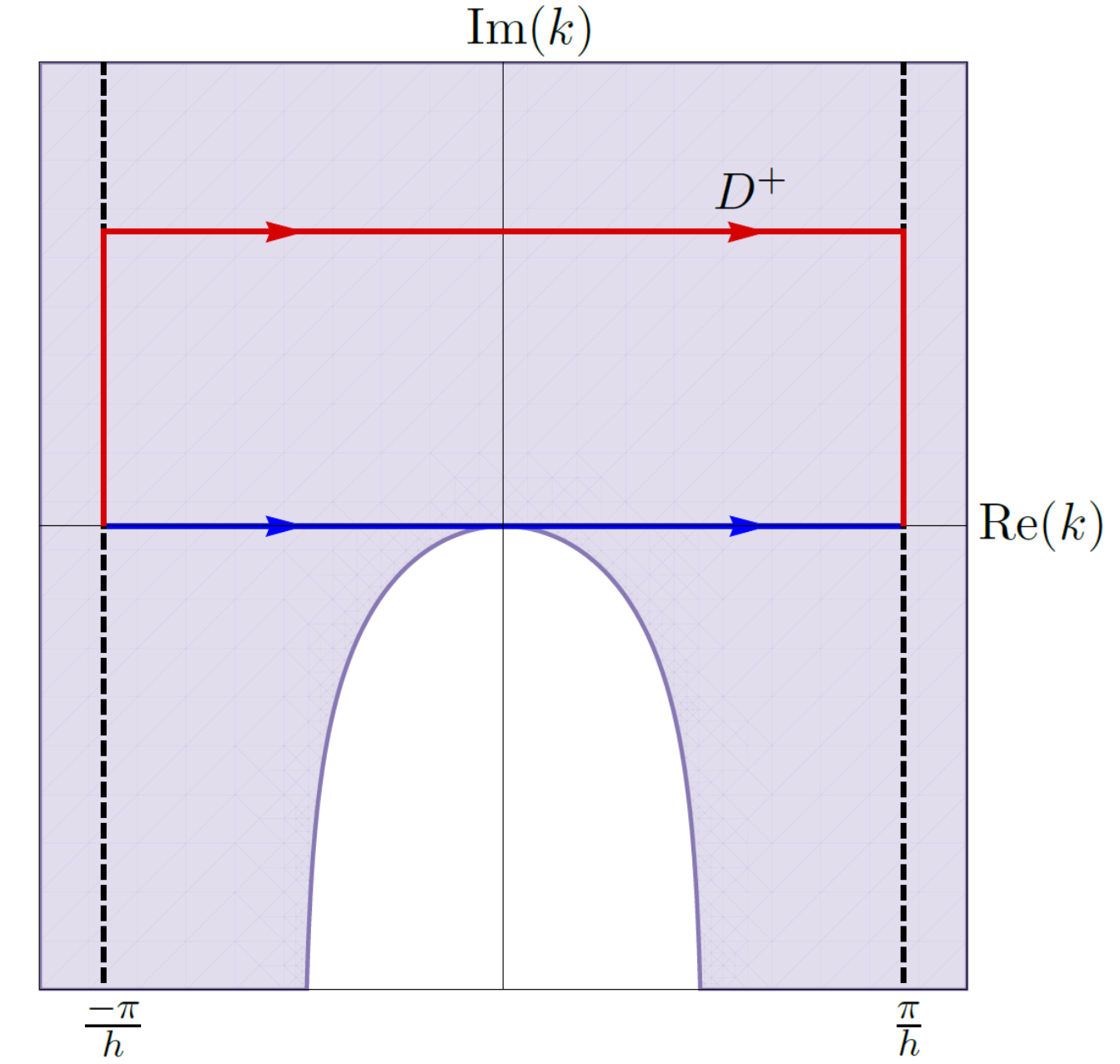}
			\caption{The shaded region depicts where $\text{Re}(-W) \leq 0$ and $e^{-WT}$ is bounded, for the dispersion relation \eqref{W_advec1_forward}.}
			\label{advec_forward_W}
		\end{center}
	\end{figure}

	We start with the continuous problem for the advection equation $q_t = c\, q_{x}$ with wave-speed $c > 0$:
		\begin{equation}\begin{dcases}
		q_t = c \, q_{x},& 0 < x < L,\, t > 0, \\
		q(x,0) = \phi(x),& 0 < x < L,\\
		q(L,t) = v^{(0)}(t),& t > 0.
		\label{advec1_prob}
	\end{dcases}\end{equation}
	For well posedness, the IBVP requires the initial condition and a Dirichlet boundary condition at $x = L$. Since information travels from right to left, a forward discretization of $q_x(x,t)$ is natural, and we consider
	\begin{align}
		\dot{q}_n(t) = c\,\frac{q_{n+1}(t) - q_{n}(t)}{h}.
		\label{advec1_forward}
	\end{align}
	As for the continuous UTM, the local relation is determined by writing the problem in divergence form. For this semidiscrete problem, we replace $\partial_x$ with the shift operator $\Delta$, and \eqref{advec1_forward} is rewritten as the one-parameter family of problems
	\begin{equation}
		\partial_t \left(e^{-iknh} e^{Wt} q_n \right) = \frac{c}{h}\Delta \left(e^{-ik(n-1)h} e^{Wt} q_{n} \right),
		\label{LR_advec1_forward}
	\end{equation}
	with dispersion relation
	\begin{equation}
		W(k) = c\, \frac{1 - e^{ikh}}{h}.
		\label{W_advec1_forward}
	\end{equation}
	The symmetries of the dispersion relation are those transformations $k \rightarrow \nu(k)$ that leave $W(k)$ invariant, \textit{i.e.}, $W(\nu) = W(k)$. Here, \eqref{W_advec1_forward} only has the trivial symmetry $\nu_0(k) = k$, up to periodic copies. From the local relation \eqref{LR_advec1_forward}, we obtain the global relation by taking a time transform over $t \in [0,T]$ and a finite sum from $n = 0$ (because $q_0(t)$ is not known) to $n = N$ (because $q_{N+1}(t)$ is known):
	\begin{align}
		&&\hspace{-55pt} \sum_{n=0}^{N} h \int_0^T \left[ \partial_t \left(e^{-iknh} e^{Wt} q_n \right) - \frac{c}{h}\Delta \left(e^{-ik(n-1)h} q_{n} \right)e^{Wt} \right] dt &= 0 \notag \\
		\Rightarrow \,&&\hspace{-55pt} e^{WT} \hat{q}(k,T) - \hat{q}(k,0) - c \left[ - e^{ikh} f_{0} + e^{-ik(L-h)} g_{0} \right] &= 0,
		\label{GR_advec1_forward}
	\end{align}
	valid for $k \in \mathbb{C}$. Solving for $\hat{q}(k,T)$ and inverting using the inverse transform \eqref{inv_fourier_cont_SD},
	\begin{align}\begin{split}
		q_n(T) &= \frac{1}{2\pi} \int_{-\pi/h}^{\pi/h} e^{iknh} e^{-WT} \hat{q}(k,0)\,dk + \frac{c}{2\pi} \int_{-\pi/h}^{\pi/h} e^{iknh} e^{-WT} \left[ - e^{ikh} f_{0} + e^{-ik(L-h)} g_{0} \right]\,dk.
		\label{soln1_advec_forward}
	\end{split}\end{align}
	 Since the Fourier transform $\hat{q}(k,0)$ consists of a finite sum, both integrands of \eqref{soln1_advec_forward} are defined for all $k \in \mathbb{C}$. We refer to the expression above as the ``solution,'' since $f_0(W,T)$ in the second integral term is not known, unlike $g_0(W,T)$. For $n = 0, \ldots, N$, $e^{ik(n+1)h}$ decays in the upper-half plane and $e^{-WT}$ is bounded in the shaded regions, including on the boundary, of Figure \ref{advec_forward_W}. This figure also shows the integration path for ``solution'' \eqref{soln1_advec_forward} from $-\pi/h$ to $\pi/h$ on the real line. Note that the sign of $c$ is essential in determining the region of exponential growth of the integrand, \textit{i.e.}, the white region in Figure \ref{advec_forward_W}.
	
	Following the approach discussed for half-line semidiscrete problems \cite{SDUTM_HL}, we show that \eqref{soln1_advec_forward} does not depend on $f_0(W,T)$, \textit{i.e.}, no Dirichlet boundary data at $x = 0$ is required. First, we substitute the definition of $f_0(W,T)$ in order to collect the $k$ dependence:
	\begin{align*}
		\frac{c}{2\pi} \int_{-\pi/h}^{\pi/h} e^{ik(n+1)h} e^{-WT} f_0\,dk &= \frac{c}{2\pi} \int_{-\pi/h}^{\pi/h} e^{ik(n+1)h} e^{-WT} \left[\int_{0}^T e^{Wt} q_0(t)\,dt \right]\,dk = \int_{0}^T A(n,T-t) q_0(t)\,dt,
	\end{align*} 
	with $T - t > 0$, and
	$$A(n,T) = \frac{c}{2\pi} \int_{-\pi/h}^{\pi/h} e^{ik(n+1)h} e^{-WT} \,dk.$$
	Let $R >0$. We define the line segment
	$$D^{+} = \left\{k \in \mathbb{C} \, \Big| \, \frac{-\pi}{h} \leq \text{Re}(k) \leq \frac{\pi}{h} \, \text{ and } \, \text{Im}(k) = R \right\},$$
	with left-to-right orientation, a horizontal straight-line path above the real line, from $k = -\pi/h + i R$ to $k = \pi/h + i R$. Next, we introduce a closed contour that consists of four straight segments: the original real-line path, the new path $D^{+}$, and two vertical segments that connect the endpoints of the real-line path with those of $D^{+}$ (see Figure \ref{advec_forward_W}). The contributions to the integral from these vertical paths cancel due to periodicity. Hence,
	$$A(n,T) = \frac{c}{2\pi} \int_{D^+} e^{ik(n+1)h} e^{-W(T-t)} \,dk,$$
	by Cauchy's Theorem. Taking $R \rightarrow +\infty$ implies $\text{Im}(k) \rightarrow \infty$ in the integrand. Because of the exponential decay above the real line, $A(n,T) = 0$ and
	$$\frac{c}{2\pi} \int_{-\pi/h}^{\pi/h} e^{ik(n+1)h} e^{-WT} f_0\,dk = 0.$$
	It follows that the solution to the finite interval IBVP with the forward discretization \eqref{advec1_forward} depends only on the initial condition and Dirichlet data at $x = L$:
	\begin{align}
		q_n(T) &= \frac{1}{2\pi} \int_{-\pi/h}^{\pi/h} e^{iknh} e^{-WT} \hat{q}(k,0)\,dk + \frac{c}{2\pi} \int_{-\pi/h}^{\pi/h} e^{ik(nh - L+h)} e^{-WT} g_{0}\,dk.
		\label{soln_advec_forward}
	\end{align}
	Note that $e^{-WT}$ grows in the lower-half plane, and we cannot deform and remove the integral contribution from $g_0(W,T)$ for the same $n$. Thus, \eqref{soln_advec_forward} is the final representation of the solution.
	
	For reference, we solve the IBVP \eqref{advec1_prob} using the continuous UTM, following Steps \eqref{utm_1} -- \eqref{utm_8} from Section \ref{cont_UTM_steps}. In short, we find the dispersion relation $\tilde{W}(k) = - i c k$, with only the trivial symmetry $\tilde{\nu}_0(k) = k$, and the global relation 
	$$\hat{q}(k,0) - e^{\tilde{W}T} \hat{q}(k,T) + c \left( - F_0 + e^{-i k L} G_0 \right)= 0,\quad k \in \mathbb{C},$$
	where
	$$\hat{q}(k,t) = \int_{0}^{L} e^{-ikx} q(x,t) \, dx ,\quad k \in \mathbb{C}$$
	and the time transforms are
	\begin{align*}
		F_j(\tilde{W},T) &= \int_{0}^T e^{\tilde{W}t} \left. \frac{\partial^j q}{d x^j} \right|_{x = 0} \, dt, \quad k \in \mathbb{C}, \\
		G_j(\tilde{W},T) &= \int_{0}^T e^{\tilde{W}t} \left. \frac{\partial^j q}{d x^j} \right|_{x = L} \, dt, \quad k \in \mathbb{C}.
	\end{align*}
	After inverting the Fourier transform and showing there is no dependence on $F_0(\tilde{W},T)$, the solution is represented as 
	\begin{align}
		q(x,T) &= \frac{1}{2 \pi} \int_{-\infty}^{\infty} e^{ikx} e^{-\tilde{W}T} \hat{q}(k,0) \, dk + \frac{c}{2\pi} \int_{-\infty}^{\infty} e^{ik(x - L)} e^{-\tilde{W}T} G_{0}\,dk.
		\label{soln_advec1_cont}
	\end{align}
	Taking the limit as $h \rightarrow 0$ of \eqref{soln_advec_forward}, we recover \eqref{soln_advec1_cont}, where the limits of integration approach $\pm \infty$ at rate $1/h$, $\lim_{h \rightarrow 0}W(k) = - c i k = \tilde{W}$, and $\lim_{h \rightarrow 0} e^{ik(nh - L+h)} = e^{ik(x-L)}$ with $n h = x_n \rightarrow x$ and
	$$\lim_{h \rightarrow 0} g_j (W,T) =\lim_{h \rightarrow 0} \int_{0}^T e^{Wt} q(L+jh,t) \, dt = \int_{0}^T e^{\tilde{W}t} q(L,t) \, dt = G_0(\tilde{W},T),$$
	for any fixed $j$.
	
	After substituting the definitions of $\hat{q}(k,0)$ and $G_0$ in \eqref{soln_advec1_cont}, we recover the classical, traveling wave solution:
	\begin{equation}
		q(x,T) = \begin{dcases}
			\phi(x+cT)\, , \quad &0 < x < L - cT, \\
			v^{(0)}\left(\tfrac{x-L}{c} + T\right)\,, \quad &L - cT < x < L.
		\end{dcases}
		\label{soln_advec1_cont_classical}
	\end{equation}


	\subsubsection{\textbf{Series Representation}}\label{advec_forward_finiteinterval_series} To facilitate numerical computation, we simplify the solution \eqref{soln_advec_forward} by substituting the definitions of $\hat{q}(k,0)$ and $g_0(W,T)$. For the initial-condition integral term with $\phi(x_n) \equiv \phi_n$,
	\begin{align*}
		\frac{1}{2 \pi} \int_{-\pi/h}^{\pi/h} e^{iknh} e^{-WT}\hat{q}(k,0)\,dk &= \sum_{m = 0}^N \left[ \frac{h}{2 \pi} \int_{-\pi/h}^{\pi/h} e^{ik(n-m)h} e^{-WT}\,dk \right] \phi_m. 
 	\end{align*} 
 	The integral can be evaluated with $z=e^{ikh}$ and $W(z) = c(1-z)/h$:
 	\begin{align*}
 		\hspace{-20pt}\frac{h}{2 \pi} \int_{-\pi/h}^{\pi/h} e^{ik(n-m)h} e^{-WT}\,dk &= \frac{e^{-cT/h}}{2 \pi i} \oint_{|z|=1} z^{n-m-1} e^{\tfrac{cT}{h}z}\, dz = e^{-cT/h} \,\underset{z = 0}{\text{Res}} \left\{ z^{n-m-1} e^{\tfrac{cT}{h}z}\right\} = \frac{e^{-cT/h}}{\left(m - n\right) !} \left(\frac{cT}{h}\right)^{m - n}.
 	\end{align*}
 	Since $m \geq n$, the initial condition integral gives
 	\begin{align*}
 		\hspace{-25pt}\frac{1}{2 \pi} \int_{-\pi/h}^{\pi/h} e^{iknh} e^{-WT}\hat{q}(k,0)\,dk = \sum_{m = n}^N \frac{e^{-cT/h}}{\left(m - n\right) !} \left(\frac{cT}{h}\right)^{m - n} \phi_m = e^{-cT/h} \sum_{m = 0}^{N-n} \left(\frac{cT}{h}\right)^{m} \frac{\phi_{n+m}}{m !}.
 	\end{align*}
 	The boundary integral from \eqref{soln_advec_forward} is treated similarly by substituting the definition for $g_0(W,T)$ in terms of $v^{(0)}(t)$:
	\begin{align*}
		\frac{c}{2 \pi} \int_{-\pi/h}^{\pi/h} e^{ik(nh-L+h)} e^{-WT} g_0\,dk &= c\int_0^T \left[ \frac{1}{2 \pi} \int_{-\pi/h}^{\pi/h} e^{ik(n -N)h} e^{-W(T-t)} \,dk \right]  v^{(0)}(t) \, dt,
	\end{align*}
	after substituting $L = (N+1)h$.	Again with $z = e^{ikh}$, 
	\begin{align*}
		\frac{1}{2 \pi} \int_{-\pi/h}^{\pi/h} e^{ik(n -N)h} e^{-W(T-t)} \,dk &= \frac{e^{-c(T-t)/h}}{2 \pi i h} \oint_{|z|=1} z^{n -N - 1} e^{\tfrac{c(T-t)}{h}z} \,dz \\
		&= \frac{e^{-c(T-t)/h}}{h} \,\underset{z = 0}{\text{Res}} \left\{ z^{n -N - 1} e^{\tfrac{c(T-t)}{h}z} \right\}\\
		&= \frac{e^{-c(T-t)/h}}{h} \left(\frac{c(T-t)}{h}\right)^{N - n} \frac{1}{(N - n) !}.
	\end{align*}
	Hence,
	$$\frac{c}{2 \pi} \int_{-\pi/h}^{\pi/h} e^{ik(nh-L+h)} e^{-WT} g_0\,dk  = \frac{c}{h} \int_0^T e^{-c(T-t)/h} \left(\frac{c(T-t)}{h}\right)^{N - n} \frac{ v^{(0)}(t)}{ (N - n) !}  \, dt. $$
	Thus, \eqref{soln_advec_forward} is rewritten as
	\begin{align}
		q_n(T) &= e^{-cT/h} \sum_{m = 0}^{N-n} \left(\frac{cT}{h}\right)^{m} \frac{\phi_{n+m}}{m !}  + \frac{c}{h} \int_0^T e^{-c(T-t)/h} \left(\frac{c(T-t)}{h}\right)^{N - n} \frac{ v^{(0)}(t)}{ (N - n) !}  \, dt. 
		\label{soln1_advec_forward_num}
	\end{align}

	Unlike for the classical continuous solution \eqref{soln_advec1_cont_classical}, the initial and boundary conditions in the semidiscrete solution \eqref{soln1_advec_forward_num} are inseparable and both contribute at every mesh point. For computational purposes, the factorials become too large to store in finite precision as we refine $h$, and one last rewrite is in order. For the sum in \eqref{soln1_advec_forward_num}, 
	\begin{align*}
		\left(\frac{cT}{h}\right)^{m} \frac{1}{m !} & = \exp \left[m \ln \left(\frac{cT}{h} \right) - \ln\left[ \Gamma \left(m+1\right) \right] \right],
	\end{align*}  
	so that combining with $e^{-cT/h}$ gives
	$$e^{-cT/h} \sum_{m = 0}^{N-n} \left(\frac{cT}{h}\right)^{m} \frac{\phi_{n+m}}{m !} = \sum_{m = 0}^{N-n} \exp \left[m \ln \left(\frac{cT}{h} \right) - \ln\left[ \Gamma \left(m+1\right) \right] - \frac{cT}{h} \right] \phi_{n+m}.$$
	After a similar rewrite for the integral term, \eqref{soln1_advec_forward_num} becomes
	\begin{align}
	\begin{split}
		q_n(T) &=  \sum_{m = 0}^{N-n} \exp \left[m \ln \left(\frac{cT}{h} \right) - \ln\left[ \Gamma \left(m+1\right) \right] - \frac{cT}{h} \right] \phi_{n+m} \\
		&\quad\, + c\int_0^T \exp \left[(N - n) \ln \left(\frac{c(T-t)}{h} \right) - \ln\left[ \Gamma \left(N - n+1\right) \right] - \frac{c(T-t)}{h} - \ln(h) \right]  v^{(0)}(t)  \, dt.
		\label{soln_advec_forward_num}
	\end{split}
	\end{align}
	Using \MATLAB, we can make use of the built-in \texttt{integral()} and \texttt{gammaln()} functions.
	
	From the stencil \eqref{advec1_forward}, we know that \eqref{soln_advec_forward}, and hence \eqref{soln_advec_forward_num}, is a first-order accurate approximation to the solution $q(x,T)$ of the IBVP \eqref{advec1_prob}. We can reveal more information about the behavior and structure of this approximate solution by determining its modified equation \cite{randy,mod_pde}. Suppose $q_n(T)$ samples the solution of a PDE with dependent variable $p(x,T)$, such that $q_n(T) \equiv p\left(x_n,T\right)$. Substituting into the forward stencil \eqref{advec1_forward} and Taylor expanding gives
	\begin{align*}
		\dot{q}_n(t) &= \frac{c}{h} \left[ q_{n+1}(t) - q_{n}(t) \right] \\
		\Rightarrow \quad p_t(x_n,t) &= \frac{c}{h} \left[ p(x_n+h,t) - p(x_n,t) \right] \\
		\Rightarrow \quad p_t &= c\, p_x + \frac{c\,p_{xx}}{2} h + \frac{c\, p_{xxx}}{6} h^2 + \mathcal{O}\left(h^3\right).
	\end{align*}
	Keeping the $\mathcal{O}(h)$ term, we find that \eqref{soln_advec_forward} is a second-order accurate solution approximation to the advection-diffusion PDE:
	\begin{equation}
		p_t = c \,p_x + \frac{c\,h}{2} p_{xx},
		\label{advec_forward_modified_eqn}
	\end{equation}
	so we expect solution profiles of \eqref{soln_advec_forward} to travel at the correct speed $c$, while dissipating in time. Since $c >0$, the diffusion coefficient $c\,h/2$ is positive. If we allow $c < 0$ or if we apply the same forward stencil to the PDE $q_t = - a \,q_x$ with $a >0$, we obtain a similar convection-diffusion modified PDE, but with a negative diffusion coefficient, which gives an ill-posed problem with exponentially growing solutions.
	
\begin{figure}[tb]	
		\raggedleft
		\begin{subfigure}[t]{.45\textwidth}
			\centering
  			\includegraphics[width=1.1\linewidth]{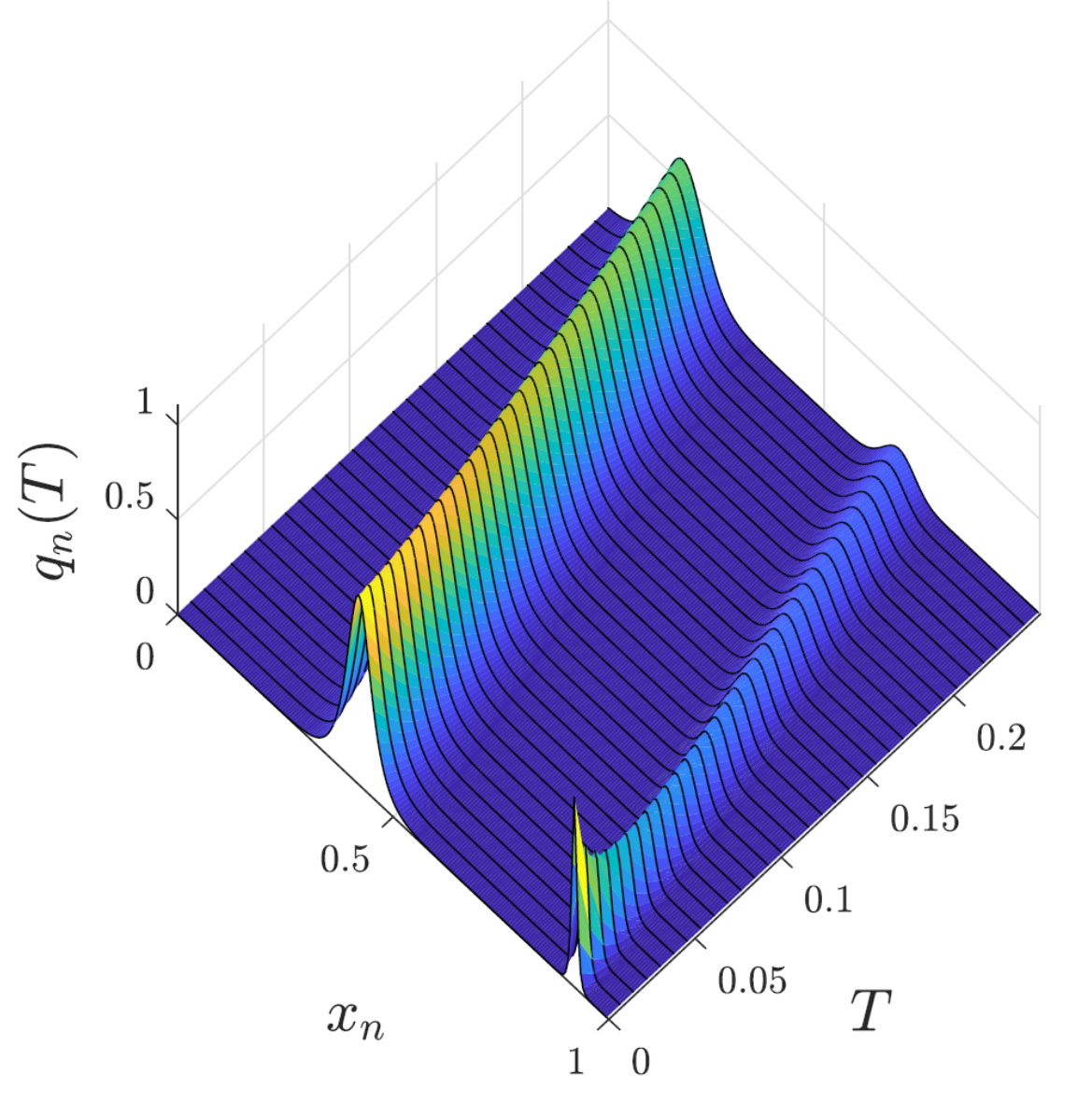}
  			\caption{}
  			\label{advec1_UTM1_FI_solnplot}
		\end{subfigure}\hfill 
		\begin{subfigure}[t]{.45\textwidth}
			\centering
  			\includegraphics[width=1\linewidth]{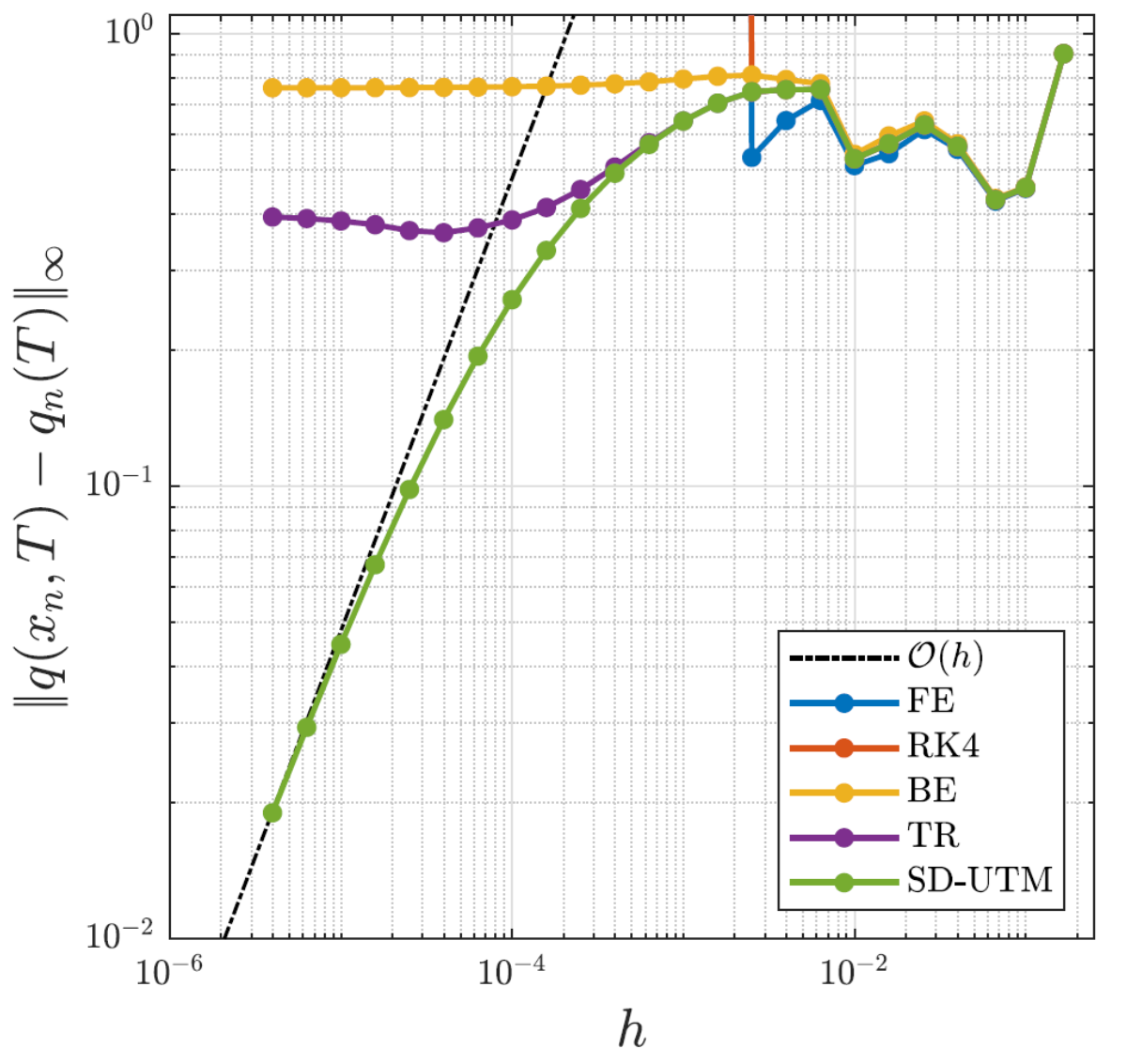}
  			\caption{}
  			\label{advec1_UTM1_FI_errorplot}
		\end{subfigure}	
		\caption{(a) The semidiscrete solution \eqref{soln_advec_forward_num} evaluated at various $T$ with $h = 0.005$. (b) Error plot of the semidiscrete solution \eqref{soln_advec_forward_num} and finite-difference schemes relative to the exact solution as $h \rightarrow 0$ with $T = 0.25$ and $\Delta t = 2.5 \times 10^{-3}$.}
		\label{advec1_UTM1_FI}
	\end{figure}

	
	As an explicit example, we compute the numerical solution of 
	\begin{equation}
	\begin{dcases}
		q_t = q_{x},& 0 < x < 1,\, t > 0, \\
		q(x,0) = \phi(x) = \sech\left[200 (x-0.925)\right] + \sech\left[40 (x-0.425)\right],& 0 < x < 1,\\
		q(0,t) = v^{(0)}(t) = \phi(1 + t), &t > 0, 
	\end{dcases}
	\label{advec1_numerical1_FI}
	\end{equation}
	with $c = 1$ and $L = 1$, so that the boundary condition acts as the continuation of the initial condition from outside the interval for $t > 0$. The initial condition consists of two peaks of equal heights, except the leading peak is wider than the trailing peak. We solve \eqref{advec1_numerical1_FI} with the aforementioned finite-difference numerical methods and the SD-UTM solution \eqref{soln_advec_forward_num} via a first-order spatial forward discretization. Figure \ref{advec1_UTM1_FI} shows the semidiscrete solution $q_n(t)$ (left panel) and a log-log error plot (right panel) of the $\infty$-norm of $q_n(0.25)-q(x_n, 0.25)$, as a function of $h$, where the finite-difference schemes use a fixed time step $\Delta t \equiv \Delta t = 2.5 \times 10^{-3}$.
	 The $(x_n,t)$-plot in Figure \ref{advec1_UTM1_FI_solnplot} shows that both peaks decrease in amplitude and widen as time increases, predicted by the modified PDE \eqref{advec_forward_modified_eqn}, with the narrow peak quickly dissipating compared to the wider peak. The error plot in Figure \ref{advec1_UTM1_FI_errorplot} shows that SD-UTM works well compared to the traditional numerical methods. The explicit methods become unstable after their CFL conditions are violated \cite{randy} and the implicit methods' errors are asymptotic to the temporal truncation errors as $h \rightarrow 0$ for a fixed time step $\Delta t$. Figure \ref{advec1_UTM1_FI_errorplot} demonstrates there is no CFL condition for the SD-UTM to succeed, while, for example, FE only works well for large $h$ values where $\Delta t/h \leq 1$. The implicit methods do not have such restrictions, but have truncation error $\mathcal{O}\left(h\right) + \mathcal{O}\left(\Delta t^p \right)$, where $p = 1$ for BE and $p = 2$ for TR. For a fixed $\Delta t$, only the spatial error decreases as $ h \rightarrow 0$, while $\mathcal{O}\left(\Delta t^p \right)$ remains and eventually dominates. The asymptotic limits of BE and TR in Figure \ref{advec1_UTM1_FI_errorplot} as $h \rightarrow 0$ reveal this temporal truncation error. Figure \ref{advec1_UTM1_FI_errorplot} implies that the SD-UTM solution \eqref{soln_advec_forward_num} has a slow convergence rate to the continuous solution for this IBVP, likely due to the sharp peak, reaching the expected $\mathcal{O}\left(h\right)$ for $h < 10^{-5}$. 

\begin{figure}[tb]	
		\raggedleft
		\begin{subfigure}[t]{.45\textwidth}
			\centering
  			\includegraphics[width=1\linewidth]{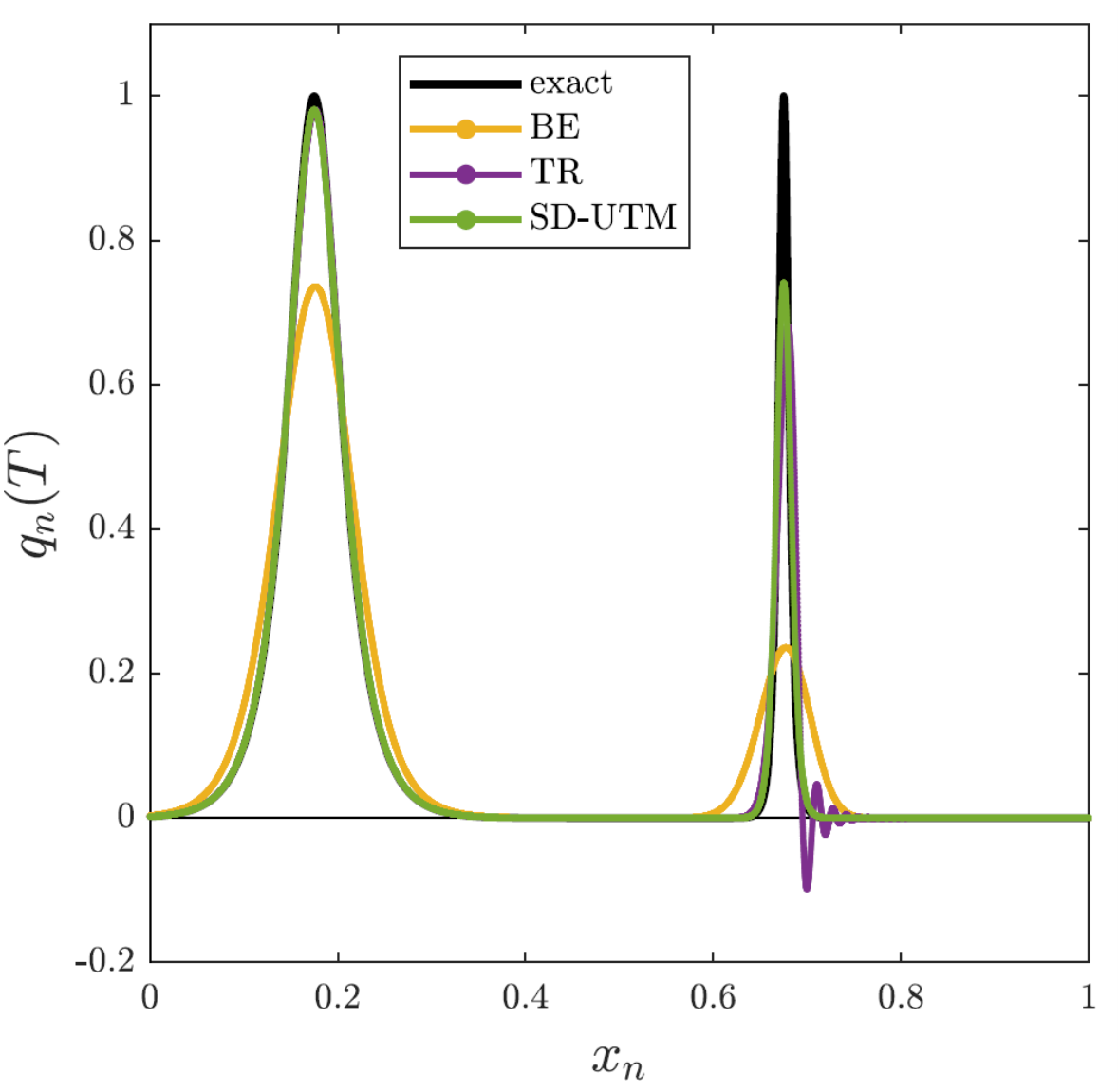}
  			\caption{}
  			\label{advec1_UTM1_FI_solnplot}
		\end{subfigure}\hfill 
		\begin{subfigure}[t]{.45\textwidth}
			\centering
  			\includegraphics[width=1\linewidth]{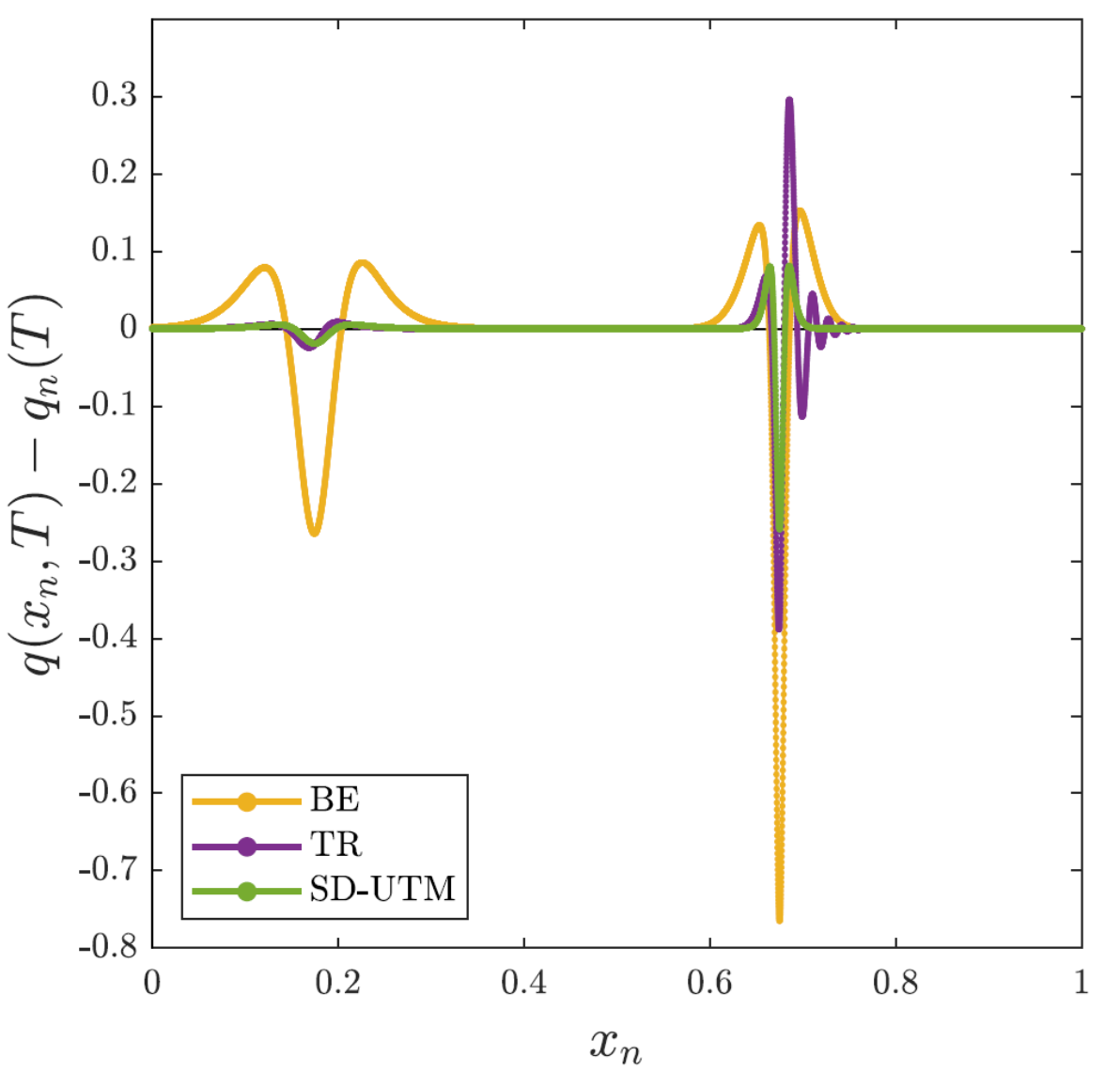}
  			\caption{}
  			\label{advec1_UTM1_FI_errorplot}
		\end{subfigure}	
		\caption{(a) The numerical solutions to IBVP \eqref{advec1_numerical1_FI} at $T = 0.25$ with $h = 10^{-4}$ for all the methods and $\Delta t = 2.5 \times 10^{-3}$ for the finite-difference methods. (b) The difference between the exact solution and the numerical solutions under the same conditions as (a).}
		\label{advec1_numerical1_t}
	\end{figure}

	 Figure \ref{advec1_numerical1_t} compares the exact solution with the numerical solution profiles for all methods, except FE and RK4 due to their instabilities, with $h = 10^{-4}$. For the SD-UTM, we simply compute the solution \eqref{soln_advec_forward_num} at $T = 0.25$ once, while the finite-difference solutions time-step to $T = 0.25$ with step size $\Delta t = 2.5 \times 10^{-3}$.
	From the solution plot, it appears that every method is dissipative, including FE and RK4 (not shown), and TR is also dispersive. With the SD-UTM, both peaks drop in amplitude and diffuse, while the dispersive nature of TR is apparent near the sharper peak. Despite the lack of a dispersive tail like the other implicit method, it appears that BE is more dissipative than all the other tested methods. In summary, the SD-UTM performs better than the finite-difference methods presented. Even though dissipation is evident, dispersion is not.


\subsection{Higher-Order One-Sided Discretization of $\bms{q_t = c\, q_{x}}$}\label{advec_forward_finiteinterval_2}
	All forward discretizations produce $f_j(W,T)$ terms with a coefficient $C_j\,e^{i \gamma_j kh}$ for some $C_j \in \mathbb{C}$ and $\gamma_j \in \mathbb{N}$ in the global relation. Coupled with a polynomial dispersion relation $W(z)$, we can remove all integral terms containing any $f_j(W,T)$ from ``solutions'' using the steps above. On the opposite side of the interval, more terms with $g_j(W,T)$ are introduced, but with the help of symmetries or additional boundary conditions given by the PDE, we can remove the unknown $g_j(W,T)$ terms. The steps in the semidiscrete UTM become more intricate and tedious, yet remain systematic. 
	
\begin{figure}[tb]
	\raggedleft
	\begin{subfigure}[t]{.45\textwidth}
		\begin{center}
		\def\svgwidth{2.75in}
			\includegraphics[width=1\linewidth]{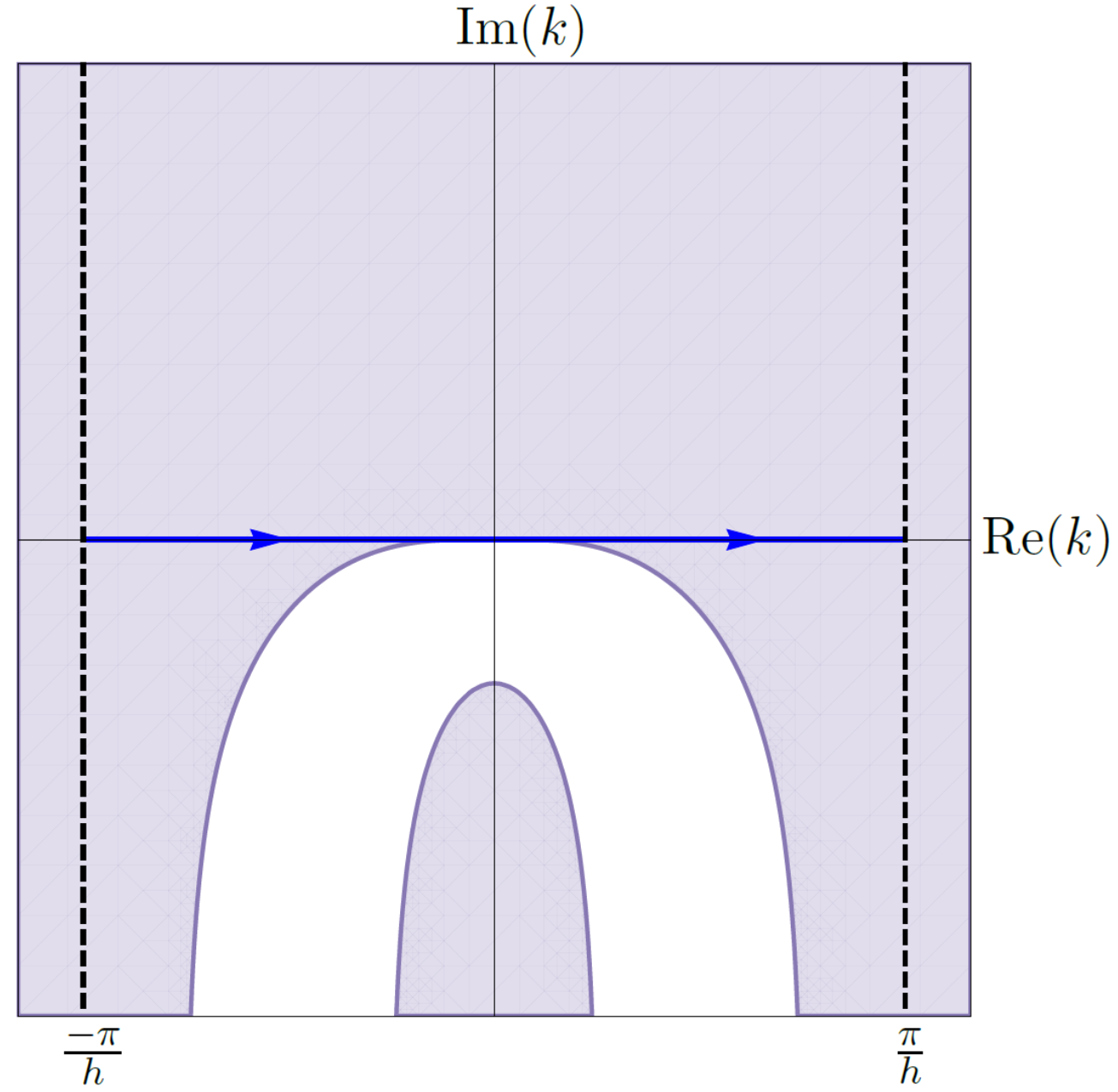}
			\caption{}
			\label{W_advec1_backward2}
		\end{center}
	\end{subfigure}\hfill 
	\begin{subfigure}[t]{.45\textwidth}
		\begin{center}
		\def\svgwidth{2.75in}
			\includegraphics[width=1.\linewidth]{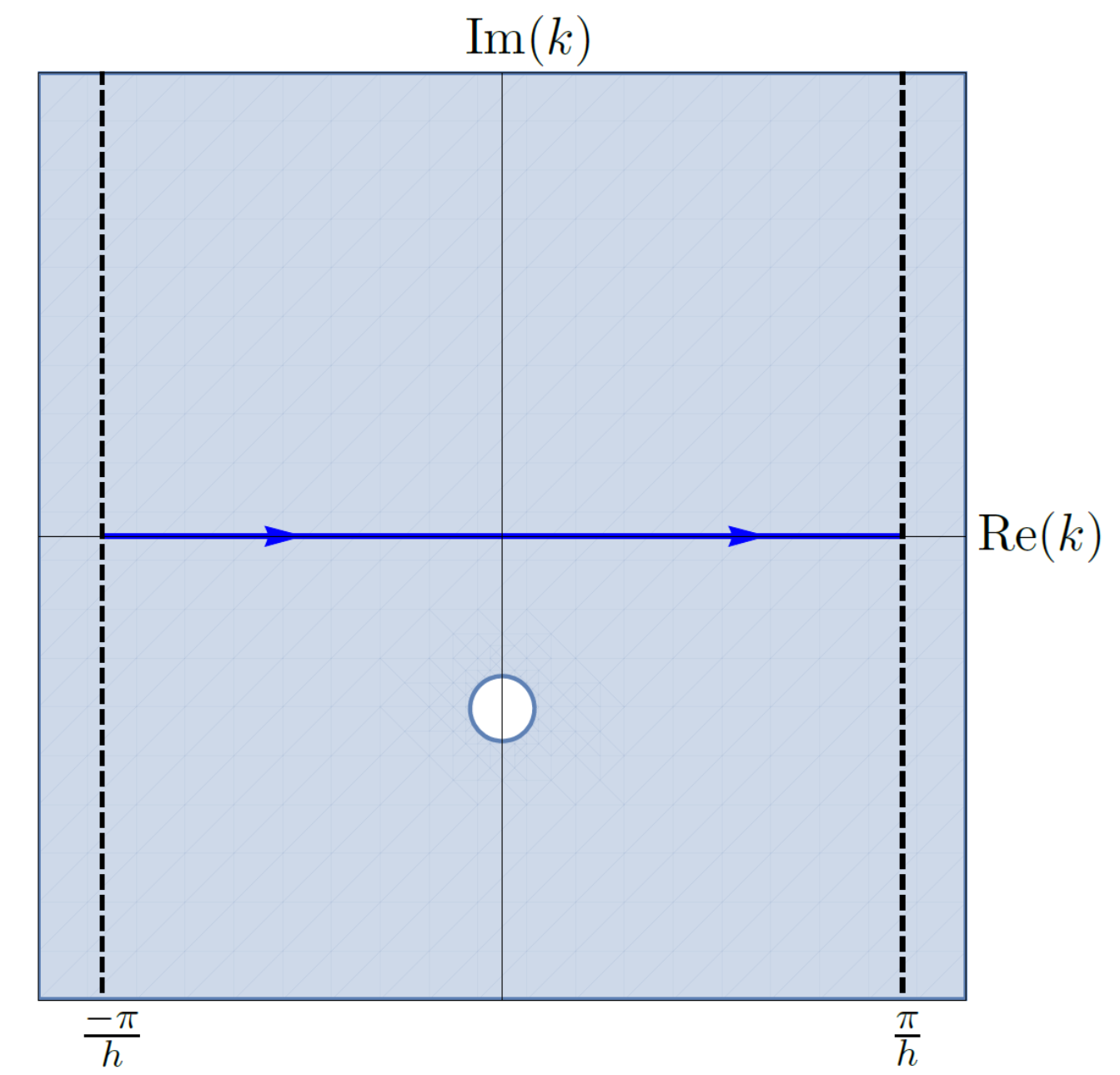}
			\caption{}
			\label{symm_advec1_backward2}
		\end{center}
	\end{subfigure}
	\caption{(a) The shaded regions depict where $\text{Re}(-W) \leq 0$ and $e^{-WT}$ is bounded, for the dispersion relation \eqref{W_advec1_forward2}. (b) The shaded regions depict where the global relation with $ k \rightarrow \nu_1(k)$ is valid, \textit{i.e.}, $\text{Im}(\nu_1)\leq 0$.} 
	\label{advec1_backward2_fig}
	\end{figure}

	We consider a second-order stencil for a forward discretization of $q_x(x,t)$:
	\begin{align}
		\dot{q}_n(t) = c\,\frac{-3 q_{n}(t) + 4 q_{n+1}(t) - q_{n+2}(t)}{2h}.
		\label{advec1_forward2}
	\end{align}
	We find the local relation
	\begin{equation}
		\partial_t \left(e^{-iknh} e^{Wt} q_n \right) = \frac{c}{2h}\Delta \left(4 e^{-ik(n-1)h} e^{Wt} q_{n} - e^{-ik(n-1)h} e^{Wt} q_{n+1} - e^{-ik(n-2)h} e^{Wt} q_{n}\right),
		\label{LR_advec1_forward2}
	\end{equation}
	with dispersion relation 
	\begin{equation}
		W(k) = c\, \frac{3 - 4 e^{ikh} + e^{2ikh}}{2h}.
		\label{W_advec1_forward2}
	\end{equation}
	Taking a time transform and a finite sum from $n = 0$ to $n = N$, the global relation is
	\begin{align}
		e^{WT} \hat{q}(k,T) - \hat{q}(k,0) - \frac{c}{2} \left( f(k,T) + e^{-ikL} g(k,T) \right) &= 0, \quad k \in \mathbb{C},
		\label{GR_advec1_forward2}
	\end{align}
	where 
	\[\begin{dcases*}
		f(k,T) = e^{2 i k h}f_0 - 4 e^{i k h} f_0 + e^{i k h} f_1,\\
		g(k,T) = 4 e^{i k h} g_0 - e^{2 i k h} g_0 - e^{i k h}g_1. 
	\end{dcases*}\]	
	Solving for $\hat{q}(k,T)$ and using the inverse transform,
	\begin{align}\begin{split}
		q_n(T) &= \frac{1}{2\pi} \int_{-\pi/h}^{\pi/h} e^{iknh} e^{-WT} \hat{q}(k,0)\,dk + \frac{c}{2\pi} \int_{-\pi/h}^{\pi/h} e^{iknh} e^{-WT}  \left[ \frac{f(k,T) + e^{-ikL} g(k,T)}{2} \right]\,dk.
		\label{soln1_advec_forward2}
	\end{split}\end{align}
	We can deform and remove $f(k,T)$ from ``solution'' \eqref{soln1_advec_forward2}, so that
	\begin{align}\begin{split}
		q_n(T) &= \frac{1}{2\pi} \int_{-\pi/h}^{\pi/h} e^{iknh} e^{-WT} \hat{q}(k,0)\,dk + \frac{c}{2\pi} \int_{-\pi/h}^{\pi/h} e^{ik(nh-L+h)} e^{-WT} \left[ \frac{\left(4 - e^{ i k h} \right) g_0 - g_1}{2} \right]\,dk.
		\label{soln2_advec_forward2}
	\end{split}\end{align}
	Figure \ref{W_advec1_backward2} illustrates this, since $e^{-WT}$ is bounded in the whole upper-half plane. Thus, the second integral term of \eqref{soln2_advec_forward2} depends only on the transformed Dirichlet data and data at the unknown ghost point $q_{N+2}(T) = q(L+h,T)$ through $g_1(W,T)$.  
	
	The dispersion relation \eqref{W_advec1_forward2} has the nontrivial symmetry 
	$$\nu_1(k) = \frac{\ln \left(4-e^{i k h}\right)}{i h},$$
	up to periodic copies. The global relation \eqref{GR_advec1_forward2} with $k \rightarrow \nu_1(k)$ is valid for all $k \in \mathbb{C}$, except for a bounded region in the lower-half plane shown in Figure \ref{symm_advec1_backward2}. We solve this global relation for the unknown $g_1(W,T)$ to find
	$$g_1 = e^{i k h} g_0  -  e^{i k h} \left(4-e^{i k h}\right)^{N+1} f_0 + \left(4-e^{i k h}\right)^{N+1} f_1 + \frac{2 \left(4-e^{i k h}\right)^N }{c} \left[ \hat{q}(\nu_1,0) - e^{W T}  \hat{q}(\nu_1,T) \right],$$
	where $e^{i \nu_1 L} = \big(4 - e^{ikh}\big)^{N+1}$ with $L = (N+1)h$.	Substituting into \eqref{soln2_advec_forward2} does not only reintroduce dependence on $f_0(W,T)$ and $f_1(W,T)$ that can no longer be deformed away, but it also introduces a nonzero contribution from $\hat{q}(\nu_1,T)$, the transform of the solution itself.
	
	In \cite{SDUTM_HL}, we present an alternative route to obtain a valid solution representation. Returning to the continuous problem \eqref{advec1_prob}, the PDE itself gives first and second-derivative boundary conditions from the Dirichlet condition:
	\begin{subequations}
	\begin{equation}
		q_x(L,t) = \frac{\dot{v}(t)}{c}, \quad \dot{v}(t) = \frac{d}{dt}v^{(0)}(t),
		\label{advec1_neumann_bc}
	\end{equation}
	\begin{equation}
		q_{xx}(L,t) = \frac{\ddot{v}(t)}{c^2}, \quad \ddot{v}(t) = \frac{d^2}{dt^2}v^{(0)}(t).
		\label{advec1_neumann_bc2}
	\end{equation}
	\end{subequations}
	We discretize the derivative conditions using centered second-order accurate stencils and apply time transforms:
	\begin{subequations}
	\begin{equation}
		\frac{g_{1} - g_{-1}}{2h} = \frac{\dot{V}}{c}, \quad\quad \dot{V}(W,T) = \int_0^T e^{Wt} \dot{v}(t)\,dt,
		\label{advec1_neumann_bc_g1}
	\end{equation}
	\begin{equation}
		\frac{g_{1} - 2 g_0 + g_{-1}}{h^2} = \frac{\ddot{V}}{c^2}, \quad\quad \ddot{V}(W,T) = \int_0^T e^{Wt} \ddot{v}(t)\,dt.
		\label{advec1_neumann_bc_g2}
	\end{equation}
	\end{subequations}
	Solving \eqref{advec1_neumann_bc_g1} and \eqref{advec1_neumann_bc_g2} for $g_1(W,T)$ and $g_{-1}(W,T)$ gives
	\begin{align}\begin{split}
		q_n(T) &= \frac{1}{2\pi} \int_{-\pi/h}^{\pi/h} e^{iknh} e^{-WT} \hat{q}(k,0)\,dk + \frac{c}{2\pi} \int_{-\pi/h}^{\pi/h} e^{ik(nh-L+h)} e^{-WT} \left[ \frac{\left(3 - e^{ i k h} \right) }{2}g_0 \right]\,dk \\
		&\quad\, - \frac{c}{2\pi} \int_{-\pi/h}^{\pi/h} e^{ik(nh-L+h)} e^{-WT} \left[ \frac{h  }{2c} \dot{V} + \frac{h^2 }{4c^2}\ddot{V}  \right]\,dk.
		\label{soln_advec_forward2}
	\end{split}\end{align}
	The additional steps of including \eqref{advec1_neumann_bc_g1} and \eqref{advec1_neumann_bc_g2} allow \eqref{soln_advec_forward2} to maintain $\mathcal{O}(h^2)$ accuracy. The modified PDE corresponding to the second-order discretization \eqref{advec1_forward2} is the dispersive PDE $p_t = c p_x - (c h^2/3) p_{xxx}$, and \eqref{soln_advec_forward2} is its third-order approximation (the omitted higher-order term in the modified PDE is $\mathcal{O}(h^3)$). The stencils \eqref{advec1_neumann_bc_g1} and \eqref{advec1_neumann_bc_g2} are both fourth-order approximations to their respective modified PDEs with nonzero $\mathcal{O}(h^2)$ coefficients. As before, the semidiscrete solution \eqref{soln_advec_forward2} converges to \eqref{soln_advec1_cont} as $h \rightarrow 0$ and correctly loses dependence on the Neumann boundary condition in the continuum limit.


 	\subsubsection{\textbf{Series Representation}} We rewrite solution \eqref{soln_advec_forward2} by substituting the definitions of $\hat{q}(k,0)$, $g_0(W,T)$, and $\dot{V}_0(W,T)$. For the initial-condition integral term,
	\begin{align}
 		\frac{1}{2 \pi} \int_{-\pi/h}^{\pi/h} e^{iknh} e^{-WT}\hat{q}(k,0)\,dk =  e^{-3cT/(2h)} \sum_{m = 0}^{N-n} \sum_{k = 0}^{m/2} \frac{ 4^{m-2k} (-1)^k }{(m-2k)!\, k!} \left(\frac{cT}{2h} \right)^{m-2k} \phi_{m+n}.
 	\label{IC_soln_advec_forward2}
 	\end{align}
 	For the boundary integrals from \eqref{soln_advec_forward2},
	\begin{align*}
		\frac{c}{2\pi} \int_{-\pi/h}^{\pi/h} e^{ik(nh-L+h)} e^{-WT} \left[ \frac{\left(3 - e^{ i k h} \right) }{2}\right.&\left.g_0 - \frac{h  }{2c} \dot{V} - \frac{h^2 }{4c^2}\ddot{V} \right]\,dk \\
		&= 3 B_0\big(n,T,v^{(0)}\big) - B_1\big(n,T,v^{(0)}\big) - \frac{h}{c} B_0(n,T,\dot{v}) - \frac{h^2}{2c^2}B_0(n,T,\ddot{v}),
	\end{align*}
	with
	\begin{align}
		B_j(n,T,v) &= \frac{c}{4\pi} \int_{-\pi/h}^{\pi/h} e^{ik(nh-L+h + j h)} e^{-WT} V(W,T) \,dk  \notag\\
		&= \frac{a c}{2 h} \sum_{k = 0}^{(N-j - n)/2}  \frac{4^{N-j - n-2k} (-1)^k }{(N-j - n-2k)!\, k!} \int_0^T e^{-3c(T-t)/(2h)} \left(\frac{c(T-t)}{2h} \right)^{N-j - n-2k} v(t) \,dt ,
		\label{BC_soln_advec_forward2}
	\end{align}
	and 
	$$V(W,T) = \int_0^T e^{Wt} v(t)\, dt,$$
	after substituting definitions and expanding. Combining the initial-condition term \eqref{IC_soln_advec_forward2} and the boundary-condition term \eqref{BC_soln_advec_forward2} allows for a different representation of \eqref{soln_advec_forward2}, and we can use optimized built-in functions in \MATLAB and other languages. As for \eqref{soln1_advec_forward_num}, the initial and boundary conditions contribute at every interior mesh point.

\begin{remark}
	The choice of equations to remove additional unknowns is not necessarily unique. The main challenge is finding equations that are linearly independent and give the desired order of accuracy. For example, instead of the centered discretization \eqref{advec1_neumann_bc_g1} for $q_x(L,t)$, we apply the second-order forward stencil:
	\begin{align}
		\frac{-3q_{N+1}(t) + 4 q_{N+2}(t) - q_{N+3}(t)}{2h} = \frac{\dot{v}(t)}{c} \quad\quad \Rightarrow \quad\quad \frac{-3 g_{0} + 4 g_{1} - g_2}{2h} = \frac{\dot{V}}{c}.
		\label{advec1_neumann_bc_g3}
	\end{align}
	This choice requires a second equation, different from \eqref{advec1_neumann_bc_g2}, that does not introduce any new unknowns. We discretize $q_{xx}(L,t)$ using the \textit{first}-order forward stencil: 
	\begin{align}
		\frac{q_{N+1}(t) - 2 q_{N+2}(t) + q_{N+3}(t)}{h^2} = \frac{\ddot{v}(t)}{c^2} \quad\quad \Rightarrow \quad\quad \frac{g_0 - 2 g_1 + g_2}{h^2} = \frac{\ddot{V}}{c^2}.
		\label{advec1_neumann_bc_g4}
	\end{align}
	Interestingly, both pairs of discretizations, \eqref{advec1_neumann_bc_g1} -- \eqref{advec1_neumann_bc_g2} and \eqref{advec1_neumann_bc_g3} -- \eqref{advec1_neumann_bc_g4}, give the same expression for the unknown $g_1(W,T)$ and the same second-order accurate solution \eqref{soln_advec_forward2}. It is noteworthy that $g_1(W,T)$ in \eqref{advec1_neumann_bc_g1} and \eqref{advec1_neumann_bc_g3} arises as $g_1/h$. Similarly rewriting the $g_1(W,T)$ term in \eqref{advec1_neumann_bc_g2} and \eqref{advec1_neumann_bc_g4}, we have
	$$\frac{g_{1} - 2 g_0 + g_{-1}}{h} = h \frac{\ddot{V}}{c^2} + \mathcal{O}(h^3), \quad \frac{g_0 - 2 g_1 + g_2}{h} = h \frac{\ddot{V}}{c^2} + \mathcal{O}(h^2),$$
	respectively, so that $g_1(W,T)$ is solved (at least) to $\mathcal{O}(h^2)$ as in \eqref{advec1_neumann_bc_g1} and \eqref{advec1_neumann_bc_g3}. To remove $g_{1}(W,T)$ from \eqref{soln2_advec_forward2} without introducing new unknowns, we can discretize the Neumann condition \eqref{advec1_neumann_bc} using the standard forward stencil to find
	\begin{align}\begin{split}
		q_n(T) &= \frac{1}{2\pi} \int_{-\pi/h}^{\pi/h} e^{iknh} e^{-WT} \hat{q}(k,0)\,dk + \frac{c}{2\pi} \int_{-\pi/h}^{\pi/h} e^{ik(nh-L+h)} e^{-WT} \left[ \frac{\left(3 - e^{ i k h} \right) }{2}g_0 - \frac{h}{2c} \dot{V} \right]\,dk.
		\label{soln_advec_forward2_low}
	\end{split}\end{align}
	Because of the discretization of the Neumann condition, the accuracy of \eqref{soln_advec_forward2_low} is $\mathcal{O}(h)$ instead of the expected $\mathcal{O}(h^2)$. In fact, from this discretization, the modified PDE $p_x(L,t) = \dot{v}(t) + (h/2)p_{xx}(L,t)$ implies local dissipation near the $x = L$ boundary. Even so, solution \eqref{soln_advec_forward2_low} converges to \eqref{soln_advec1_cont} and loses dependence on the Neumann boundary condition as $h \rightarrow 0$. Note that an integral term with $h^2 \ddot{V}/(4c^2)$ is the only difference between \eqref{soln_advec_forward2} and \eqref{soln_advec_forward2_low}.
\end{remark}

	\begin{remark}
		Unlike for the half-line problem, on a finite interval, the semidiscretized IBVPs for the advection equation $q_t = -c \, q_x$ $(c > 0)$ are similar to those for $q_t = +c \, q_x$, except we now apply backward stencils to $q_x(x,t)$ instead of forward ones.
	\end{remark}
	
	\begin{remark}
		For the half-line problem \cite{SDUTM_HL}, the centered discretization also yields a suitable SD-UTM solution that maintains $\mathcal{O}(h^2)$ accuracy through use of the nontrivial symmetry $\nu_1(k) = - k - \pi/h$ from the dispersion relation 
	\begin{equation}
		W(k) = c\,\frac{e^{ikh} - e^{-ikh}}{2h}= \frac{- c\,\sin(kh)}{i h}.
		\label{W_advec2_centered}
	\end{equation}
	In the finite interval problem, however, this is not the case. For the IBVP
	\begin{equation}\begin{dcases}
		q_t = -c\,q_{x},& 0 < x < L,\, t > 0, \\
		q(x,0) = \phi(x),& 0 < x < L,\\
		q(0,t) = u^{0}(t),& t > 0,
		\label{advec2_prob}
	\end{dcases}\end{equation}
	with $c > 0$ and a centered discretization, the global relation is
	\begin{align}
		e^{WT} \hat{q}(k,T) - \hat{q}(k,0) - c\,\left[ \frac{e^{-i k h} f_0 + f_1 - e^{-i k L}\left(e^{-i k h} g_0 + g_1\right)}{2} \right] &= 0,\quad k \in \mathbb{C},
		\label{GR_advec2_centered}
	\end{align}
	so that the ``solution'' contains three unknowns: $f_1(W,T)$ and both $g_j(W,T)$ terms. Figure \ref{advec_cent_W} implies we cannot argue away dependence on all $g_j(W,T)$ terms, since we have regions of exponential growth in both the upper and lower halves of the complex $k$-plane. 
	\begin{figure}[tb]
		\begin{center}
			\def\svgwidth{2.75in}
			\includegraphics[width=0.45\linewidth]{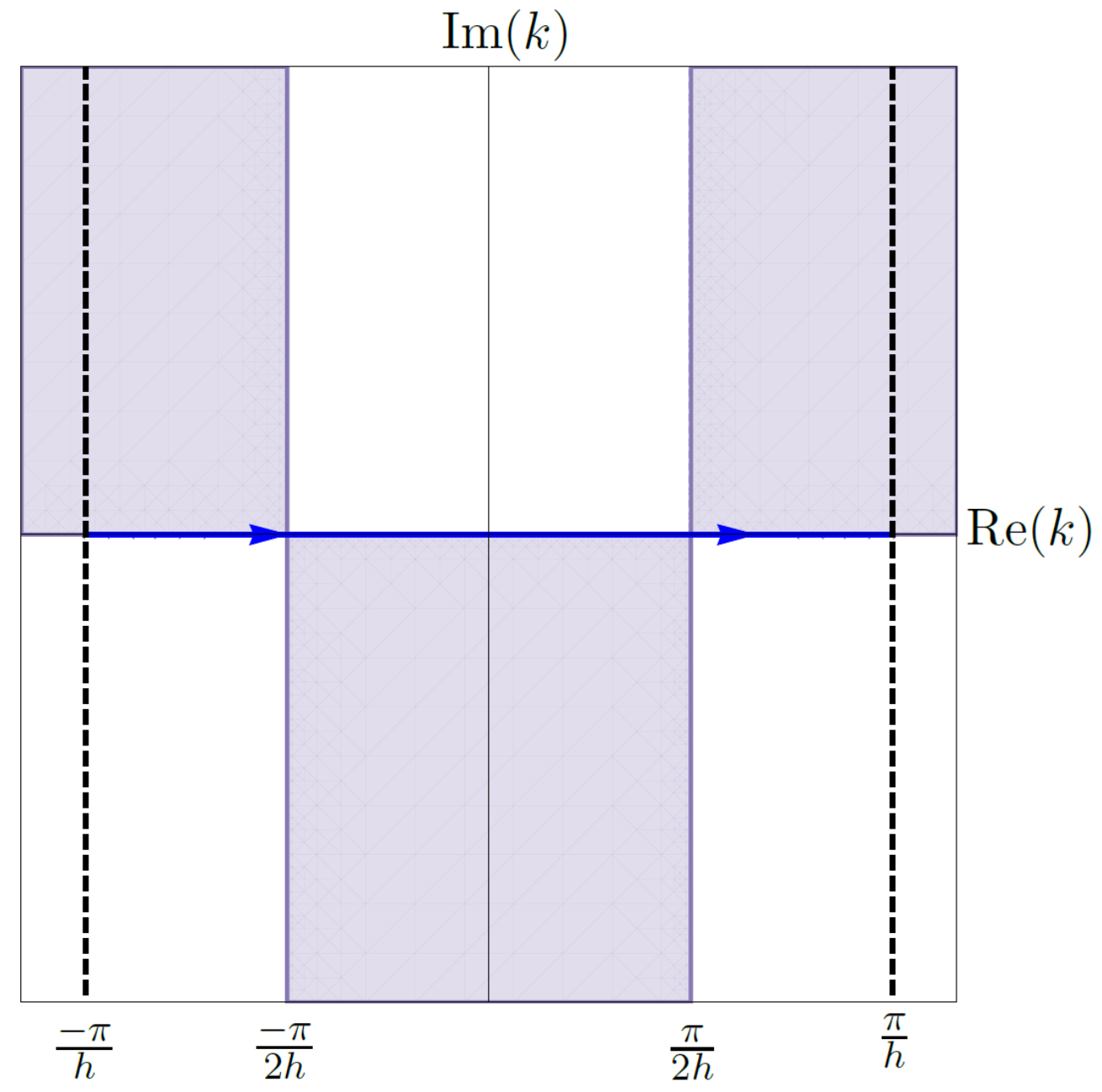}
			\caption{The shaded regions depict where $\text{Re}(-W) \leq 0$ and $e^{-WT}$ is bounded, for the dispersion relation \eqref{W_advec2_centered}.}
			\label{advec_cent_W}
		\end{center}
	\end{figure}
	Deforming the integration paths onto the boundaries of the shaded regions in Figure \ref{advec_cent_W} (see \cite{fokas_paper,fokas_collab,fokas_book,bernard_fokas} and future sections on higher-order discretizations), the global relation \eqref{GR_advec2_centered} with $k$ and $k \rightarrow \nu_1$ provides two equations to remove one $f_j(W,T)$ and one $g_j(W,T)$ terms, say $f_1(W,T)$ and $g_1(W,T)$. We require a third equation that relates $g_0(W,T)$ to at least one of the other $f_j(W,T)$ or $g_1(W,T)$ terms. Unless we have periodic boundary conditions, there is no such relation that does not introduce more unknowns. Hence, the SD-UTM shows that a solution to the centered-discretized IBVP \eqref{advec2_prob} does not exist. As in the higher-order discretization in Section \ref{advec_forward_finiteinterval_2}, we could derive and discretize the Neumann and second-derivative boundary conditions, $q_x(0,t) = -\dot{u}(t)/c$ and $q_{xx}(0,t) = \ddot{u}(t)/c^2$ respectively, given by the PDE from the available Dirichlet condition $u^{(0)}(t)$, with $\dot{u}(t) = d u^{(0)}(t)/dt$ and $\ddot{u}(t) = d ^2u^{(0)}(t)/dt^2$. However, this approach only serves to remove $f_1(W,T)$ from the ``solution,'' impairing the global relation equations with $k$ and $k \rightarrow \nu_1$ to remove the remaining $g_0(W,T)$ and $g_1(W,T)$. We reach a similar conclusion for the centered-discretized IBVP \eqref{advec1_prob} with a Dirichlet condition at $x = L$. Lastly, the excess unknown terms occur near the interval's boundary where there is no prescribed condition, so this result does not change whether or not we take into consideration the known boundary points (the starting and ending index) in the Fourier transform $\hat{q}(k,t)$ definition \eqref{fourier_cont_SD}.
	\end{remark}


\section{The Heat Equation}

\subsection{Centered Discretization of $\bms{q_t = q_{xx}}$ with Dirichlet boundary conditions}\label{heat_centered_finiteinterval_sec}
	Consider the problem
		\begin{equation}\begin{dcases}
		q_t = q_{xx},& 0 < x > L,\, t > 0, \\
		q(x,0) = \phi(x),& 0 < x < L,\\
		q(0,t) = u^{(0)}(t),& t > 0, \\
		q(L,t) = v^{(0)}(t),& t > 0,
		\label{heat_prob}
	\end{dcases}\end{equation}
	with Dirichlet boundary conditions on both sides of the interval. We write the centered-discretized heat equation as
	\begin{equation}
		\dot{q}_n(t) = \frac{q_{n+1}(t) - 2 q_n(t) + q_{n-1}(t)}{h^2},
		\label{heat_centered}
	\end{equation}
	with local relation
	\begin{align}
		\partial_t \left(e^{-iknh} e^{Wt} q_n \right) &= \frac{1}{h^2}\Delta \left(e^{-ik(n-1)h} e^{Wt} q_{n} - e^{-iknh} e^{Wt} q_{n-1} \right),
		\label{LR_heat_centered}
	\end{align}
	and dispersion relation
	\begin{equation}
		W(k) = \frac{2 - e^{ikh} - e^{-ikh}}{h^2} = \frac{2\left( 1 - \cos(kh) \right)}{h^2}.
		\label{W_heat_centered}
	\end{equation}
	The global relation is obtained by summing from $n = 1$ to $n = N$ and integrating in time:
	\begin{align}
		e^{WT} \hat{q}(k,T) - \hat{q}(k,0) - \left[ \frac{e^{-i k h} f_0 - f_1 + e^{-ik L}\left(e^{i k h} g_0 - g_{-1}  \right)}{h} \right] &= 0, \quad k \in \mathbb{C}.
		\label{GR_heat_centered}
	\end{align}
	Inverting, we obtain the ``solution'' formula
	\begin{align}\begin{split}
		q_n(T) &= \frac{1}{2 \pi} \int_{-\pi/h}^{\pi/h} e^{iknh} e^{-WT}\hat{q}(k,0)\,dk + \frac{1}{2 \pi} \int_{-\pi/h}^{\pi/h} e^{iknh} e^{-WT}\left[ \frac{e^{-ikh} f_{0} - f_1}{h} \right]\,dk \\
		&\quad\,+ \frac{1}{2 \pi} \int_{-\pi/h}^{\pi/h} e^{ik(nh-L)} e^{-WT}\left[ \frac{e^{i k h} g_0 - g_{-1}}{h} \right]\,dk,
		\label{soln1_heat_centered}
	\end{split}\end{align}
	which depends on the unknowns $f_1(W,T)$ and $g_{-1}(W,T)$. The dispersion relation \eqref{W_heat_centered} has the trivial $\nu_0(k) = k$ and nontrivial $\nu_1(k) = -k$ symmetries, up to periodic copies. To remove both unknowns, we need to deform the integration path of the second integral with the $f_{j}(W,T)$ terms away from the integration path of the $g_{j}(W,T)$ terms. This results in two equations to solve for two unknowns: 
	\begin{equation}\begin{dcases}
		0 = e^{WT} \hat{q}(k,T) - \hat{q}(k,0) - \left[ \frac{e^{-i k h} f_0 - f_1 + e^{-ik L}\left(-g_{-1} + e^{i k h} g_0 \right)}{h} \right], \\
		0 = e^{WT} \hat{q}(-k,T) - \hat{q}(-k,0) - \left[ \frac{e^{i k h} f_0 - f_1 + e^{ik L}\left(-g_{-1} + e^{-i k h} g_0 \right)}{h} \right],
		\label{heat_centered_system}
	\end{dcases}\end{equation}
	both valid for $k \in \mathbb{C}$. 
	
\begin{figure}[tb]	
		\raggedleft
		\begin{subfigure}[t]{.45\textwidth}
			\begin{center}
			\def\svgwidth{2.75in}
			\vspace{10pt}\includegraphics[width=1\linewidth]{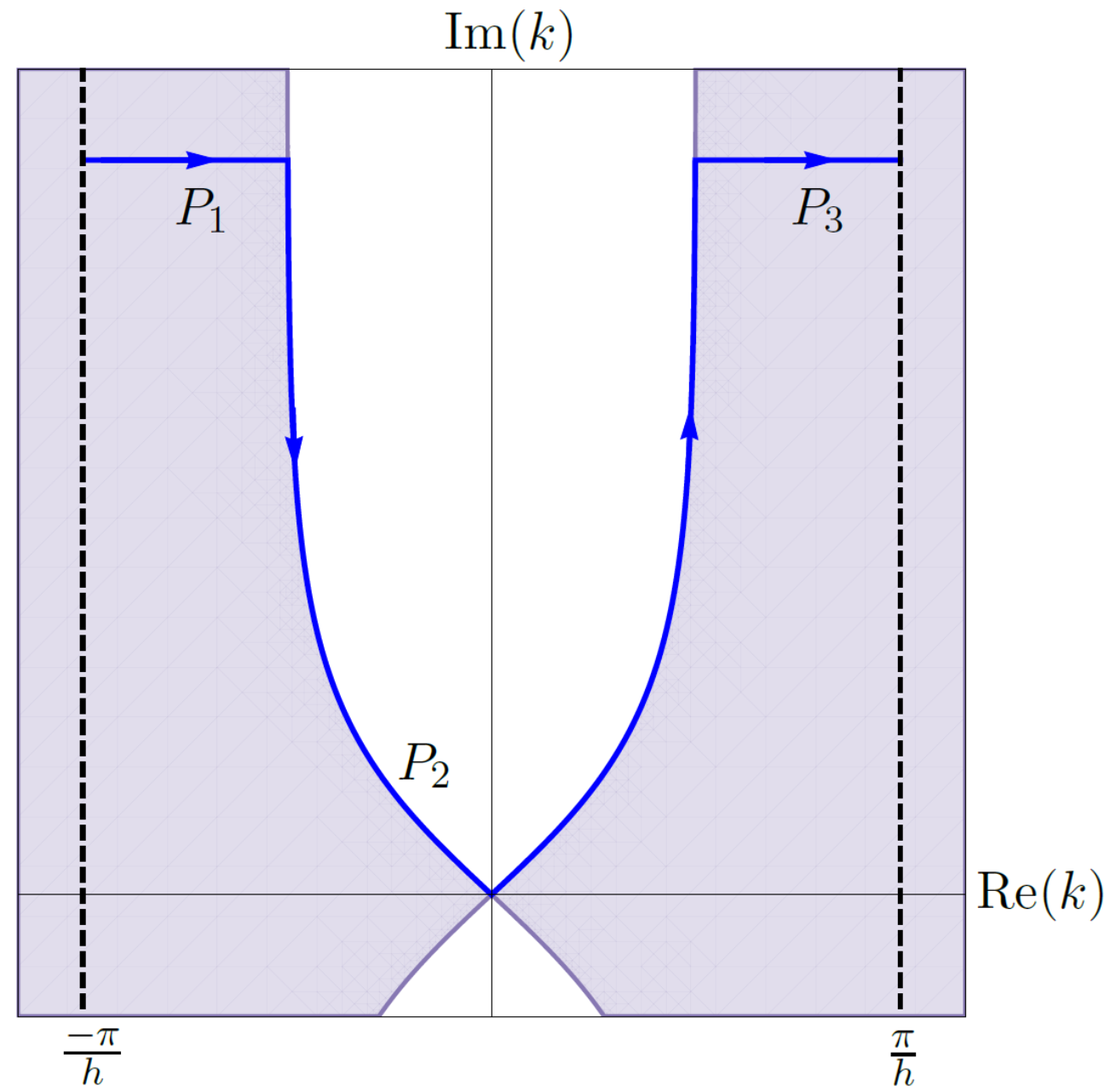}\vspace{10pt}
			\caption{}
			\label{paths_heat_centered_FI_P}
		\end{center}
		\end{subfigure}\hfill 
		\begin{subfigure}[t]{.45\textwidth}
			\begin{center}
			\def\svgwidth{2.75in}
			\vspace{10pt}\includegraphics[width=1\linewidth]{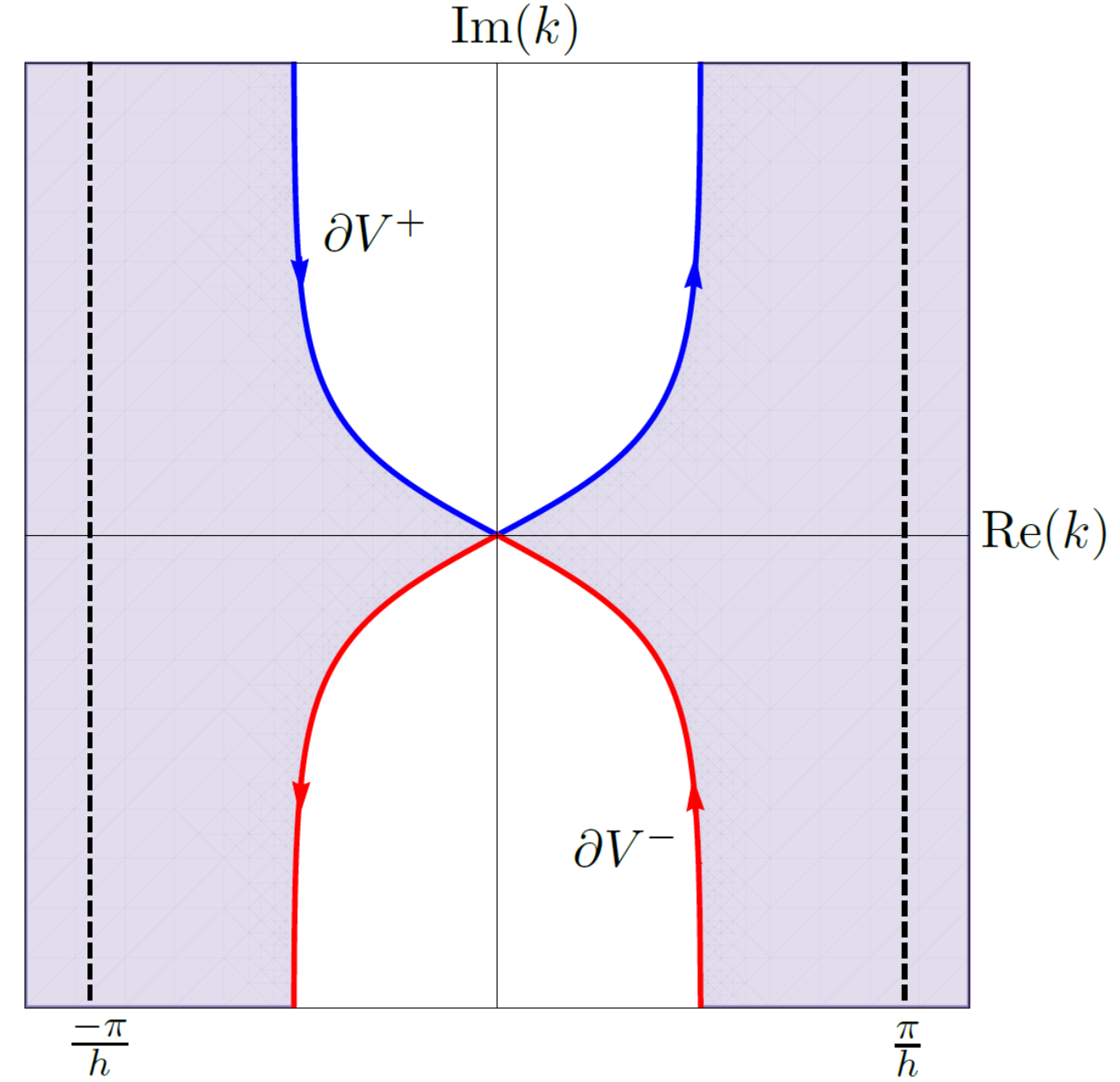}\vspace{10pt}
			\caption{}
			\label{heat_centered_paths_V}
		\end{center}
		\end{subfigure}	
		\caption{(a) The shaded regions depict where $\text{Re}(-W) \leq 0$ and $e^{-WT}$ is bounded, for the dispersion relation \eqref{W_heat_centered}. The integration paths that constitute $P$ are also shown. (b) The integration paths $\partial V^{\pm}$.}
		\label{heat_centered_paths_fig}
	\end{figure}

\begin{figure}[tb]	
		\raggedleft
		\begin{subfigure}[t]{.45\textwidth}
			\begin{center}
			\def\svgwidth{2.75in}
			\vspace{10pt}\includegraphics[width=1\linewidth]{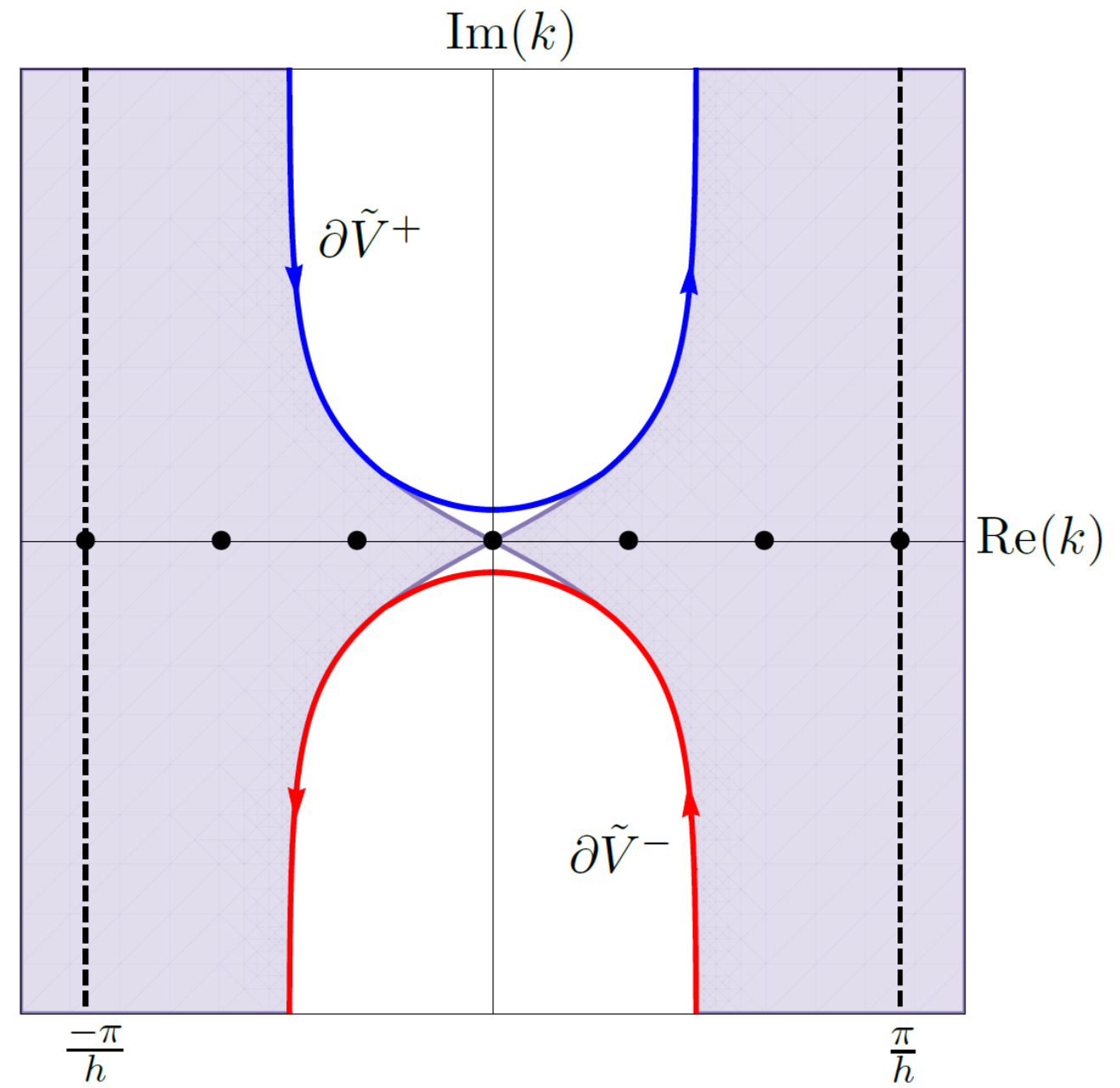}\vspace{10pt}
			\caption{}
			\label{new paths_heat_centered_FI}
		\end{center}
		\end{subfigure}\hfill 
		\begin{subfigure}[t]{.45\textwidth}
			\begin{center}
			\def\svgwidth{2.75in}
			\vspace{10pt}\includegraphics[width=1\linewidth]{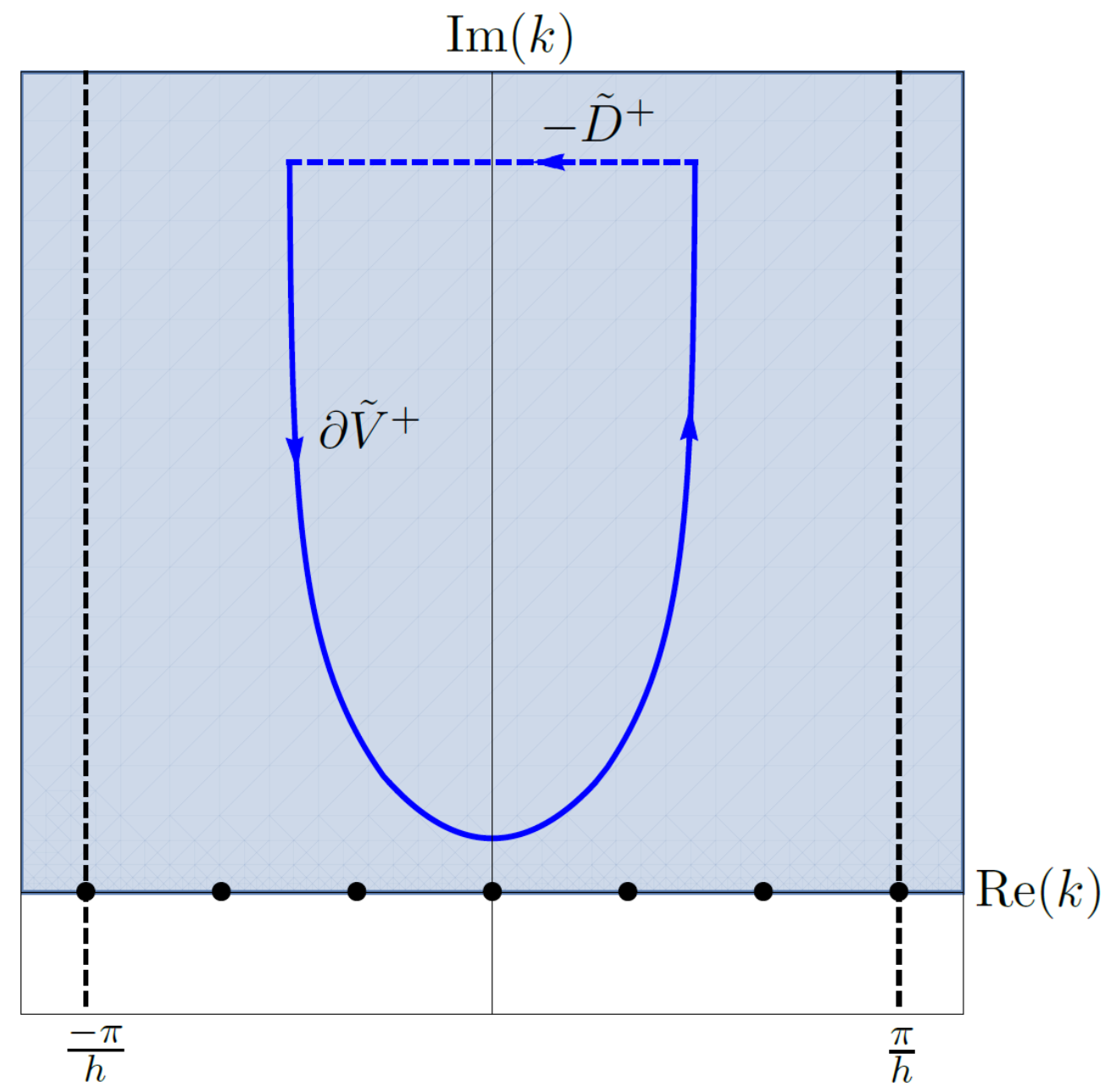}\vspace{10pt}
			\caption{}
			\label{new paths_heat_centered_FI_D}
		\end{center}
		\end{subfigure}	
		\caption{(a) The integration paths $\partial \tilde{V}^{\pm}$ deformed away from the origin. (b) Deforming $\partial \tilde{V}^{+}$ to $\tilde{D}^{+}$ in the upper-half plane. Without decay from $e^{-WT}$, the shaded region depicts where $e^{i k n h}$ decays.}
		\label{new paths_heat_centered_FI_D_fig}
	\end{figure}

	Let us deform the $f_j(W,T)$ terms to the upper-half plane in order to abide with well-posedness as $h \rightarrow 0$. We introduce 
	$$V^{\pm} = \left\{k \in \mathbb{C}^{\pm} \, \Big| \, \text{Re}(-W) \leq 0 \right\}.$$
	With hindsight, we define the integration path $P = P_1 + P_2 + P_3$, shown in Figure \ref{paths_heat_centered_FI_P}, where the two horizontal paths $P_{1,3}$ are at height $\text{Im}(k) = R > 0$ above the real line and $P_2$ is on the boundary of $V^{+}$ up to $\text{Im}(k) = R$.
	Using periodicity, we deform the second integral of ``solution'' \eqref{soln1_heat_centered} to $P$, so that
	$$\frac{1}{2 \pi} \int_{-\pi/h}^{\pi/h} e^{iknh} e^{-WT}\left( \frac{e^{-ikh} f_{0} - f_1}{h} \right)\,dk = \frac{1}{2 \pi} \int_{P} e^{iknh} e^{-WT}\left( \frac{e^{-ikh} f_{0} - f_1}{h} \right)\,dk.$$
	Since $P_{1,3}$ are in regions of exponential decay, we let $R \rightarrow \infty$, so that the integrals on $P_{1,3}$ vanish, the endpoints of $P_2$ are extended to $+ i \infty$ approaching the vertical asymptotes $\text{Re}(k) = \pm \pi/(2h)$, and $\lim_{R\rightarrow \infty} P_2 = \partial V^{+}$ (the entire boundary of $V^+$ in the upper-half plane). Figure \ref{heat_centered_paths_V} shows $\partial V^{\pm}$, where $\partial V^{-}$ is obtained in a similar fashion for the third integral of ``solution'' \eqref{soln1_heat_centered} containing the $g_j(W,T)$ terms. Hence,
	\begin{align}
	\begin{split}
		q_n(T) &= \frac{1}{2 \pi} \int_{-\pi/h}^{\pi/h} e^{iknh} e^{-WT}\hat{q}(k,0)\,dk + \frac{1}{2 \pi} \int_{\partial V^+} e^{iknh} e^{-WT}\left(\frac{e^{-ikh} f_{0} - f_1}{h} \right)\,dk\\
		&\quad\, - \frac{1}{2 \pi} \int_{\partial V^-} e^{ik(nh-L)} e^{-WT}\left(\frac{e^{i k h} g_0 - g_{-1}}{h} \right)\,dk.
		\label{soln2_heat_centered_FI}
	\end{split}
	\end{align} 
	Note the path directions, which have introduced a minus sign on the third integral. Now that the $f_j(W,T)$ and $g_j(W,T)$ terms are on different integration paths, we solve the two global relation equations \eqref{heat_centered_system} for $f_{1}(W,T)$ and $g_{-1}(W,T)$, obtaining
	\begin{equation}\hspace{-35pt}\begin{dcases}
		\frac{e^{-i k h} f_0 - f_1}{h} = \frac{ h\left[ \hat{q}(-k,0) - e^{2 i k L} \hat{q}(k,0) + e^{2 i k L} e^{WT} \hat{q}(k,T)  -  e^{WT} \hat{q}(-k,T)\right] + 2 i \sin(kh) \left(f_0 - e^{i k L} g_0 \right)}{h \left(e^{2 i k L } - 1\right)}, \\[5pt]
		\frac{e^{-ik L}\left(-g_{-1} + e^{i k h} g_0 \right)}{h} = \frac{ h \left[\hat{q}(k,0) - \hat{q}(-k,0) + e^{WT}\hat{q}(-k,T) - e^{WT}\hat{q}(k,T)\right]  - 2 i \sin(kh) \left(f_0 - e^{i k L } g_0 \right)}{h \left(e^{2 i k L } -1\right)}.
		\label{heat_centered_system_solved}
	\end{dcases}\end{equation}
	Both left-hand sides are analytic in $k$, thus the roots of the denominator, $k_{\ell} = \pi \ell/L$, $\ell \in \mathbb{Z}$, are removable singularities, including at the ends of the interval $[-\pi/h ,\pi/h]$ and at the origin $k_{0} = 0$. Since we are only interested in this interval, we can restrict $\ell$ to $-(N+1) \leq \ell \leq  N+1$, using $L = (N+1)h$. For any finite $h$, the number of singularities is finite and increasing as $h \rightarrow 0$. Our integration paths are off the real line except at the origin. To avoid passing through the removable singularity $k_{0}$, it is convenient to deform $\partial V^{\pm}$ to $\partial \tilde{V}^{\pm}$, which is entirely off the real line as depicted in Figure \ref{new paths_heat_centered_FI}. Of course, $\partial \tilde{V}^{\pm}$ now crosses into the unshaded regions where $e^{-WT}$ grows, but this growth is bounded on this segment. 
	On $\partial \tilde{V}^{\pm}$, ``solution'' \eqref{soln2_heat_centered_FI} becomes
	\begin{align}
	\begin{split}
		q_n(T) &= \frac{1}{2 \pi} \int_{-\pi/h}^{\pi/h} e^{iknh} e^{-WT}\hat{q}(k,0)\,dk + \frac{1}{2 \pi} \int_{\partial \tilde{V}^+} e^{iknh} e^{-WT}\left(\frac{e^{-ikh} f_{0} - f_1}{h} \right)\,dk\\
		&\quad\, - \frac{1}{2 \pi} \int_{\partial \tilde{V}^-} e^{ik(nh-L)} e^{-WT}\left(\frac{e^{i k h} g_0 - g_{-1}}{h} \right)\,dk.
		\label{soln3_heat_centered_FI}
	\end{split}
	\end{align}
	Using \eqref{heat_centered_system_solved},
	\begin{align}
	\begin{split}
		q_n(T) &= \frac{1}{2 \pi} \int_{-\pi/h}^{\pi/h} e^{iknh} e^{-WT}\hat{q}(k,0)\,dk \\
		&\quad\, + \frac{1}{2 \pi} \int_{\partial \tilde{V}^+} e^{i k nh} e^{-WT} \left[ \frac{ \hat{q}(-k,0) - e^{2 i k L} \hat{q}(k,0) }{e^{2 i k L} - 1} + \frac{ 2 i \sin(kh) \left(f_0 - e^{i k L} g_0 \right)}{h \left(e^{2 i k L } - 1\right)} \right] \,dk \\
		&\quad\, - \frac{1}{2 \pi} \int_{\partial \tilde{V}^-} e^{ik nh} e^{-WT} \left[ \frac{ \hat{q}(k,0) - \hat{q}(-k,0)}{e^{2 i k L} - 1} - \frac{2 i \sin(kh) \left(f_0 - e^{i k L } g_0 \right)}{h \left(e^{2 i k L }  -1\right)} \right]\,dk \,\,+\,\, S(n),
		\label{soln4_heat_centered_FI}
	\end{split}
	\end{align}
	where
	\begin{align*}
		S(n) &= \frac{1}{2 \pi} \int_{\partial \tilde{V}^+} e^{i k nh} \left[ \frac{e^{2 i k L} \hat{q}(k,T)  - \hat{q}(-k,T) }{e^{2 i k L} - 1} \right] \,dk - \frac{1}{2 \pi} \int_{\partial \tilde{V}^-} e^{ik nh} \left[ \frac{ \hat{q}(-k,T) - \hat{q}(k,T)}{e^{2 i k L} - 1}\right]\,dk.
	\end{align*}
	
	We wish to determine the contributions from $S(n)$. Since the exponential $e^{-WT}$ is not present, our aim is to close the contours using a path at infinity. First, we truncate the infinite paths $\partial \tilde{V}^{\pm}$ and close the curves by introducing
		$$\tilde{D}^{\pm} = \left\{k \in \mathbb{C} \, \Big| \, \frac{-\pi}{2 h} \leq \text{Re}(k) \leq \frac{\pi}{2 h} \, \text{ and } \, \text{Im}(k) =\pm R \right\},$$
		with $R > 0$, so that $\tilde{D}^{\pm}$ is a horizontal interval above/below $\partial \tilde{V}^{\pm}$ in the complex $k$-plane, as shown in Figure \ref{new paths_heat_centered_FI_D} for $\tilde{D}^+$, and $\lim_{R \rightarrow \infty} \tilde{D}^{\pm} = \pm \partial \tilde{V}^{\pm}$.
	
	Now,
	$$ S(n) = \lim_{R \rightarrow \infty} \left(\frac{1}{2 \pi} \int_{\tilde{D}^+} e^{i k nh} \left[ \frac{e^{2 i k L} \hat{q}(k,T)  - \hat{q}(-k,T) }{e^{2 i k L} - 1} \right] \,dk + \frac{1}{2 \pi} \int_{\tilde{D}^-} e^{ik nh} \left[ \frac{ \hat{q}(-k,T) - \hat{q}(k,T)}{e^{2 i k L} - 1}\right]\,dk \right). $$
	Notice the sign change for $\tilde{D}^{-}$. For the first integral on $\tilde{D}^+$, $\hat{q}(k,T)$ grows exponentially as $R \rightarrow \infty$. However, we rewrite the first term as
	\begin{align*}
		\hspace{-30pt}\frac{1}{2 \pi} \int_{\tilde{D}^+} \frac{e^{i k (nh + 2 L)} }{e^{2 i k L} - 1} \hat{q}(k,T)  \,dk = \frac{1}{2 \pi} \int_{\tilde{D}^+} \frac{e^{i k (nh + 2 L)} }{e^{2 i k L} - 1} \left[h \sum_{m=1}^{N} e^{-ikmh} q_m(T)  \right]  \,dk = \frac{h}{2 \pi} \sum_{m=1}^{N} q_m(T) \left[ \int_{\tilde{D}^+} \frac{e^{i k (n - m + 2N + 2)h} }{e^{2 i k L} - 1} \,dk \right],
	\end{align*}
	after using $L = (N+1)h$. For all $n$ and $m$, $n - m + 2N + 2 > 0$. Letting $R \rightarrow \infty$ on $\tilde{D}^+$ implies $e^{i k (n - m + 2N + 2)h} \rightarrow 0$ and $e^{2 i k L} \rightarrow 0$. In this limit, the integrand approaches zero and we recover the original integration paths, so that 
	$$\frac{1}{2 \pi} \int_{\partial \tilde{V}^+} \frac{e^{i k (nh + 2 L)} }{e^{2 i k L} - 1} \hat{q}(k,T)  \,dk = 0,$$
	for all $n$. The second term with $\hat{q}(-k,T)$ on $\tilde{D}^+$ similarly goes to zero as $R \rightarrow \infty$. The third term is rewritten as
	\begin{align*}
		\frac{1}{2 \pi} \int_{\tilde{D}^-} \frac{ e^{ik nh}}{e^{2 i k L} - 1} \hat{q}(-k,T)\,dk &= \frac{h}{2 \pi} \sum_{m=1}^{N} q_m(T) \left[ \int_{\tilde{D}^-} \frac{e^{i k (n + m)h} }{e^{2 i k L} - 1} \,dk\right].
	\end{align*}
	On $\tilde{D}^-$, $e^{i k (n + m)h}$ and $e^{2i k L}$ grow as $R \rightarrow \infty$ for all $n$ and $m$, but
	$$\frac{e^{i k (n + m)h} }{e^{2 i k L} - 1} \sim \frac{e^{i k (n + m)h} }{e^{2 i k L}} = e^{i k (n + m - 2N - 2) h}.$$
	Since $n < N+1$ and $m < N+1$, $n + m - 2N - 2 < 0$, and the integrand approaches zero, as $R \rightarrow \infty$. Thus, 
	$$\frac{1}{2 \pi} \int_{\partial \tilde{V}^-} \frac{ e^{ik nh}}{e^{2 i k L} - 1} \hat{q}(-k,T)\,dk = 0,$$
	for all $n$. Similarly, the fourth term is zero. Hence, $ S(n) = 0$ and the final representation for the solution to the finite-interval problem for the heat equation \eqref{heat_centered} with Dirichlet boundary conditions is
	\begin{align}
	\begin{split}
		q_n(T) &= \frac{1}{2 \pi} \int_{-\pi/h}^{\pi/h} e^{iknh} e^{-WT}\hat{q}(k,0)\,dk \\
		&\quad\, + \frac{1}{2 \pi} \int_{\partial \tilde{V}^+} e^{i k nh} e^{-WT} \left[ \frac{ \hat{q}(-k,0) - e^{2 i k L} \hat{q}(k,0) }{e^{2 i k L} - 1} + \frac{ 2 i \sin(kh) \left(f_0 - e^{i k L} g_0 \right)}{h \left(e^{2 i k L } - 1\right)} \right] \,dk \\
		&\quad\, - \frac{1}{2 \pi} \int_{\partial \tilde{V}^-} e^{ik nh} e^{-WT} \left[ \frac{ \hat{q}(k,0) - \hat{q}(-k,0)}{e^{2 i k L} - 1} - \frac{2 i \sin(kh) \left(f_0 - e^{i k L } g_0 \right)}{h \left(e^{2 i k L }  -1\right)} \right]\,dk.
		\label{soln_heat_centered_FI}
	\end{split}
	\end{align}
	
\begin{figure}[tb]
		\begin{center}
			\def\svgwidth{2.75in}
			\vspace{10pt}\includegraphics[width=0.45\linewidth]{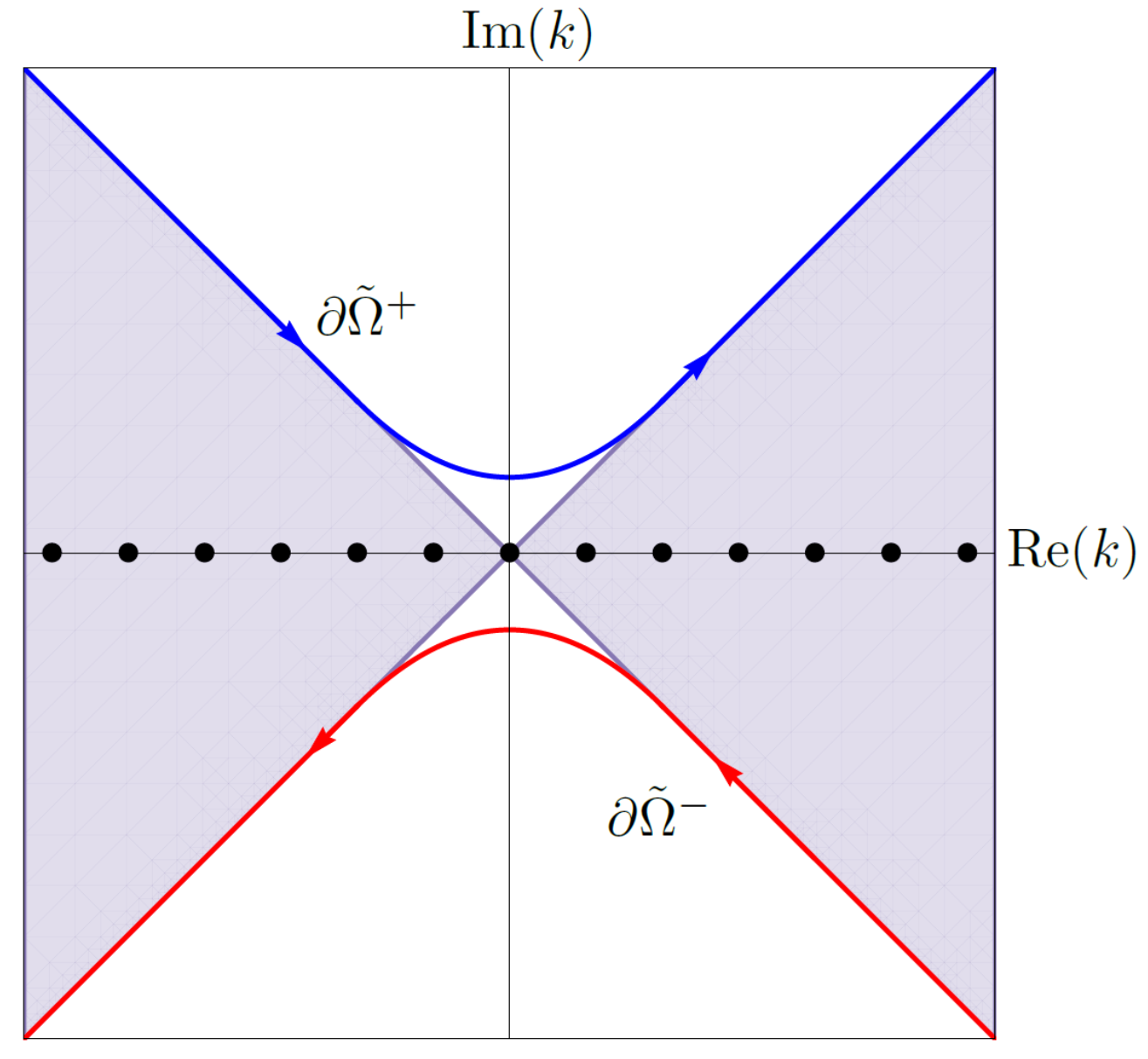}\vspace{10pt}
			\caption{The shaded regions depict where $\text{Re}(-\tilde{W}) \leq 0$ and $e^{\tilde{W}T}$ is bounded, with the boundaries $\partial \tilde{\Omega}^{\pm}$ approaching $\partial \Omega^{\pm}$ asymptotically as $|k| \rightarrow \infty$.}
			\label{new_paths_heat_cont_FI}
		\end{center}
	\end{figure}

	For the IBVP \eqref{heat_prob}, using the continuous UTM \cite{bernard_fokas}, we find the dispersion relation $\tilde{W}(k) = k^2$ and define 
	$$ \Omega^{\pm} = \left\{k \in \mathbb{C}^{\pm} \, \Big| \, \text{Re}(-k^2) \leq 0 \right\},$$
	to use the integration paths $\partial \tilde{\Omega}^{\pm}$ as illustrated in Figure \ref{new_paths_heat_cont_FI}. As in the semidiscrete case, there is a (removable) singularity at the origin, so $\partial \Omega^{\pm}$ is deformed to $\partial \tilde{\Omega}^{\pm}$.
	The solution is
	\begin{align}\begin{split}
		q(x,T) &= \frac{1}{2 \pi} \int_{-\infty}^{\infty} e^{ikx} e^{-\tilde{W}T}\hat{q}(k,0)\,dk \\
		&\quad\, + \frac{1}{2 \pi} \int_{\partial \tilde{\Omega}^+} e^{i k x} e^{-\tilde{W}T} \left[ \frac{ \hat{q}(-k,0) - e^{2 i k L} \hat{q}(k,0) + 2 i k \left(F_0 - e^{ikL} G_0\right)}{e^{2 i k L} - 1} \right] \,dk \\
		&\quad\, - \frac{1}{2 \pi} \int_{\partial \tilde{\Omega}^-} e^{ik x} e^{-\tilde{W}T} \left[ \frac{ \hat{q}(k,0) - \hat{q}(-k,0) - 2 i k \left(F_0 - e^{ikL} G_0\right)}{e^{2 i k L} - 1} \right]\,dk.
		\label{soln_heat_cont_FI}
	\end{split}\end{align}
	Taking the continuum limit, it is clear that $\lim_{h\rightarrow 0} W(k) = \tilde{W}(k)$ and thus $\lim_{h\rightarrow 0}\partial \tilde{V}^{\pm} = \partial \tilde{\Omega}^{\pm}$, as well. In addition, the coefficient of the boundary terms from either $\partial \tilde{V}^{\pm}$ integral converge to $ 2 i k /(e^{2 i k L } - 1)$ with $\lim_{h\rightarrow 0} \left(f_0 - e^{i k L} g_0 \right) = F_0 - e^{i k L} G_0 $, so that the SD-UTM solution \eqref{soln_heat_centered_FI} converges to the continuous UTM solution \eqref{soln_heat_cont_FI}.
	
	Unlike for the heat equation on the half-line \cite{SDUTM_HL} and the previous advection IBVPs, the integration paths of the last two integral terms in \eqref{soln_heat_centered_FI} are off the real line, avoiding the integrands' simple poles. To numerically evaluate the integrals, we design any contour path that is within the shaded regions of Figure \ref{new paths_heat_centered_FI} and off the real line, with endpoints that have real part $\pm \pi/h$. For computational purposes, having paths off of the boundaries $\partial \tilde{V}^{\pm}$ is preferred in order to have some exponential decay. Since we are taking $h \ll 1$ in practical settings, we must orient our contours so that in the continuum limit, the semidiscrete solution converges to the continuous one. 
	

	\subsubsection{\textbf{Series Representation}}\label{heat_centered_series} 
\begin{figure}[tb]
		\begin{center}
			\def\svgwidth{2.75in}
			\vspace{10pt}\includegraphics[width=0.45\linewidth]{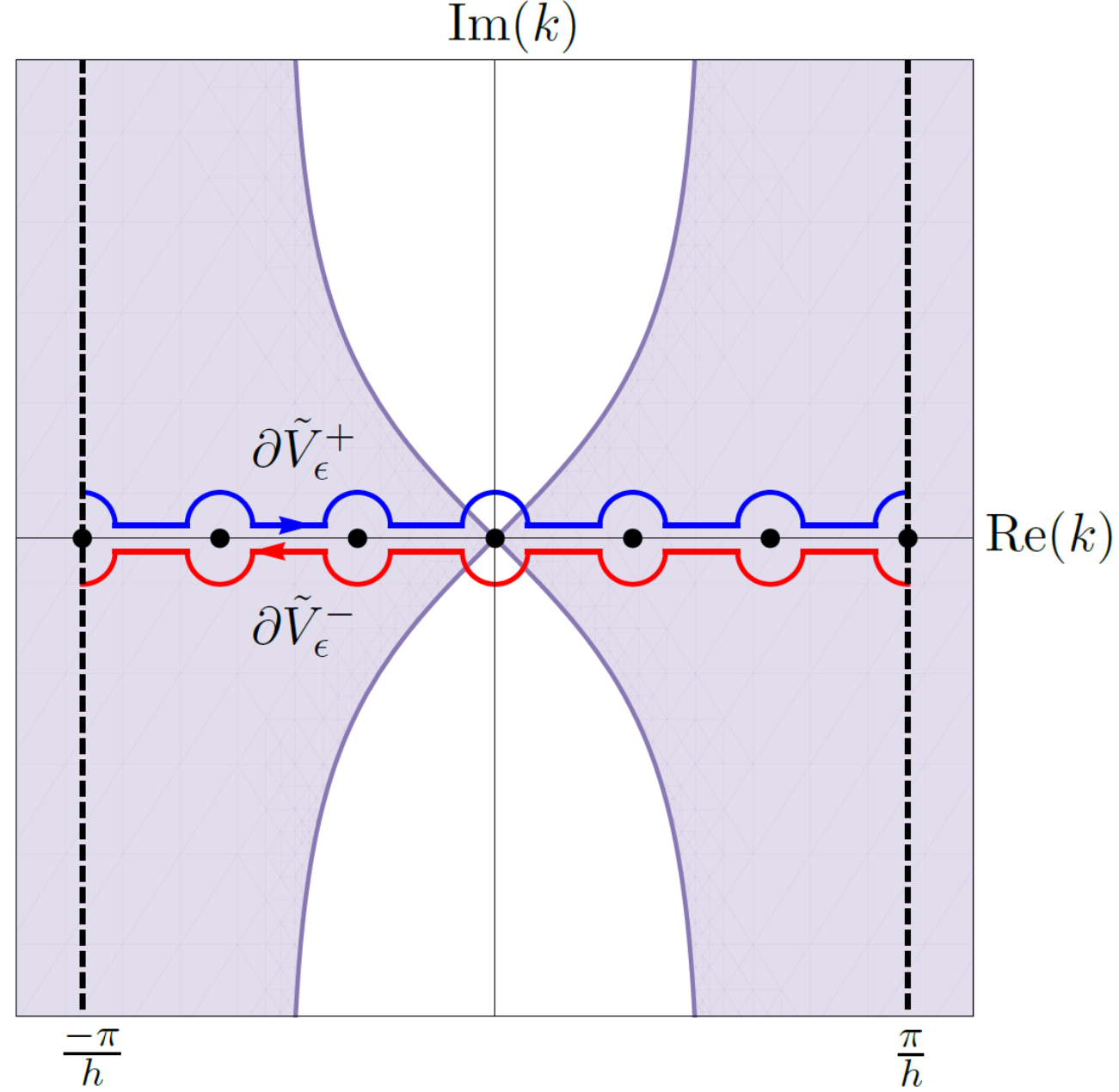}\vspace{10pt}
			\caption{The semi-circle integration paths $\partial \tilde{V}^{\pm}_{\epsilon}$ around singularities on the real line.}
			\label{eps_paths_heat_centered_FI}
		\end{center}
	\end{figure}

	To bypass complex integration paths, we derive a series representation equivalent to \eqref{soln_heat_centered_FI}, as before. For any $n$, we deform the first integral term to $\partial \tilde{V}^+$, since $\hat{q}(k,0)$ is valid for all $k \in \mathbb{C}$ and the bounds on the integral with respect to $k$ are finite. Combining with the other initial condition terms on $\partial \tilde{V}^+$, solution \eqref{soln_heat_centered_FI} becomes
	\begin{align}
		\begin{split}
		q_n(T) &= \frac{-1}{2 \pi} \left( \int_{\partial \tilde{V}^+} + \int_{\partial \tilde{V}^-} \right) A(n,T,k)\,dk + \frac{i}{ \pi h} \left( \int_{\partial \tilde{V}^+} + \int_{\partial \tilde{V}^-} \right) B(n,T,k)\,dk,
		\label{soln_heat_centered_FI_series1}
	\end{split}\end{align}
	with
	$$A(n,T,k) = e^{iknh}  e^{-WT}\left[ \frac{ \hat{q}(k,0) - \hat{q}(-k,0)}{e^{2 i k L} - 1} \right] \quad \text{ and } \quad B(n,T,k) = e^{iknh} e^{-WT} \left[ \frac{ \sin(kh) \left(f_0 - e^{i k L} g_0 \right)}{e^{2 i k L } - 1}  \right]. $$
	Next, we deform the paths back to the real line, excluding the $2(N+1) + 1$ singularities on $[-\pi/h, \pi/h]$ using half-circles with radius $\epsilon$ smaller than half the distance between singularities, see Figure \ref{eps_paths_heat_centered_FI}.
	The horizontal line segments for both $\partial \tilde{V}^{\pm}_{\epsilon}$ are on the real line, but are drawn above and below for illustrative purposes. We obtain
	\begin{align}
	\begin{split}
		q_n(T) &= \frac{-1}{2 \pi} \left(\int_{\partial \tilde{V}^+_{\epsilon}} + \int_{\partial \tilde{V}^-_{\epsilon}} \right) A(n,T,k)\,dk + \frac{i}{ \pi h} \left(\int_{\partial \tilde{V}^+_{\epsilon}} + \int_{\partial \tilde{V}^-_{\epsilon}} \right)B(n,T,k)\,dk.
		\label{eps_soln_heat_centered_FI}
		\end{split}
	\end{align}
	Taking the limit of \eqref{eps_soln_heat_centered_FI} as $\epsilon \rightarrow 0$, the integrals give rise to residue contributions and principal value integrals. For instance, the first term above becomes
	\begin{align*}
		\frac{-1}{2 \pi} \int_{\partial \tilde{V}^+_{\epsilon}} A(n,T,k) \, dk &= \frac{i}{4} \left[ \,\underset{k = -\pi/h}{\text{Res}} A(n,T,k) + \,\underset{k = \pi/h}{\text{Res}} A(n,T,k)\right] + \frac{i}{2}  \sum_{\ell = -N}^{N} \,\underset{k = k_{\ell}}{\text{Res}} \,A(n,T,k) \\
		& \quad\, - \frac{1}{2 \pi} \dashint_{-\pi/h}^{\pi/h} A(n,T,k) \, dk.
	\end{align*}
	We have two contributions of the quarter-circle contours from the singularities on the edges of the integration path plus contributions from the half-circle contours from the inner singularities. Similarly, for the second term in \eqref{eps_soln_heat_centered_FI},
	\begin{align*}
		\frac{-1}{2 \pi} \int_{\partial \tilde{V}^-_{\epsilon}} A(n,T,k) \, dk &= \frac{i}{4 } \left[  \,\underset{k = -\pi/h}{\text{Res}} A(n,T,k) +  \,\underset{k = \pi/h}{\text{Res}} A(n,T,k)\right] + \frac{i}{2}  \sum_{\ell = -N}^{N} \,\underset{k = k_{\ell}}{\text{Res}} \,A(n,T,k) \\
		&\quad\, + \frac{1}{2 \pi} \dashint_{-\pi/h}^{\pi/h} A(n,T,k) \, dk.
	\end{align*}
	Combining the two rewrites above cancels the integrals, so that 
	\begin{align*}
		\hspace{-10pt}\frac{-1}{2 \pi} \left( \int_{\partial \tilde{V}^+_{\epsilon}} + \int_{\partial \tilde{V}^-_{\epsilon}} \right) A(n,T,k) \, dk &=  i \left[\frac{1}{2} \,\underset{k = -\pi/h}{\text{Res}} A(n,T,k) +  \sum_{\ell = -N}^{N} \,\underset{k = k_{\ell}}{\text{Res}} \,A(n,T,k) + \frac{1}{2} \,\underset{k = \pi/h}{\text{Res}} A(n,T,k)\right].
	\end{align*}
	Similarly, the $B(n,T,k)$ integrals in \eqref{eps_soln_heat_centered_FI} are rewritten using residue contributions, so that
	substituting into \eqref{eps_soln_heat_centered_FI} gives an expression with residue contributions only:
	\begin{align}
	\begin{split}
		q_n(T) &= i \left[\frac{1}{2} \,\underset{k = -\pi/h}{\text{Res}} A(n,T,k) +  \sum_{\ell = -N}^{N} \,\underset{k = k_{\ell}}{\text{Res}} \,A(n,T,k) + \frac{1}{2} \,\underset{k = \pi/h}{\text{Res}} A(n,T,k)\right]    \\
		&\quad\, + \frac{1}{h} \left[ \,\underset{k = -\pi/h}{\text{Res}} B(n,T,k) + 2 \sum_{\ell = -N}^{N} \,\underset{k = k_{\ell}}{\text{Res}} \,B(n,T,k) +  \,\underset{k = \pi/h}{\text{Res}} B(n,T,k)\right].
		\label{res_soln_heat_centered_FI}
	\end{split}
	\end{align}
	
	Next, we determine these residues, starting with the $A(n,T,k)$ residues for $\ell = -N-1, \ldots, N+1$, where the lone residues at $k = \pm \pi/h$ or $\ell = \pm (N+1)$ follow trivially. Since we only have simple poles at $k = k_{\ell}$,
	\begin{align*}
		\underset{k = k_{\ell}}{\text{Res}} A(n,T,k) &= e^{ik_{\ell}nh}  e^{-W\left(k_{\ell}\right)T}\left[ \frac{\hat{q}\left(k_{\ell},0\right) - \hat{q}\left(-k_{\ell},0\right)  }{2 i L } \right].
	\end{align*} 
	Using the definitions of $\hat{q}(\pm k,0)$,
	\begin{align*}
		\hat{q}\left(k_{\ell},0\right) - \hat{q}\left(-k_{\ell},0\right) &= -2 i h \sum_{m=1}^{N} \sin\left(\frac{\pi \ell mh}{L} \right)\phi_m = - i L b_{\ell}, \quad\quad b_{\ell} = \frac{2 h}{L} \sum_{m=1}^{N} \sin\left(\frac{\pi \ell mh}{L} \right) \phi_m.
	\end{align*} 
	Note that $b_0 = b_{\pm(N+1)}= 0$, $b_{-\ell} = - b_{\ell}$, $k_{-\ell} =- k_{\ell}$, and $W(k_{\ell}) = W( - k_{\ell})$. Using these observations,
	\begin{align*}
		\sum_{\ell = -N}^{N} \,\underset{k = k_{\ell}}{\text{Res}} A(n,T,k) &=  \frac{-1}{2 }  \sum_{\ell = -N}^{N} e^{ik_{\ell}nh}  e^{-W\left(k_{\ell}\right)T} b_{\ell}= -i \sum_{\ell = 1}^{N}  e^{-W\left(k_{\ell}\right)T} \sin\left(\frac{\pi \ell nh}{L}\right) b_{\ell} .
	\end{align*}
	Similarly,
	\begin{align*}
		\underset{k = k_{\ell}}{\text{Res}} B(n,T,k) & = e^{ik_{\ell}nh} e^{-W_{\ell}T} \sin\left(\frac{\pi \ell h}{L}\right) \frac{ H(W_{\ell},T)}{2 i L}, 
	\end{align*}
	where we have introduced $W_{\ell} \equiv W( k_{\ell})$ for brevity, with 
	$$H(W_{\ell},T)= f_0(W_{\ell},T) + (-1)^{\ell+1} g_0(W_{\ell},T).$$
	As with the residues of $A(n,T,k)$, one can show there are no contributions from the poles at the endpoints and at the origin. After rearranging, 
	\begin{align*}
		\sum_{\ell = -N}^{N} \,\underset{k = k_{\ell}}{\text{Res}} B(n,T,k) &= \frac{1}{L} \sum_{\ell = 1}^{N}  e^{-W_{\ell}T} \sin\left(\frac{\pi \ell h}{L}\right) \sin\left(\frac{\pi \ell nh}{L}\right)H(W_{\ell},T).
	\end{align*}
	
	Returning to \eqref{res_soln_heat_centered_FI}, we recover the classical series solution:
	\begin{align}
		 q_n(T) &=  \sum_{\ell = 1}^{N}  e^{-W_{\ell}T} \sin\left(\frac{\pi \ell nh}{L}\right) \left[ b_{\ell}   + \frac{2}{L h} \sin\left(\frac{\pi \ell h}{L}\right) H(W_{\ell},T) \right] .
		\label{class_soln_heat_centered_FI}
	\end{align}
	From a numerical point of view, this solution representation is favored over the integral representation \eqref{soln_heat_centered_FI}.

\begin{figure}[tb]	
		\raggedleft
		\begin{subfigure}[t]{.45\textwidth}
			\centering
  			\includegraphics[width=1\linewidth]{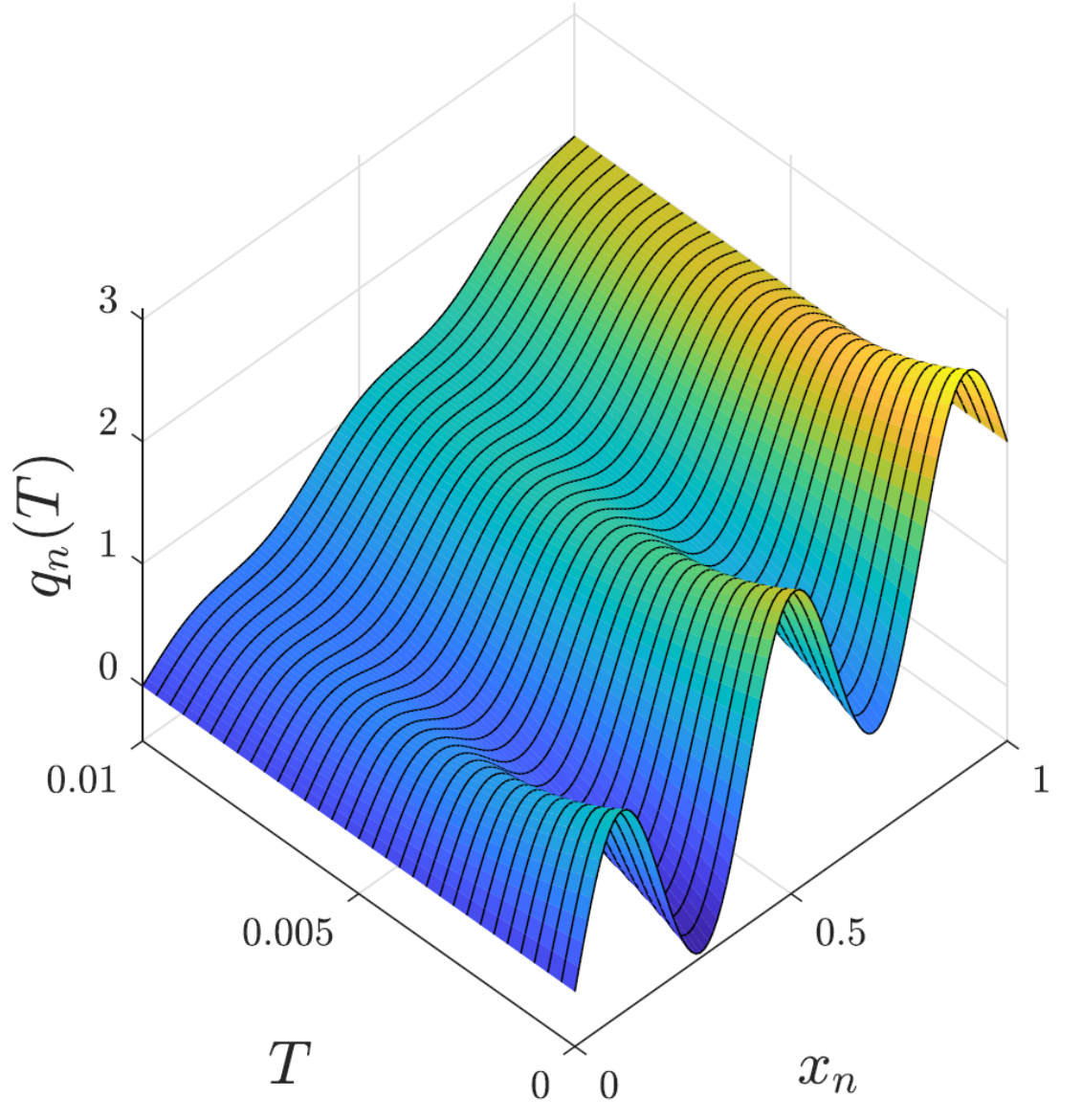}
  			\caption{}
  			\label{heat_xt_UTM2_FI}
		\end{subfigure}\hfill 
		\begin{subfigure}[t]{.45\textwidth}
			\centering
  			\includegraphics[width=1\linewidth]{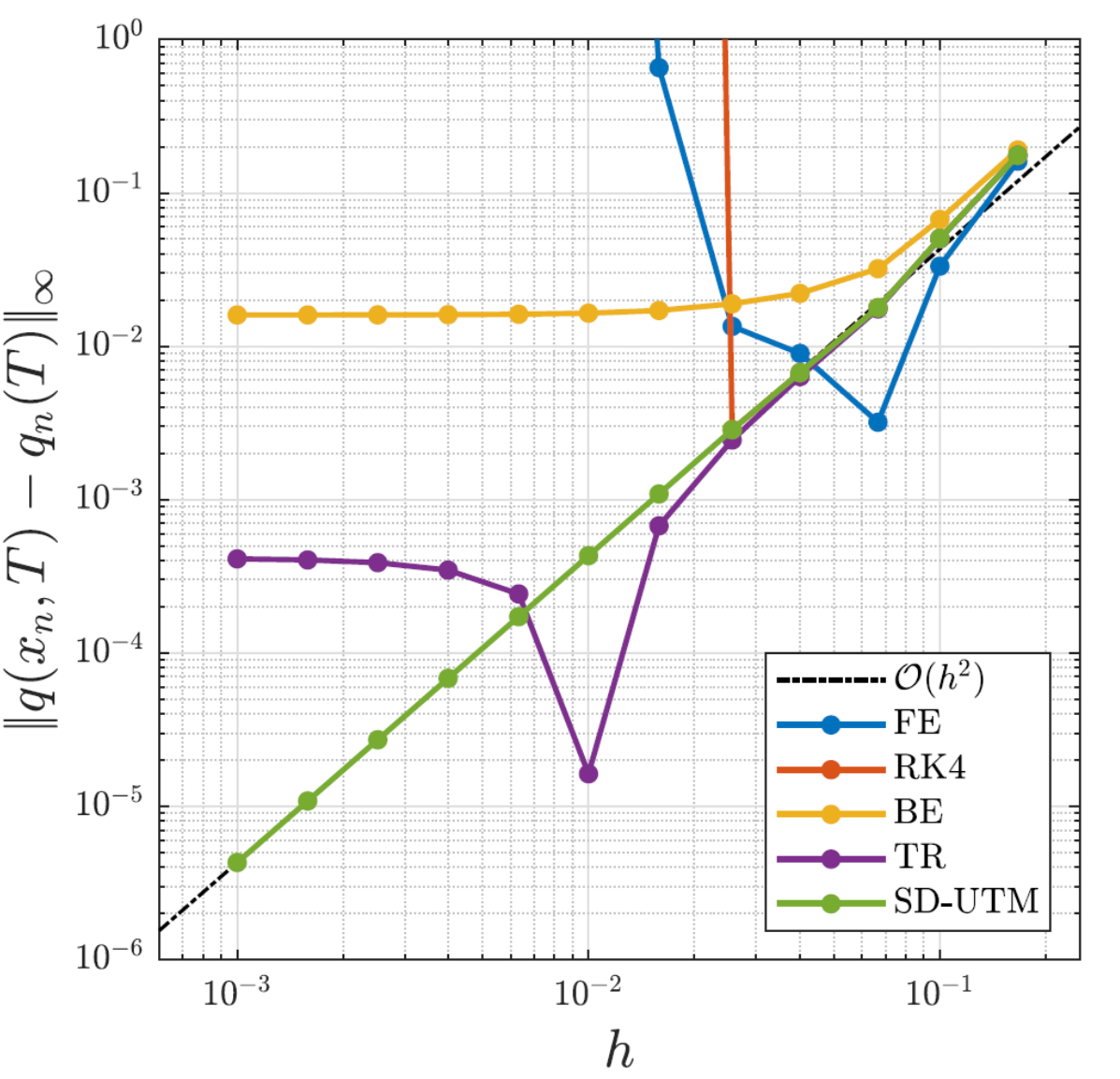}
  			\caption{}
  			\label{heat_errorplot_UTM2_FI}
		\end{subfigure}	
		\caption{(a) The semidiscrete solution \eqref{class_soln_heat_centered_FI} evaluated at various $t$ with $h = 0.01$. (b) Error plot of the semidiscrete solution \eqref{class_soln_heat_centered_FI} and finite-difference schemes relative to the exact solution as $h \rightarrow 0$ with $T = 0.01$ and $\Delta t  = 6.25 \times 10^{-4}$.}
		\label{heat_UTM2_FI}
	\end{figure}

	As an example, the exact solution to the IBVP
\begin{equation}
\begin{dcases}
	q_t = q_{xx},& 0 < x < 1,\, t > 0, \\
	q(x,0) = \phi(x) = 2 x + \sin(5 \pi x),& 0 < x < 1,  \\
	q(0,t) = u^{(0)}(t) = 0,& t > 0,\\
	q(1,t) = v^{(0)}(t) = 2,& t > 0,\\
\end{dcases}
\label{heat_numerical2_HL}
\end{equation}	
	is $q(x,t) = 2 x + \sin(5\pi x) e^{-25 \pi^2 t}$. We have chosen time-independent boundary conditions for simplicity only. We solve the IBVP \eqref{heat_numerical2_HL} with the centered finite-difference methods and the series SD-UTM solution \eqref{class_soln_heat_centered_FI}. Deriving the modified PDE from the centered stencil \eqref{heat_centered}, we find that \eqref{soln_heat_centered_FI}, and hence \eqref{class_soln_heat_centered_FI}, is a fourth-order accurate approximation to the solution of the dissipative PDE 
	\begin{equation}
		p_t = p_{xx} +\frac{h^2}{12} p_{4x}.
		\label{heat_centered_modified_eqn}
	\end{equation} 
	The presence of the higher-order dissipation term $p_{4x}$ causes high-frequency oscillations to be further damped for $t >0$. The original heat equation is also dissipative, but solution \eqref{soln_heat_centered_FI} might overdamp in scenarios where the initial data contains high-frequency oscillations or the boundary condition oscillates in time with large amplitude. Although the dissipation coefficient of $p_{4x}$ is $\mathcal{O}(h^2)$, the overdamping nature can be troublesome for a practical $h \ll 1$ as $t$ increases. This can be counteracted by decreasing $h$. With the SD-UTM solution \eqref{soln_heat_centered_FI}, the left plot of Figure \ref{heat_UTM2_FI} shows the exponential decay as time increases and the right plot shows the expected $\mathcal{O}(h^2)$ error as $h \rightarrow 0$. 
	 The SD-UTM solution outperforms the traditional numerical methods in the continuum limit, where the explicit methods become unstable and the implicit methods are asymptotic to their respective temporal truncation error as $h \rightarrow 0$. The finite-difference solutions have CFL conditions that must be satisfied for stability, but Figure \ref{heat_errorplot_UTM2_FI} shows that there is no such restriction for the SD-UTM to succeed. The dips in this error plot are due to mesh points being placed near stationary points in the solution. With these boundary conditions, we know that $\lim_{t\rightarrow \infty}q(x,t) = 2x$. Hence, mesh points on the initial condition that are near $q = 2x$ tend to stay there as $t$ increases.


\subsection{Centered Discretization of $\bms{q_t = q_{xx}}$ with Neumann boundary conditions} \label{neumann_FI}
	
	We consider the same centered-discretized heat equation \eqref{heat_centered}, now with Neumann boundary conditions at both ends of the interval:
	\begin{equation}\begin{dcases}
		q_t = q_{xx},& 0< x < L,\, t > 0, \\
		q(x,0) = \phi(t),& 0< x < L,\\
		q_x(0,t) = u^{(1)}(t),& t > 0, \\
		q_x(L,t) = v^{(1)}(t),& t > 0.
		\label{heat_prob_N}
	\end{dcases}\end{equation}
	We can discretize the Neumann data with many different stencils, but we show that the SD-UTM restricts which of these are available to be paired with \eqref{heat_centered}.
	
	With the centered discretization \eqref{heat_centered}, we retain the local relation \eqref{LR_heat_centered} and dispersion relation \eqref{W_heat_centered} with nontrivial symmetry $\nu_1(k) = -k$. We cannot use the global relation \eqref{GR_heat_centered}, because we assumed Dirichlet boundary data to obtain it. Without information at $n = 0$ or $n = N+1$, we define the forward transform to start and end at these points:
	$$\hat{q}\left(k,t\right) =h \sum_{n=0}^{N+1} e^{-i k n h} q_{n}(t),$$
	directly affecting the global relation. From the local relation \eqref{LR_heat_centered}, 
	\begin{align}
		e^{WT} \hat{q}(k,T) - \hat{q}(k,0) - \left[ \frac{f_{-1} - e^{ikh} f_{0} + e^{-i k L}\left( g_1 - e^{-ikh} g_0\right)}{h} \right] &= 0, \quad k \in \mathbb{C}.
		\label{GR_heat_centered_N}
	\end{align}
	Solving for $\hat{q}(k,T)$ and inverting, we obtain
	\begin{align}\begin{split}
		q_n(T) &= \frac{1}{2 \pi} \int_{-\pi/h}^{\pi/h} e^{iknh} e^{-WT}\hat{q}(k,0)\,dk \\
		&\quad\, + \frac{1}{2 \pi} \int_{-\pi/h}^{\pi/h} e^{iknh} e^{-WT}\left[ \frac{f_{-1} - e^{ikh} f_{0} + e^{-i k L}\left( g_1 - e^{-ikh} g_0\right)}{h} \right]\,dk,
		\label{soln1_heat_centered_N}
	\end{split}\end{align}
	with unknowns $f_{-1}(W,T)$, $f_{0}(W,T)$, $g_{0}(W,T)$, and $g_{1}(W,T)$. In order to use the global relation \eqref{GR_heat_centered_N} with $k$ and $k\rightarrow -k$, we separate and deform the integration paths of the $f_j(W,T)$ terms from the $g_j(W,T)$ terms. As in \eqref{soln3_heat_centered_FI} with Dirichlet boundary conditions, we deform the boundary terms in ``solution'' \eqref{soln1_heat_centered_N} off the real line to $\partial \tilde{V}^{\pm}$:
	\begin{align}\begin{split}
		q_n(T) &= \frac{1}{2 \pi} \int_{-\pi/h}^{\pi/h} e^{iknh} e^{-WT}\hat{q}(k,0)\,dk + \frac{1}{2 \pi} \int_{\partial \tilde{V}^+} e^{iknh} e^{-WT} \left[ \frac{f_{-1} - e^{ikh} f_{0}}{h} \right]\,dk \\
		&\quad\, - \frac{1}{2 \pi} \int_{\partial \tilde{V}^-} e^{ik(nh-L)} e^{-WT} \left[ \frac{g_1 - e^{-ikh} g_0}{h} \right]\,dk.
		\label{soln2_heat_centered_N}
	\end{split}\end{align}
	Through the time transforms, the global relation \eqref{GR_heat_centered_N} contains boundary nodal information $n = -1, 0$ at the interval's left boundary and $n = N+1, N+2$ at the right, leading us to backward-discretize $q_x \left(0,t\right)$ and forward-discretize $q_x \left(L,t\right)$ with $\mathcal{O}(h)$ stencils to obtain
	\begin{equation}
		\frac{q_0(t) - q_{-1}(t)}{h} = u^{(1)}(t) \quad\quad \text{ and } \quad\quad \frac{q_{N+2}(t) - q_{N+1}(t)}{h} = v^{(1)}(t).
		\label{heat_N_stencils}
	\end{equation}
	Second-order discretizations, both in stencil width and accuracy, introduce additional unknowns, despite having information of all odd derivatives at either boundary. Without even derivatives, one can show these discretizations are linearly dependent through the method of undetermined coefficients and the Casoratian \cite{casoratian}, the discrete analogue of the Wronskian. Taking the time integrals of the $\mathcal{O}(h)$ discretizations above, we have four equations to remove four unknowns $f_{-1}(W,T)$, $f_0(W,T)$, $g_{0}(W,T)$, and $g_{1}(W,T)$:
	\begin{equation}\begin{dcases}
		e^{WT} \hat{q}(k,T) - \hat{q}(k,0) - \left[ \frac{f_{-1} - e^{ikh} f_{0} + e^{-i k L}\left( g_1 - e^{-ikh} g_0\right)}{h} \right] = 0, \\
		e^{WT} \hat{q}(-k,T) - \hat{q}(-k,0) - \left[ \frac{f_{-1} - e^{-ikh} f_{0} + e^{i k L}\left( g_1 - e^{ikh} g_0\right)}{h} \right] = 0,\\
		\frac{f_0 - f_{-1}}{h} = U^{(1)}, \\
		\frac{g_1 - g_{0}}{h} = V^{(1)},
		\label{heat_N_system}
	\end{dcases}\end{equation}
	with
	$$U^{(1)}(W,T) = \int_{0}^T e^{W t} u^{(1)}(t) \, dt, \quad\quad \text{ and } \quad\quad V^{(1)}(W,T) = \int_{0}^T e^{W t} v^{(1)}(t) \, dt.$$
	Solving \eqref{heat_N_system},
	\begin{equation}\hspace*{-30pt}\begin{dcases}
		\frac{f_{-1} - e^{ikh} f_{0}}{h} = \frac{1}{e^{2 i k (L +h)} - 1} \bigg[e^{2i k (L + h)}e^{W T} \hat{q}(k,T) - e^{2i k ( L +h)} \hat{q}(k,0) + e^{i k h} e^ {WT} \hat{q}(-k,T) \\
		\hspace{160pt} - e^{i k h} \hat{q}(-k,0) + \left(1+e^{i k h}\right) U^{(1)} -\left(e^{i k (L +h)}+e^{i k (L +2 h)}\right) V^{(1)}\bigg], \\
		\frac{e^{-ikL}\left(g_1 - e^{-ikh} g_0\right)}{h} = \frac{1}{e^{2 i k (L +h)} - 1} \bigg[ -e^{ikh} e^{WT} \hat{q}(-k,T) + e^{i k h} \hat{q}(-k,0) - e^{W T} \hat{q}(k,T) + \hat{q}(k,0) \\ 
		\hspace{245pt} - \left(1+e^{i k h}\right) U^{(1)}+ \left(e^{i k (L +h)} + e^{i k (L +2 h)}\right) V^{(1)} \bigg],
		\label{heat_N_system_solved}
		\end{dcases}\end{equation}
	with (removable) singularities at $k_{\ell} = \pi \ell/(L+h)$ for $-(N+2) \leq \ell \leq N+2$. Since the boundary integrals of \eqref{soln2_heat_centered_N} are off the real line, substituting \eqref{heat_N_system_solved} into \eqref{soln2_heat_centered_N} gives
	\begin{align}
	\begin{split}
		\hspace*{-25pt}q_n(T) &= \frac{1}{2 \pi} \int_{-\pi/h}^{\pi/h} e^{iknh} e^{-WT}\hat{q}(k,0)\,dk \\
		&\quad\, - \frac{1}{2 \pi} \int_{\partial \tilde{V}^+} e^{i k nh} e^{-WT} \left[  \frac{e^{2 i k (L +h)} \hat{q}(k,0) + e^{i k h}\hat{q}(-k,0) }{e^{2 i k (L +h)} - 1} - \frac{ \left(1+e^{i k h}\right)\left(U^{(1)} -e^{i k (L +h)} V^{(1)}\right)}{e^{2 i k (L +h)} - 1} \right] \,dk \\
		&\quad\, - \frac{1}{2 \pi} \int_{\partial \tilde{V}^-} e^{ik nh} e^{-WT} \left[ \frac{e^{i k h} \hat{q}(-k,0) + \hat{q}(k,0)}{e^{2 i k (L +h)} - 1} -  \frac{ \left(1+e^{i k h}\right)\left( U^{(1)}- e^{i k (L +h)} V^{(1)}\right)}{e^{2 i k (L +h)} - 1} \right]\,dk.
		\label{soln_heat_centered_N}
	\end{split}
	\end{align}
	after checking that the integrals with $\hat{q}(k,T)$ and $\hat{q}(-k,T)$ vanish by deforming the paths to $\tilde{D}^{\pm}$ and taking $R \rightarrow \infty$.
	
	The solution representation for IBVP \eqref{heat_prob_N} using the continuous UTM \cite{bernard_fokas} is
		\begin{align}
	\begin{split}
		q(x,T) &= \frac{1}{2 \pi} \int_{-\pi/h}^{\pi/h} e^{ikx} e^{-\tilde{W}T}\hat{q}(k,0)\,dk \\
		&\quad\, - \frac{1}{2 \pi} \int_{\partial \tilde{\Omega}^+} e^{i k x} e^{-\tilde{W}T} \left[ \frac{e^{2i k L} \hat{q}(k,0) + \hat{q}(-k,0) }{e^{2 i k L} - 1} - \frac{2\left( F_1 - e^{i k L} G_1\right)}{e^{2 i k L} - 1} \right] \,dk \\
		&\quad\, - \frac{1}{2 \pi} \int_{\partial \tilde{\Omega}^-} e^{ik x} e^{-\tilde{W}T} \left[ \frac{\hat{q}(-k,0) + \hat{q}(k,0)}{e^{2 i k L} - 1} - \frac{ 2\left(F_1- e^{i k L} G_1\right)}{e^{2 i k L} - 1}\right]\,dk,
		\label{soln_heat_cont_N}
	\end{split}
	\end{align}
	where $\tilde{\Omega}^{\pm}$ is illustrated in Figure \ref{new_paths_heat_cont_FI}. Referencing \eqref{soln_heat_centered_N}, the continuum limits of the coefficients of $\hat{q}(\pm k,0)$, $U^{(1)}$, and $V^{(1)}$ converge to their continuous counterparts, where 
	\begin{align*}
		\lim_{h\rightarrow 0} U^{(1)} &= \lim_{h \rightarrow 0} \int_0^T e^{Wt} u^{(1)}(t)\,dt = \int_0^T e^{\tilde{W}t} u^{(1)}(t)\,dt = F_1,\\
		\lim_{h\rightarrow 0} V^{(1)} &= \lim_{h \rightarrow 0} \int_0^T e^{Wt} v^{(1)}(t)\,dt = \int_0^T e^{\tilde{W}t} v^{(1)}(t)\,dt = G_1.
	\end{align*}
	

	\subsubsection{\textbf{Series Representation}} 
	Following similar steps to those in Section \ref{heat_centered_series}, we deform $\partial \tilde{V}^{\pm}$ to $\partial \tilde{V}_{\epsilon}^{\pm}$, so that the integral representation of \eqref{soln_heat_cont_N} is rewritten in terms of residue contributions:
	\begin{align}
	\begin{split}
		q_n(T) &= i \left( \sum_{\ell = -N - 1}^{N+1}\underset{k = k_{\ell}}{\text{Res}} A(n,T,k) - \sum_{\ell = -N - 1}^{N+1}\underset{k = k_{\ell}}{\text{Res}} B(n,T,k) \right),
		\label{FI_Neumann_res}
	\end{split}
	\end{align}
	with
	$$\hspace{-25pt}A(n,T,k) = e^{iknh}  e^{-WT}\left[  \frac{e^{i k h} \hat{q}(-k,0) + \hat{q}(k,0)}{e^{2 i k (L +h)} - 1} \right] \quad \text{ and } \quad B(n,T,k) = e^{iknh} e^{-WT} \left[\frac{ \left(1+e^{i k h}\right)\left( U^{(1)}- e^{i k (L +h)} V^{(1)}\right)}{e^{2 i k (L +h)} - 1}  \right]. $$
	There are no contributions from the residues at the endpoints $k_{\pm(N+2)} = \pm \pi/h$ as in Section \ref{heat_centered_series}, but we find a nonzero contribution at $k_{0} = 0$. Specifically,
	\[\begin{dcases*}
		\sum_{\ell = -N - 1}^{N+1}\underset{k = k_{\ell}}{\text{Res}} \left\{A(n,T,k)\right\} = \frac{L }{ i \left(L +h\right)} \left[  \frac{b_0}{2} + \sum_{\ell = 1}^{N+1} e^{-W_{\ell}T} b_{\ell} \cos\left(\frac{\pi \ell \left(n + \tfrac{1}{2}\right)h}{L + h}\right) \right],\\[2.5pt]
		\sum_{\ell = -N - 1}^{N+1}\underset{k = k_{\ell}}{\text{Res}} \left\{B(n,T,k)\right\} = \frac{2}{i (L +h)} \left[ \frac{H(W_{0})}{2} + \sum_{\ell = 1}^{N+1} e^{-W_{\ell}T} \cos\left(\frac{\pi \ell h}{2\left(L + h\right)}\right) \cos\left( \frac{\pi \ell \left(n + \frac{1}{2}\right)h}{L + h}\right) H(W_{\ell}) \right],
	\end{dcases*}\]
	where $W_{\ell} \equiv W\big(k_{\ell}\big)$,
	$$b_{\ell} = \frac{2 h}{L} \sum_{m=0}^{N+1} \cos \left(\frac{\pi \ell \left(m+\tfrac{1}{2}\right) h}{L + h}\right)\phi_m, \quad \quad H(W_{\ell},T) = U^{(1)}(W_{\ell},T) + (-1)^{\ell+1} V^{(1)}(W_{\ell},T).$$
	Returning to \eqref{FI_Neumann_res}, we find
	\begin{align}\begin{split}
		q_n(T) &= \frac{L }{L +h}  \sum_{\ell = 1}^{N+1} e^{-W_{\ell}T} \cos\left(\frac{\pi \ell \left(n + \tfrac{1}{2}\right)h}{L + h}\right) \left[ b_{\ell}  - \frac{2}{L} \cos\left(\frac{\pi \ell h}{2\left(L + h\right)}\right) H(W_{\ell})\right] \\
		&\quad\,+ \frac{L b_0}{2(L +h)} - \frac{H(W_{0})}{L +h}.
		\label{class_soln_heat_centered_N}
	\end{split}
	\end{align}

\begin{figure}[tb]	
		\raggedleft
		\begin{subfigure}[t]{.45\textwidth}
			\centering
  			\includegraphics[width=1\linewidth]{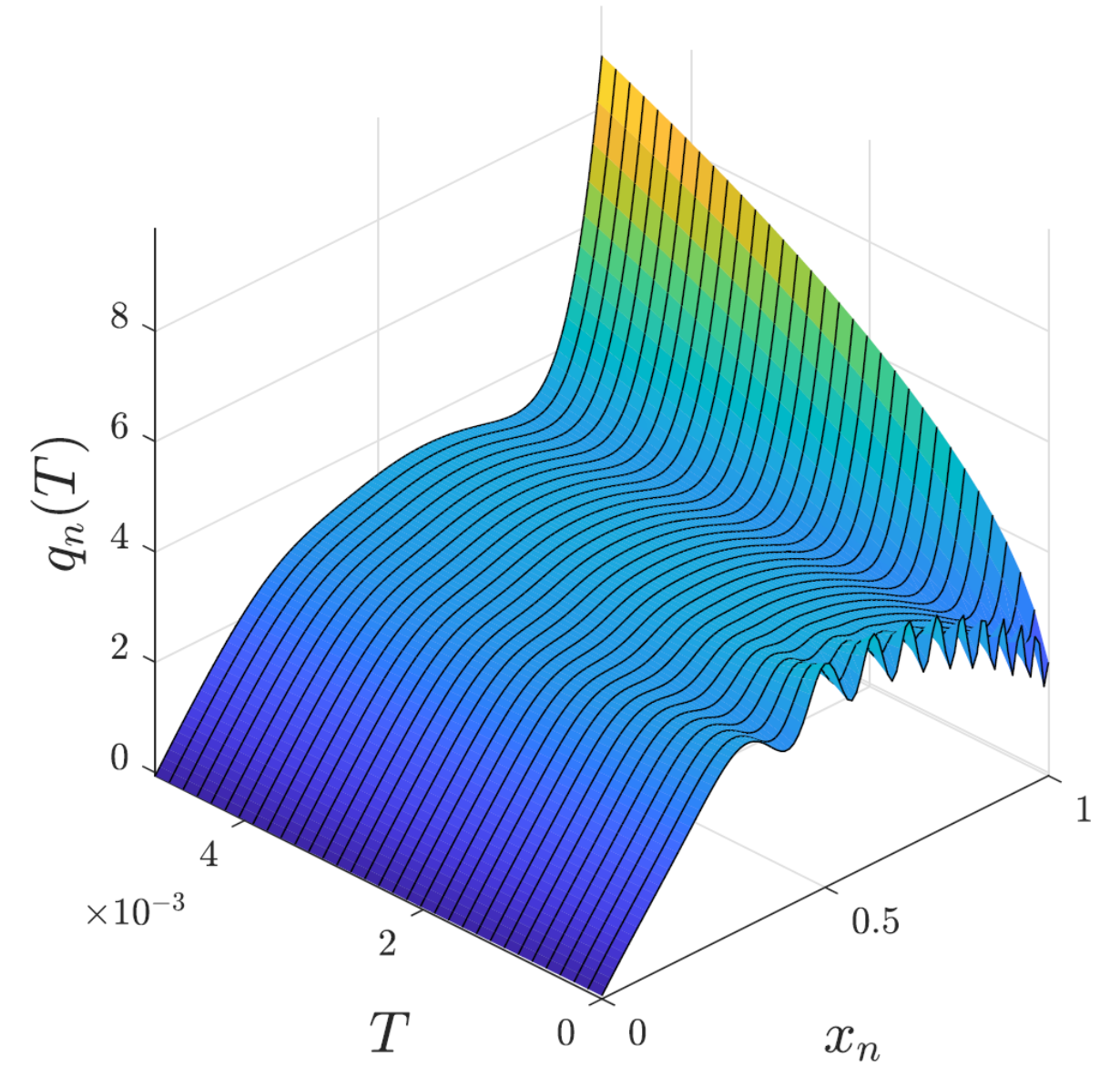}
  			\caption{}
  			\label{heat_xt_UTM2_FI_N}
		\end{subfigure}\hfill 
		\begin{subfigure}[t]{.45\textwidth}
			\centering
  			\includegraphics[width=1\linewidth]{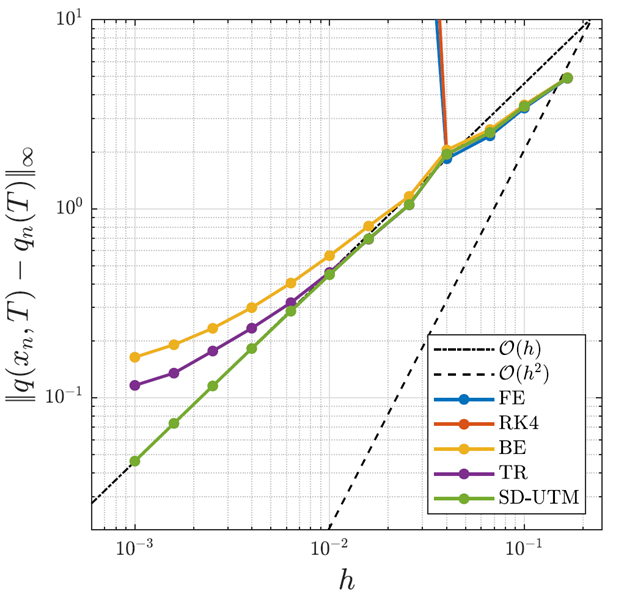}
  			\caption{}
  			\label{heat_errorplot_UTM2_FI_N}
		\end{subfigure}	
		\caption{(a) The semidiscrete solution \eqref{class_soln_heat_centered_N} evaluated at various $T$ with $h = 0.01$. (b) Error plot of the semidiscrete solution \eqref{class_soln_heat_centered_N} and finite-difference schemes relative to the exact solution as $h \rightarrow 0$ with $T = 0.005$ and $\Delta t  = 6.25 \times 10^{-4}$.}
		\label{heat_UTM2_FI_N}
	\end{figure}

	Using the series representation \eqref{class_soln_heat_centered_N} and the finite-difference schemes, we examine the solution of the IBVP
\begin{equation}
	\begin{dcases}
		q_t = q_{xx}, & 0 < x < 1,\, t > 0, \\
		q(x,0) = \phi(x) = 12x - 10x^2 + \frac{1}{2} \sin(20 \pi x^3),& 0 < x < 1,  \\
		q_x(0,t) = u^{(0)}(t) = 12,& t > 0,\\
		q_x(1,t) = v^{(0)}(t) = 30\pi - 8,& t > 0.
	\end{dcases}
	\label{heat_N_numerical1_FI}
\end{equation}	
	The exact solution to \eqref{heat_N_numerical1_FI} is given by 
	$$q(x,t) = \left(15\pi -10\right) x^2 + 12 x + \left(30\pi -20\right) t +  a_0 + \sum_{n = 1}^{\infty} a_n e^{-\left(n \pi\right)^2 t} \cos\left(n \pi x \right), $$
	with
	$$a_0 = \int_{0}^1 \left[ \frac{\sin(20 \pi x^3)}{2}  - 15\pi x^2 \right]\,dx \quad\quad \text{ and }\quad\quad a_n = 2 \int_{0}^1 \left[\frac{\sin(20 \pi x^3)}{2} - 15\pi x^2\right] \cos\left(n \pi x\right)\,dx.$$
	The given Neumann data is discretized using the first-order accurate stencils \eqref{heat_N_stencils}, which reduce the overall accuracy of the solution from the expected $\mathcal{O}(h^2)$ to $\mathcal{O}(h)$. With the centered stencil \eqref{heat_centered}, solutions \eqref{soln_heat_centered_N} and \eqref{class_soln_heat_centered_N} are fourth-order accurate approximations to the dissipative PDE \eqref{heat_centered_modified_eqn}. However, the modified equations from the Neumann boundary conditions are 
	\begin{align*}
		q_{x}(0,t) = u^{(1)}(t) - \left(\frac{h}{2}\right)q_{xx}(0,t), \quad\quad \text{ and } \quad\quad q_{x}(L,t) = v^{(1)}(t) + \left(\frac{h}{2}\right)q_{xx}(L,t),
	\end{align*} 
	implying the loss of accuracy is through the form of dissipation near the boundaries. Even so, the solution profiles from Figure \ref{heat_xt_UTM2_FI_N} depict the general diffusive behavior of the heat equation, quickly damping the high-frequency oscillations. Evaluating $a_0$ and $a_n$ numerically, we obtain the error plot in Figure \ref{heat_errorplot_UTM2_FI_N}. For consistency, the finite-difference schemes there also incorporate the Neumann boundary conditions using the first-order stencils \eqref{heat_N_stencils}.
	


	\subsection{Higher-Order Discretization of $\bms{q_t = q_{xx}}$ with Dirichlet boundary conditions}
\begin{figure}[tb]
		\begin{center}
			\def\svgwidth{2.75in}
			\vspace{10pt}\includegraphics[width=0.45\linewidth]{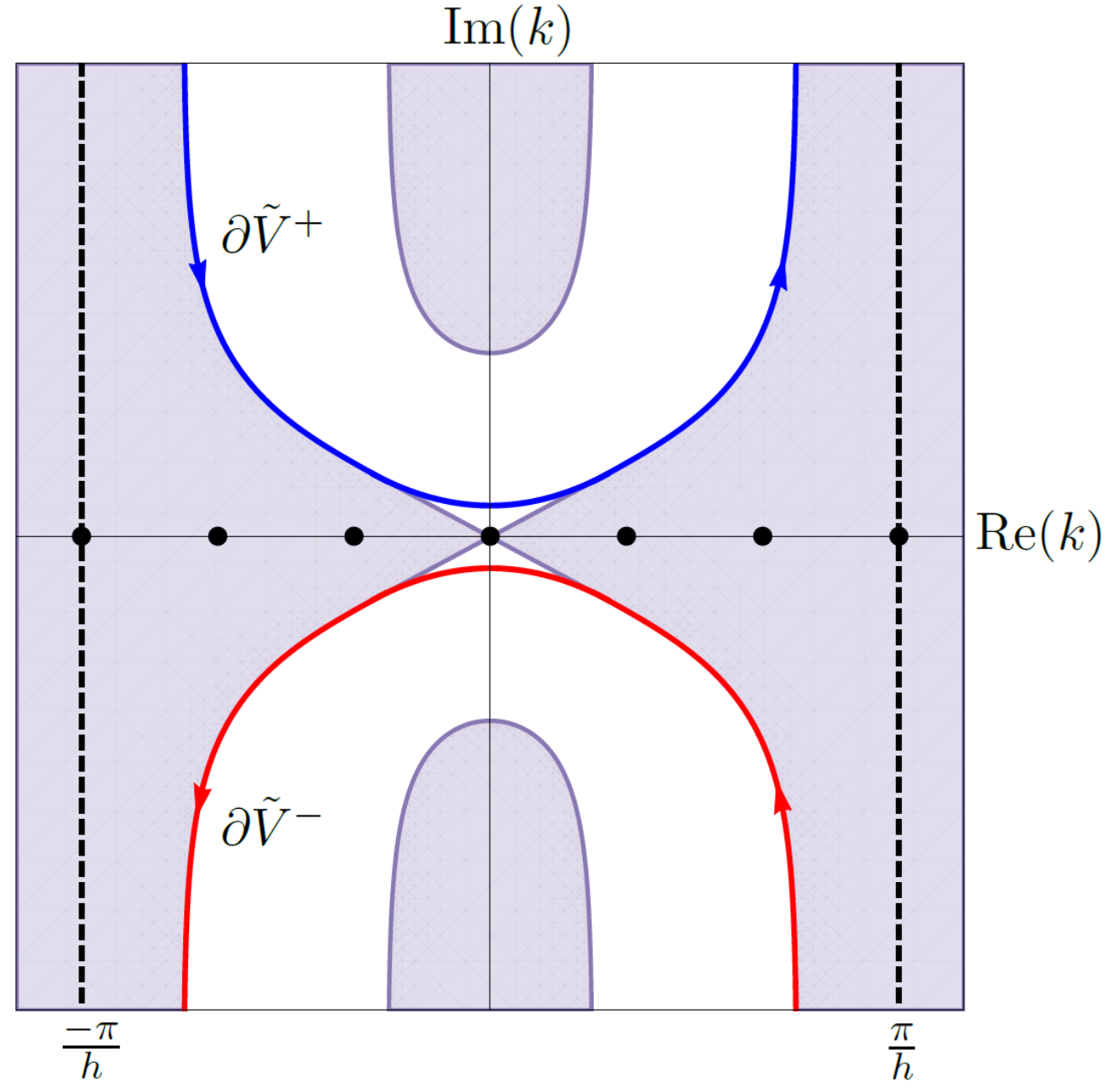}\vspace{10pt}
			\caption{The shaded regions depict where $\text{Re}(-W) \leq 0$ and $e^{-WT}$ is bounded, for the dispersion relation \eqref{W_heat_centered4}. The integration paths for $\partial \tilde{V}^{\pm}$, deformed away from the singularities on the real line, are also shown.}
			\label{paths_heat_centered4_FI}
		\end{center}
	\end{figure}

	As in Section \ref{advec_forward_finiteinterval_2}, we can apply higher-order centered discretizations to the heat equation. As before, we need extra equations in addition to the global relation formulas with $ k \rightarrow \nu_j(k)$ to eliminate unknowns. As an example, consider the Dirichlet problem \eqref{heat_prob} with the standard centered fourth-order discretization:
	\begin{equation}
		\dot{q}_n(t) =\frac{-q_{n-2}(t) + 16 q_{n-1}(t) - 30q_{n}(t) + 16q_{n+1}(t) - q_{n+2}(t)}{12h^2}.
		\label{heat_centered4}
	\end{equation}
	After some tedious steps, the global relation is
		\begin{align}
		e^{WT} \hat{q}(k,T) - \hat{q}(k,0) - \left[\frac{f(k,T) + e^{ - i k L}g(k,T)}{12h} \right]&= 0, \quad k \in \mathbb{C},
		\label{GR_heat_centered4}
	\end{align}
	where
	\[\begin{dcases*}
		f(k,T) = -e^{-i k h}f_{-1} + 16 e^{-i k h} f_0 - e^{-2 i k h} f_0 + e^{i k h} f_1 - 16 f_1 + f_2,\\
		g(k,T) = e^{-i k h} g_{-1} + 16 e^{i k h} g_0 - e^{2 i k h} g_0 - e^{i k h} g_1 + g_{-2} - 16 g_{-1},
	\end{dcases*}\]	
	with dispersion relation
	\begin{equation}
		W(k) = \frac{e^{-2ikh} - 16 e^{-ikh} + 30 - 16e^{ikh} + e^{2ikh}}{12h^2}.
		\label{W_heat_centered4}
	\end{equation}
	Solving for $\hat{q}(k,T)$ and taking the inverse transform, we obtain
	\begin{align}\begin{split}
		q_n(T) &= \frac{1}{2\pi} \int_{-\pi/h}^{\pi/h} e^{iknh} e^{-WT} \hat{q}(k,0)\,dk + \frac{1}{2\pi} \int_{-\pi/h}^{\pi/h} e^{iknh} e^{-WT}\left[\frac{f(k,T) + e^{-i k L}g(k,T)}{12h} \right]\,dk,
		\label{soln1_heat_centered4}
	\end{split}\end{align}
	which depends on the six unknowns $f_{-1}(W,T)$, $f_1(W,T)$, $f_2(W,T)$, $g_{-2}(W,T)$, $g_{-1}(W,T)$, and $g_1(W,T)$. The dispersion relation has the trivial symmetry $\nu_0(k) = k$ and three nontrivial symmetries:
	\begin{align*}
		\nu_1(k) &= -k, \\
		\nu_2(k)& = \frac{i}{h} \ln \left(\frac{e^{-i k h}}{2}  \left[16 e^{i k h}-e^{2 i k h} - 1 + \sqrt{\left(-16 e^{i k h}+e^{2 i k h}+1\right)^2-4 e^{2 i k h}}\right]\right), \\
		\nu_3(k) &= \frac{i}{h} \ln \left(\frac{e^{-i k h}}{2}  \left[16 e^{i k h}-e^{2 i k h} - 1 - \sqrt{\left(-16 e^{i k h}+e^{2 i k h}+1\right)^2-4 e^{2 i k h}}\right]\right),
	\end{align*} 
	where the branch cut for the square-root function is chosen to be on the positive real line. We separate $f(k,T)$ from $g(k,T)$ in the last integral of \eqref{soln1_heat_centered4} and, anticipating singularities, deform that integration path to $\partial \tilde{V}^{\pm}$, as shown in Figure \ref{paths_heat_centered4_FI}: 
	\begin{align}\begin{split}
		q_n(T) &= \frac{1}{2\pi} \int_{-\pi/h}^{\pi/h} e^{iknh} e^{-WT} \hat{q}(k,0)\,dk + \frac{1}{2\pi} \int_{\partial \tilde{V}^+} e^{iknh} e^{-WT}\left[\frac{f(k,T)}{12h} \right]\,dk \\
		&\quad\, - \frac{1}{2\pi} \int_{\partial \tilde{V}^-} e^{ik(nh-L)} e^{-WT}\left[ \frac{g(k,T)}{12h} \right]\,dk.
		\label{soln2_heat_centered4}
	\end{split}\end{align}
	
	In contrast to the half-line problem \cite{SDUTM_HL}, the global relations with $k \rightarrow \nu_{2,3}(k)$ are valid for all $k \in \mathbb{C}$. Using all symmetries, including $\nu_0(k)$, the global relations give four equations to remove four unknowns, specifically two $f_j(W,T)$ and two $g_j(W,T)$ contributions. To remove the remaining pair of $f_j(W,T)$ and $g_j(W,T)$ terms, we must introduce (at least) one more equation for each unknown. Using the given Dirichlet data, the heat equation itself gives all even-derivative boundary conditions, particularly $q_{xx}(0,t) = \dot{u}(t)$ and $q_{4x}(0,t) = \ddot{u}(t)$ on the left, and $q_{xx}(L,t) = \dot{v}(t)$ and $q_{4x}(L,t) = \ddot{v}(t)$ on the right, where $\dot{u} = d u^{(0)}(t)/dt$, $\ddot{u} = d^2 u^{(0)}(t)/dt^2$, $\dot{v} = d v^{(0)}(t)/dt$, and $\ddot{v} = d ^2v^{(0)}(t)/dt^2$. Discretizing the second-derivative conditions with the standard centered fourth-order stencils gives
	\begin{align}
		\frac{-q_{-2}(t)+16 q_{-1}(t)-30 q_0(t)+16 q_1(t)-q_{2}(t)}{12 h^2} &= \dot{u}(t), \label{qxx_left1} \\
		\frac{-q_{N-1}(t)+16 q_{N}(t)-30 q_{N+1}(t)+16 q_{N+2}(t)-q_{N+3}(t)}{12 h^2} &= \dot{v}(t) \label{qxx_right1}.
	\end{align}
	These stencils introduce an additional unknown, each requiring two more equations. We can find additional equations that relate nodal points with derivatives using the method of undetermined coefficients. Since we know $q_{4x}(0,t)$ and $q_{4x}(L,t)$, we derive the last pair of equations that does not introduce any more unknowns, maintains the same order of accuracy as \eqref{heat_centered4}, and is linearly independent of \eqref{qxx_left1} -- \eqref{qxx_right1}:
	\begin{align}
		\frac{q_{-2}(t) + 2 q_{-1}(t) - 6 q_0(t) + 2q_1(t) + q_{2}(t)}{6 h^2} &= \dot{u}(t) + \frac{h^2}{4} \ddot{u}(t), \label{qxx_left2} \\
		\frac{ q_{N-1}(t) + 2 q_{N}(t)-6 q_{N+1}(t)+2 q_{N+2}(t) + q_{N+3}(t)}{6 h^2} &= \dot{v}(t) + \frac{h^2}{4} \ddot{v}(t) \label{qxx_right2}.
	\end{align}
	Taking the time transforms of \eqref{qxx_left1} -- \eqref{qxx_right2}, we have a closed system of eight equations for all unknowns:
	\[\hspace*{-5pt}\left\{\begin{aligned}
		&0 = e^{WT} \hat{q}(k,T) - \hat{q}(k,0) - \left[\frac{f(k,T) + e^{ - i k L}g(k,T)}{12h} \right], &\quad &\frac{-f_{-2}+16 f_{-1}-30 f_0+16 f_1-f_{2}}{12 h^2} = \dot{U}, \\
		&0 = e^{WT} \hat{q}(-k,T) - \hat{q}(-k,0) - \left[\frac{f(-k,T) + e^{  i k L}g(-k,T)}{12h} \right], &\quad &\frac{-g_{-2}+16 g_{-1}-30 g_0+16 g_1-g_{2}}{12 h^2} = \dot{V}, \\
		&0 = e^{WT} \hat{q}(\nu_2,T) - \hat{q}(\nu_2,0) - \left[\frac{f(\nu_2,T) + e^{ - i \nu_2 L}g(\nu_2,T)}{12h} \right], &\quad &	\frac{f_{-2} + 2f_{-1} - 6 f_0 + 2 f_1 + f_{2}}{6 h^2} = \dot{U} + \frac{h^2}{4} \ddot{U},\\
		&0 = e^{WT} \hat{q}(\nu_3,T) - \hat{q}(\nu_3,0) - \left[\frac{f(\nu_3,T) + e^{ - i \nu_3 L}g(\nu_3,T)}{12h} \right], &\quad &\frac{g_{-2}+2 g_{-1}-6 g_0+2 g_1 + g_{2}}{6 h^2} = \dot{V} + \frac{h^2}{4} \ddot{V},
		\label{heat_centered4_system}
	\end{aligned}\right.\]
	where $\dot{U}(W,T)$ is the time transform of $\dot{u}(t)$, $\ddot{U}(W,T)$ is the time transform of $\ddot{u}(t)$, etc. Solving the system, we substitute our findings into \eqref{soln2_heat_centered4}:
	\begin{align}
		\hspace*{-30pt}q_n(T) &= \frac{1}{2\pi} \int_{-\pi/h}^{\pi/h} e^{iknh} e^{-WT} \hat{q}(k,0)\,dk \notag\\
		&\quad\, + \frac{1}{2\pi} \int_{\partial \tilde{V}^+} e^{iknh} e^{-WT}\left[\frac{\hat{q}(-k,0) - e^{2 i k L} \hat{q}(k,0) }{e^{2 i k L}-1} \right]\,dk + \frac{1}{2\pi} \int_{\partial \tilde{V}^-} e^{iknh} e^{-WT}\left[\frac{\hat{q}(-k,0) -  \hat{q}(k,0) }{e^{2 i k L}-1} \right]\,dk \notag\\
		&\quad\, + \frac{1}{2\pi}\left( \int_{\partial \tilde{V}^+} + \int_{\partial \tilde{V}^-} \right) e^{iknh} e^{-WT}\left[\frac{\left(e^{-2 i k h}-14 e^{-i k h}+14 e^{i k h}-e^{2 i k h}\right)\left(f_0 - e^{i k L} g_0\right)}{12 h \left(e^{2 i k L}-1\right)} \right]\,dk \label{soln_heat_centered4}\\
		&\quad\, - \frac{1}{2\pi} \left( \int_{\partial \tilde{V}^+} + \int_{\partial \tilde{V}^-} \right) e^{iknh} e^{-WT}\left[\frac{e^{-i k h} \left(e^{2 i k h} - 1\right)}{12 \left(e^{2 i k L} - 1\right)}\left( h \,\dot{U} + \frac{h^3}{12} \,\ddot{U}\right) \right]\,dk\notag\\
		&\quad\, + \frac{1}{2\pi} \left( \int_{\partial \tilde{V}^+} + \int_{\partial \tilde{V}^-} \right) e^{ik(nh+L)} e^{-WT}\left[\frac{e^{-i k h} \left(e^{2 i k h} - 1\right)}{12 \left(e^{2 i k L} - 1\right)}\left( h \,\dot{V} + \frac{h^3}{12} \,\ddot{V}\right) \right]\,dk, \notag
		\end{align}
	after showing integral terms with $\hat{q}(\pm k,T)$ vanish. Note that \eqref{soln_heat_centered4} has no dependence on $\nu_{2,3}(k)$ and the integration paths $\partial \tilde{V}^{\pm}$ are above/below the singularities given by $k_{\ell} = \pi \ell/h$, exactly as for the heat equation discretized to second order \eqref{heat_centered}. 
	
	In the continuum limit,  
	$$\lim_{h\rightarrow 0} \frac{e^{-2 i k h}-14 e^{-i k h}+14 e^{i k h}-e^{2 i k h}}{12 h} = 2 i k.$$
	Expanding the common factor between the derivative boundary conditions in \eqref{soln_heat_centered4},
	$$\frac{e^{-i k h} \left(e^{2 i k h} - 1\right)}{12 \left(e^{2 i k L} - 1\right)} = \frac{i k}{6 \left(e^{2 i k L}-1\right)}\,h + \mathcal{O}(h^3),$$
	so that the coefficients of $\dot{U}$ and $\dot{V}$ are $\mathcal{O}(h^2)$, while the coefficients for $\ddot{U}$ and $\ddot{V}$ are $\mathcal{O}(h^4)$. Hence, the SD-UTM solution loses dependence on the second and fourth-derivative boundary conditions in the continuum limit and the semidiscrete solution \eqref{soln_heat_centered4} converges to \eqref{soln_heat_cont_FI}. Although the expressions are tedious to derive, a series representation can be written down as before.
	 	
\begin{remark}
	To avoid the need for the additional equations \eqref{qxx_left2} -- \eqref{qxx_right2}, we can discretize the second-derivative boundary conditions with $\mathcal{O}(h^2)$ centered stencils:
	\begin{align*}
		\frac{q_{-1}(t)-2 q_0(t)+ q_1(t)}{h^2} &= \dot{u}(t), \quad\quad \frac{q_{N}(t)-2 q_{N+1}(t)+ q_{N+2}(t)}{h^2} = \dot{v}(t) .
	\end{align*}
	Together with the four global relation formulas, we solve for the original six unknowns in \eqref{soln2_heat_centered4}, which gives rise to a second-order accurate SD-UTM solution. As a consequence, there is a drop in accuracy from the intended $\mathcal{O}(h^4)$. This solution is exactly \eqref{soln_heat_centered4}, but without the inclusion of $h^3 \ddot{U}/12$ and $h^3 \ddot{V}/12$.
\end{remark}

\begin{remark}
	As in Section \ref{neumann_FI}, solving the fourth-order discretization \eqref{heat_centered4} with Neumann boundary conditions leads to a solution that is one order of accuracy less than that of the PDE stencil. After inverting the global relation for this IBVP, the ``solution'' depends on eight unknowns: $f_{i}(W,T)$ for $i = -2, \ldots, 1$ and $g_j(W,T)$ for $j = -1, \ldots, 2$. The four symmetries result in four global relation formulas to remove two $f_i(W,T)$ and two $g_j(W,T)$ terms after deforming off the real line to $\partial \tilde{V}^{\pm}$, as illustrated in Figure \ref{paths_heat_centered4_FI}. In order to not introduce more unknowns, four linearly independent discretizations to eliminate the remaining four unknowns are:
	$$\frac{f_1-f_{-1}}{2 h}=U^{(1)} + \frac{h^2}{6}\, \dot{U}, \quad\quad \frac{g_1-g_{-1}}{2 h}=V^{(1)} + \frac{h^2}{6}\, \dot{V},$$
	and
	$$\frac{f_{-2} -6 f_{-1}+3f_0+2 f_1}{6 h}= U^{(1)},  \quad\quad \frac{g_{-2} -6 g_{-1}+3g_0+2 g_1}{6 h}= V^{(1)},$$
	after taking time transforms. Here, 
	$$\dot{U}(W,T) = \int_0^T e^{Wt} \dot{u}(t)\,dt, \quad\quad \dot{V}(W,T) = \int_0^T e^{Wt}\dot{v}(t)\,dt,$$
	where $\dot{u}(t) = du^{(1)}(t)/dt$ and $\dot{v}(t) = dv^{(1)}(t)/dt$. The first pair of discretizations is $\mathcal{O}(h^4)$, but the second is $\mathcal{O}(h^3)$, where the $\mathcal{O}(h^3)$ terms depend on $q_{4x}(0,t)$ and $q_{4x}(L,t)$, respectively. Replacing this last pair of discretizations with a wider, more accurate stencil introduces more unknowns that no linearly independent discretizations can eliminate. Hence, the SD-UTM solution with these discretizations is third-order accurate. 
\end{remark}


\section{The Linear Schr\"{o}dinger Equation}

	Consider the free linear Schr\"{o}dinger (LS) equation
	\begin{equation}
		i q_t + \frac{1}{2} q_{xx} = 0 \quad\text{ or }\quad q_t = \frac{i}{2}q_{xx}.
		\label{LS}
	\end{equation}
	In contrast to the dissipative heat equation, this problem is dispersive.


\subsection{Centered Discretization of $\bms{q_t = (i/2)q_{xx}}$ with Dirichlet boundary conditions}\label{LS_D_section}
	The finite-interval IBVP with Dirichlet boundary conditions at both ends is
	\begin{equation}\begin{dcases}
		q_t = \tfrac{i}{2} q_{xx},& 0 < x < L, \, t > 0, \\
		q(x,0) = \phi(x),& 0 < x < L,\\
		q(0,t) = u^{(0)}(t),& t > 0,\\
		q(L,t) = v^{(0)}(t),& t > 0.
		\label{LS_prob}
	\end{dcases}\end{equation}
	Using the second-order centered discretization,
	\begin{equation}		
		\dot{q}_n(t) = \frac{i}{2}\left( \frac{q_{n+1}(t) - 2 q_n(t) + q_{n-1}(t)}{h^2}\right).
		\label{LS_centered_D}
	\end{equation}
	The local and dispersion relations are, respectively, 
	\begin{align}
		\partial_t \left(e^{-iknh} e^{Wt} q_n \right) &= \frac{i}{2 h^2}\Delta \left(e^{-ik(n-1)h} e^{Wt} q_{n} - e^{-iknh} e^{Wt} q_{n-1} \right),
		\label{LR_LS_centered}
	\end{align} 
	\begin{equation}
		W(k) = \frac{i}{2}\left( \frac{2 - e^{ikh} - e^{-ikh}}{h^2}\right) = \frac{i}{h^2} \left[ 1 - \cos(kh) \right].
		\label{W_LS_centered}
	\end{equation}
	Starting the forward transform at $n = 1$ and ending at $n = N$ gives the global relation
	\begin{align}
		e^{WT} \hat{q}(k,T) - \hat{q}(k,0) - \frac{i}{2}\left[ \frac{ e^{-ikh} f_{0} - f_1 + e^{-ikL} \left(e^{i kh}g_0 - g_{-1} \right)}{h} \right] &= 0, \quad k \in \mathbb{C}.
		\label{GR_LS_centered}
	\end{align}
	We obtain our ``solution'' formula by solving for $\hat{q}(k,T)$ and taking the inverse transform:
	\begin{align}\begin{split}
		q_n(T) &= \frac{1}{2 \pi} \int_{-\pi/h}^{\pi/h} e^{iknh} e^{-WT}\hat{q}(k,0)\,dk + \frac{1}{2 \pi} \int_{-\pi/h}^{\pi/h}\frac{i e^{iknh} e^{-WT}}{2}\left[ \frac{ e^{-ikh} f_{0} - f_1 }{h} \right]\,dk\\
		&\quad\, + \frac{1}{2 \pi} \int_{-\pi/h}^{\pi/h}\frac{i e^{ik(nh-L)} e^{-WT}}{2}\left[ \frac{ e^{i kh}g_0 - g_{-1} }{h} \right]\,dk.
		\label{soln1_LS_centered_FI}
	\end{split}\end{align}
	
\begin{figure}[tb]	
		\raggedleft
		\begin{subfigure}[t]{.45\textwidth}
			\begin{center}
			\def\svgwidth{2.75in}
			\vspace{10pt}\hspace{5pt}\includegraphics[width=1\linewidth]{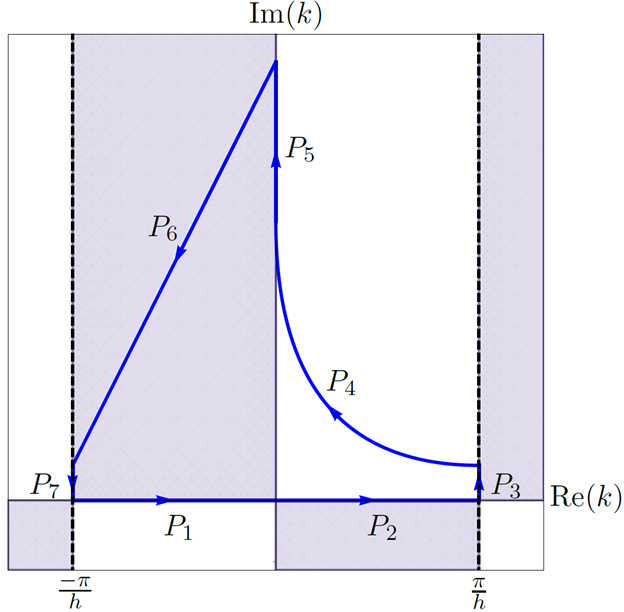}\vspace{10pt}
			\caption{}
			\label{LS_centered_W_paths_f}
		\end{center}
		\end{subfigure}\hfill 
		\begin{subfigure}[t]{.45\textwidth}
			\begin{center}
			\def\svgwidth{2.75in}
			\vspace{10pt}\includegraphics[width=1\linewidth]{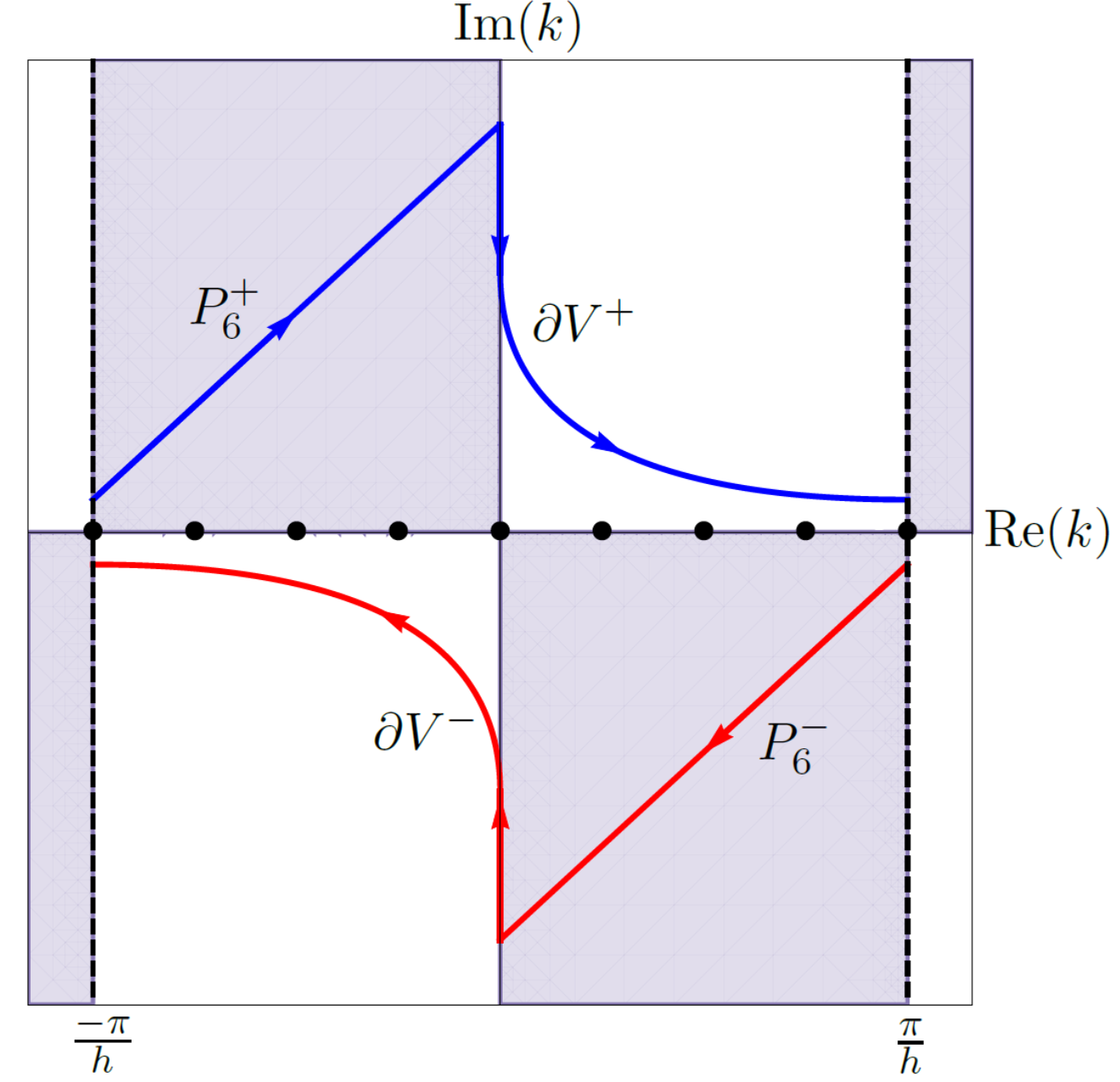}\vspace{10pt}
			\caption{}
			\label{LS_centered_paths_V}
		\end{center}
		\end{subfigure}	
		\caption{(a) The shaded regions depict where $\text{Re}(-W) \leq 0$ and $e^{-WT}$ is bounded, for the dispersion relation \eqref{W_LS_centered}. The integration paths that constitute $P$ are also shown. (b) The integration paths $\partial V^{\pm}$ and $P_6^{\pm}$.}
		\label{LS_centered_W_paths_figs}
	\end{figure}

	The dispersion relation \eqref{W_LS_centered} has the symmetries $\nu_0 (k) = k$ and $\nu_1 (k) = -k$ up to periodic copies, which can be used to remove the unknowns $f_{1}(W,T)$ and $g_{-1}(W,T)$ from \eqref{soln1_LS_centered_FI}. First, we separate the integral with the $f_j(W,T)$ terms from the integral with $g_j(W,T)$, and deform both integration paths off the real line (where we have singularities after solving for the unknowns). For the integral with the left boundary terms, consider the integration contour $P = P_1 + \ldots + P_7$, depicted in Figure \ref{LS_centered_W_paths_f}.
	Note that not all paths are straight lines. We define the integration paths' start and end points as
	\begin{align*}
		P_1 &: \text{ from } -\frac{\pi}{h} \text{ to } 0, \quad &&P_5: \text{ from } i c \text{ to } \frac{\pi}{h} i, \\ 
		P_2 &: \text{ from } 0 \text{ to } \frac{\pi}{h}, \quad &&P_6: \text{ from } \frac{\pi}{h} i \text{ to } - \frac{\pi}{h} + i \delta, \\
		P_3 &: \text{ from } \frac{\pi}{h} \text{ to } \frac{\pi}{h} + i \delta, \quad &&P_7: \text{ from } - \frac{\pi}{h} + i \delta \text{ to } - \frac{\pi}{h},\\
		P_4 &: \text{ from } \frac{\pi}{h} + i \delta \text{ to }  i c, &&
	\end{align*}
	where $\delta, c \in \mathbb{R}^+$ are nonzero constants. The point that connects $P_4$ to $P_5$ can be conveniently chosen, fixed or varying with respect to $h$, but the choice must keep the area underneath $P_4$ small, since this path is in the region of exponential growth. It is vital, however, for the curve $P_4$ to asymptotically approach the real line from above as $h \rightarrow 0$, so $\delta \ll 1$ but never zero for a finite $h$. The path $P_6$ is chosen as a straight line also for convenience, but any curved path suffices. From periodicity, $P_3 = -P_7$, so that 
	\begin{align*}
		\frac{1}{2 \pi} \int_{-\pi/h}^{\pi/h}\frac{ i e^{iknh} e^{-WT}}{2}\left( \frac{ e^{-ikh} f_{0} - f_1}{h} \right)\,dk &= \frac{1}{2 \pi} \left(\int_{-P_4} + \int_{-P_5} + \int_{-P_6} \right)\frac{ i e^{iknh} e^{-WT}}{2}\left( \frac{ e^{-ikh} f_{0} - f_1}{h} \right)\,dk\\
		&= \frac{1}{2 \pi}\int_{\partial \tilde{V}^+} \frac{ i e^{iknh} e^{-WT}}{2}\left( \frac{ e^{-ikh} f_{0} - f_1}{h} \right)\,dk,
	\end{align*}
	after defining $\partial V^+ = -P_4 -P_5$ and $\partial \tilde{V}^+ = \partial V^+ + P_6^+$ with $P_6^+ = - P_6$. Following similar arguments to deform the integral containing the right boundary information, \eqref{soln1_LS_centered_FI} becomes
	 \begin{align}
	 \begin{split}
		\hspace*{-35pt}q_n(T) &= \frac{1}{2 \pi} \int_{-\pi/h}^{\pi/h} e^{iknh} e^{-WT}\hat{q}(k,0)\,dk + \frac{1}{2 \pi} \int_{\partial \tilde{V}^+} \frac{ i e^{iknh} e^{-WT}}{2}\left( \frac{ e^{-ikh} f_{0} - f_1}{h} \right)\,dk\\
		&\quad\,- \frac{1}{2 \pi} \int_{\partial \tilde{V}^-} \frac{ i e^{ik(nh-L)} e^{-WT}}{2}\left( \frac{ e^{i kh}g_0 - g_{-1}}{h} \right)\,dk,
		\label{soln2_LS_centered_FI}
	\end{split}
	\end{align}
	where $\partial \tilde{V}^- = \partial V^- + P_6^-$, as shown in Figure \ref{LS_centered_paths_V}. Now off the real line, solving the system 
	\[\begin{dcases*}
		e^{WT} \hat{q}(k,T) - \hat{q}(k,0) - \frac{i}{2}\left[ \frac{ e^{-ikh} f_{0} - f_1 + e^{-ikL} \left( e^{i kh}g_0 - g_{-1}\right)}{h} \right] = 0, \\
		e^{WT} \hat{q}(-k,T) - \hat{q}(-k,0) - \frac{i}{2}\left[ \frac{ e^{ikh} f_{0} - f_1 + e^{ikL} \left( e^{-i kh}g_0 - g_{-1}\right)}{h} \right] = 0,
	\end{dcases*}\]
	for $f_1(W,T)$ and $g_{-1}(W,T)$ gives
	\begin{align}\begin{split}
		q_n(T) &= \frac{1}{2 \pi} \int_{-\pi/h}^{\pi/h} e^{iknh} e^{-WT}\hat{q}(k,0)\,dk \\
		&\quad\, - \frac{1}{2 \pi} \int_{\partial \tilde{V}^+} e^{iknh} e^{-WT} \left[ \frac{ e^{2 i k L} \hat{q}(k,0)  - \hat{q}(-k,0)}{e^{2 i k L} - 1} + \frac{\sin (k h) \left( f_0 -e^{i k L} g_0 \right)}{h \left(e^{2 i k L} - 1\right)}\right]\,dk\\
		&\quad\,- \frac{1}{2 \pi} \int_{\partial \tilde{V}^-} e^{iknh} e^{-WT} \left[ \frac{ \hat{q}(k,0) - \hat{q}(-k,0)}{e^{2 i k L} - 1} + \frac{ \sin (k h) \left( f_0 - e^{i k L} g_0\right)}{h\left(e^{2 i k L} - 1\right)} \right]\,dk,
		\label{soln_LS_centered_FI}
	\end{split}\end{align}
	after removing the integral terms with $\hat{q}(\pm k,T)$. The simple poles at $k_{\ell} = \pi \ell/L$, the same as in \eqref{soln_heat_centered_FI} to the second-order discretized heat equation, do not interfere with the integration paths $\partial \tilde{V}^{\pm}$. When numerically evaluating the integrals above, the paths $P_4$ in $\partial \tilde{V}^{\pm}$ are parameterized as exponentially-decaying curves toward the real line as $\text{Re}(k) \rightarrow \infty$.
	
\begin{figure}[tb]
		\begin{center}
			\def\svgwidth{2.75in}
			\vspace{10pt}\includegraphics[width=0.45\linewidth]{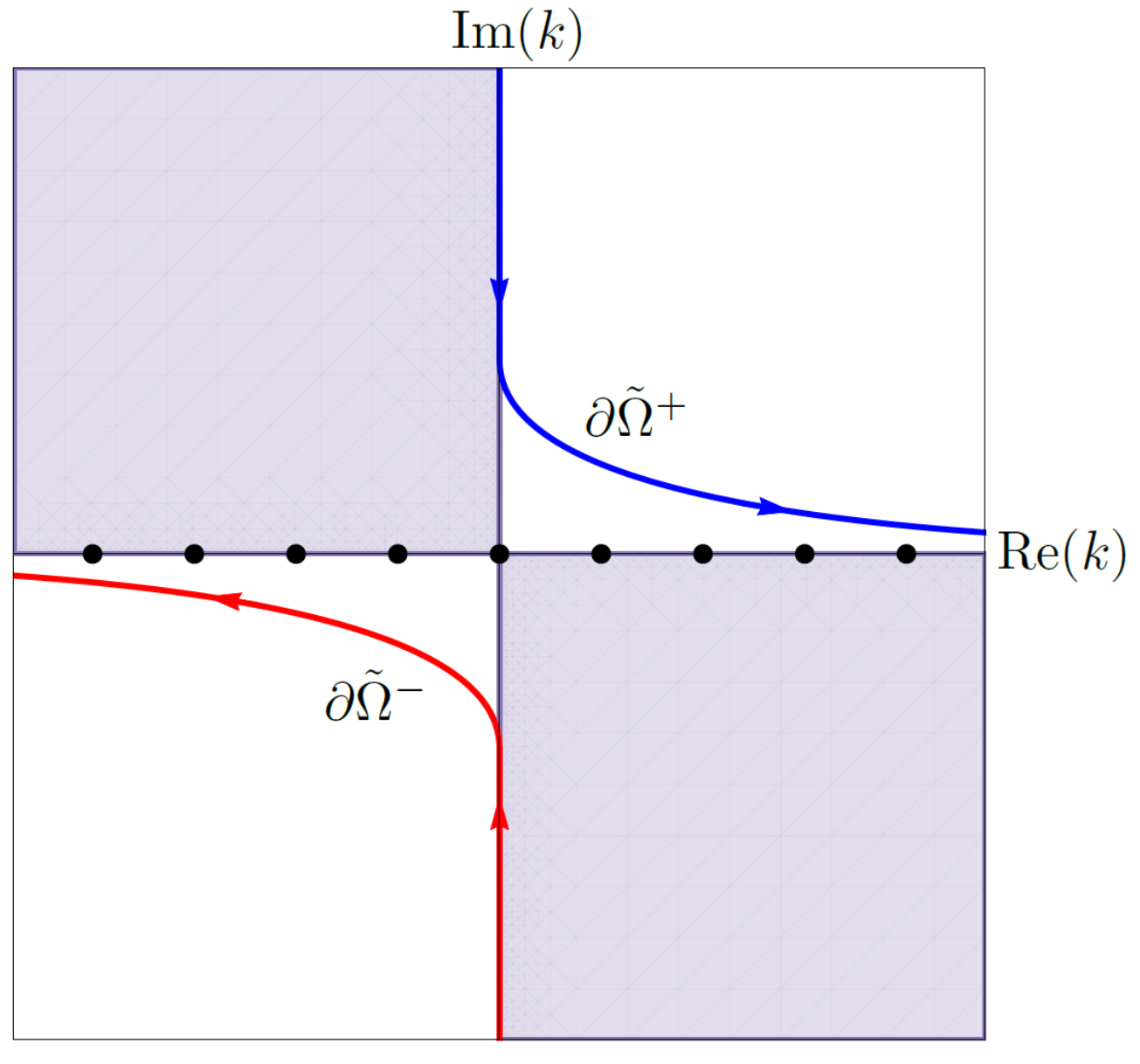}\vspace{10pt}
			\caption{The shaded regions depict where $\text{Re}(-\tilde{W}) \leq 0$ and $e^{\tilde{W}T}$ is bounded. The integration paths $\partial \tilde{\Omega}^{\pm}$ are also shown.}
			\label{Omega_path_LS_FI_new}
		\end{center}
	\end{figure}

	With the integration paths $\tilde{\Omega}^{\pm}$ illustrated in Figure \ref{Omega_path_LS_FI_new} and the dispersion relation $\tilde{W}(k) = i k^2/2$, the continuous UTM solution \cite{bernard_fokas} to \eqref{LS_prob} is	
	\begin{align}
	\begin{split}
		q(x,T) &= \frac{1}{2 \pi} \int_{-\infty}^{\infty} e^{ikx} e^{-\tilde{W}T} \hat{q}(k,0)\, dk \\
		&\quad\, - \frac{1}{2 \pi} \int_{\partial \tilde{\Omega}^+} e^{ikx} e^{-\tilde{W}T} \left[\frac{e^{2 i k L} \hat{q}(k,0) - \hat{q}(-k,0) + k \left(F_0  - e^{i k L} G_0\right) }{e^{2 i k L} - 1}\right]\,dk \\
		&\quad\, - \frac{1}{2 \pi} \int_{\partial \tilde{\Omega}^-} e^{ikx} e^{-\tilde{W}T} \left[\frac{\hat{q}(k,0) - \hat{q}(-k,0) + k \left(F_0 - e^{i k L} G_0 \right)}{e^{2 i k L} - 1} \right]\,dk.
		\label{soln_LS_cont_FI}
	\end{split}
	\end{align}
	
	Taking the continuum limit, we find $\lim_{h \rightarrow 0} W(k) = \tilde{W}(k)$ and $\lim_{h \rightarrow 0} \partial V^{\pm} = \tilde{\Omega}^{\pm}$, where $\partial \tilde{V}^{\pm} = \partial V^{\pm} + P_6^{\pm}$. For $P_6^+$ in the upper-half plane, we have
	\begin{align}
		 \hspace*{-15pt}\lim_{h \rightarrow 0} \,\frac{1}{2 \pi} \int_{-\frac{\pi}{h} + i \delta}^{\frac{\pi}{h} i} e^{iknh} e^{-WT}& \left[ \frac{ e^{2 i k L} \hat{q}(k,0) - \hat{q}(-k,0) }{e^{2 i k L} - 1} + \frac{\sin (k h) \left( f_0 -e^{i k L} g_0\right)}{h \left(e^{2 i k L} - 1\right)}\right]\,dk \notag\\
		&= \frac{1}{2 \pi} \int_{-\infty + i \delta}^{i\infty} e^{ikx} e^{-\tilde{W}T} \left[\frac{e^{2 i k L} \hat{q}(k,0) - \hat{q}(-k,0) + k \left(F_0  - e^{i k L} G_0\right) }{e^{2 i k L} - 1}\right]\,dk.
		\label{P6plus_limit}
	\end{align}
	With $R \gg 1$,
	$$\hspace{-15pt}\frac{1}{2 \pi} \int_{-R + i \delta}^{i R} e^{ikx} e^{-\tilde{W}T} \left[\frac{e^{2 i k L} \hat{q}(k,0) - \hat{q}(-k,0) }{e^{2 i k L} - 1}\right]\,dk = \frac{1}{2 \pi} \int_{0}^L \left[ \int_{-R + i \delta}^{i R} e^{ikx} e^{-\tilde{W}T} \left(\frac{e^{2 i k L} e^{-iky} - e^{iky} }{e^{2 i k L} - 1}\right)\,dk \right] q(y,0)\,dy.$$
	Because the integration path is in the shaded region of Figure \ref{Omega_path_LS_FI_new}, taking $R \rightarrow \infty$ implies the integrand 
	$$ e^{ikx} e^{-\tilde{W}T} \left[\frac{e^{i k (2 L - y)} - e^{iky} }{e^{2 i k L} - 1}\right] \rightarrow 0,$$
	since every exponent is positive and decays in the integration region. Similarly for \eqref{P6plus_limit}. Hence,
	$$\frac{1}{2 \pi} \int_{-\infty + i \delta}^{i \infty} e^{ikx} e^{-\tilde{W}T} \left[\frac{e^{2 i k L} \hat{q}(k,0) - \hat{q}(-k,0) + k \left(F_0  - e^{i k L} G_0\right) }{e^{2 i k L} - 1}\right]\,dk = 0,$$
	and in the continuum limit,
	$$\lim_{h \rightarrow 0}\, \frac{1}{2 \pi} \int_{P_6^+} e^{iknh} e^{-WT} \left[ \frac{ \hat{q}(-k,0) - e^{2 i k L} \hat{q}(k,0) }{e^{2 i k L} - 1} + \frac{\sin (k h) \left( e^{i k L} g_0 - f_0 \right)}{h \left(e^{2 i k L} - 1\right)}\right]\,dk  = 0.$$
	We reach a similar conclusion for the $P_6^-$ integral. Therefore, the SD-UTM solution \eqref{soln_LS_centered_FI} converges to its continuous counterpart \eqref{soln_LS_cont_FI}.
	

	\subsubsection{\textbf{Series Representation}} 
	Although this is a dispersive problem, the series representation of \eqref{soln_LS_centered_FI} is obtained following almost exactly the same steps as those in Section \ref{heat_centered_series} for the heat equation with Dirichlet boundary conditions. Deforming the integration paths $\partial \tilde{V}^{\pm}$ to the singularities on the real line and determining the residue contributions gives
	\begin{align}
		q_n(T) &=  \sum_{\ell = 1}^{N}  e^{-W_{\ell}T} \sin\left(\frac{\pi \ell nh}{L}\right) \left[ b_{\ell} + \frac{i}{L h} \sin\left(\frac{\pi \ell h}{L}\right) H(W_{\ell},T) \right],
		\label{class_soln_LS_centered_FI}
	\end{align}
	where $W(k_{\ell}) \equiv W_{\ell}$ and
	$$b_{\ell} = \frac{2 h}{L} \sum_{m=0}^{N+1} \sin\left(\frac{\pi \ell mh}{L} \right)\phi_m, \quad\quad H(W_{\ell},T) =  f_0(W_{\ell},T) + (-1)^{\ell+1} g_0(W_{\ell},T).$$
	
\begin{figure}[h!]	
		\raggedleft
		\begin{subfigure}[t]{.45\textwidth}
			\centering
  			\includegraphics[width=1\linewidth]{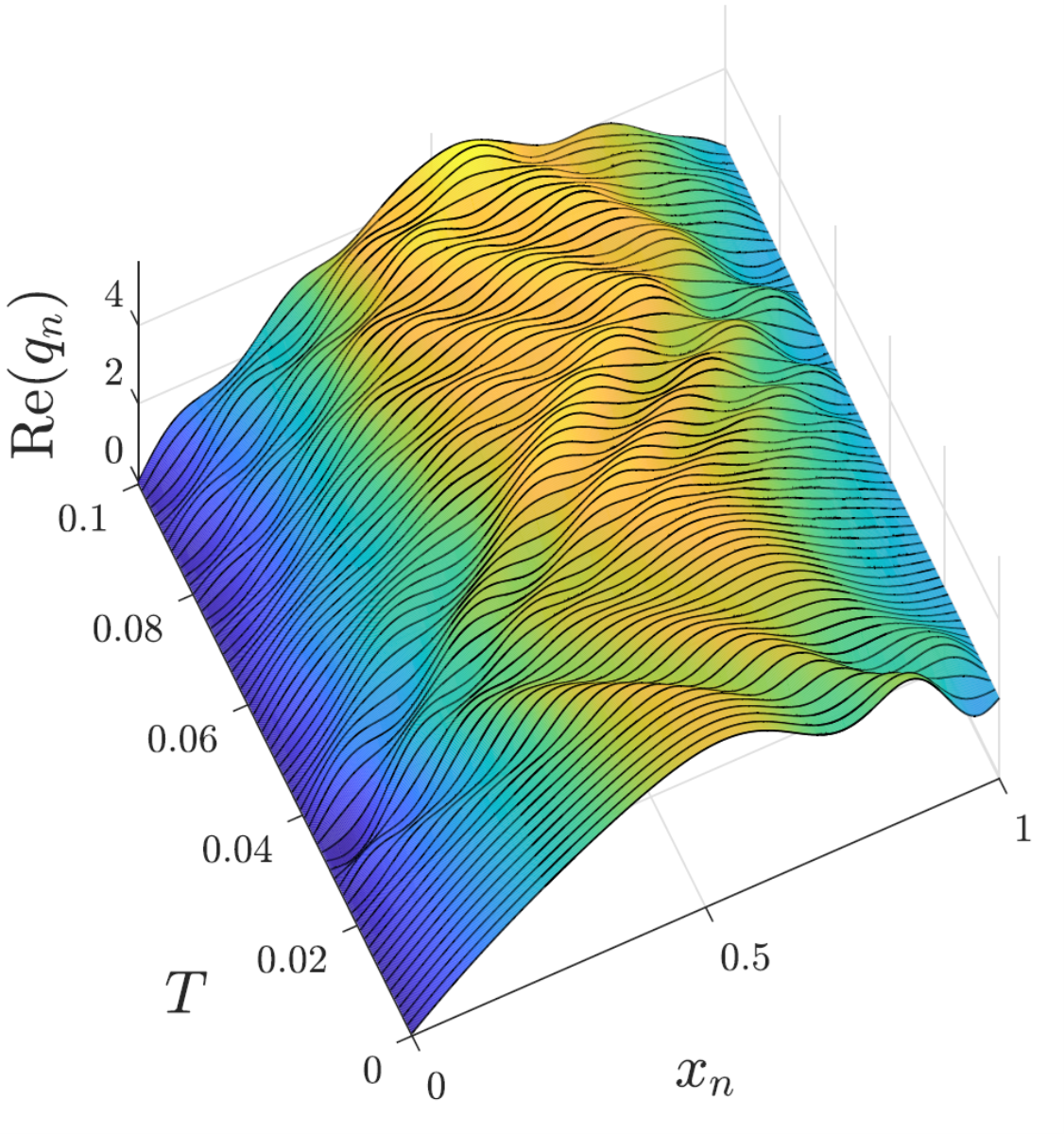}
  			\caption{}
  			\label{}
		\end{subfigure}\hfill 
		\begin{subfigure}[t]{.45\textwidth}
			\centering
  			\includegraphics[width=1\linewidth]{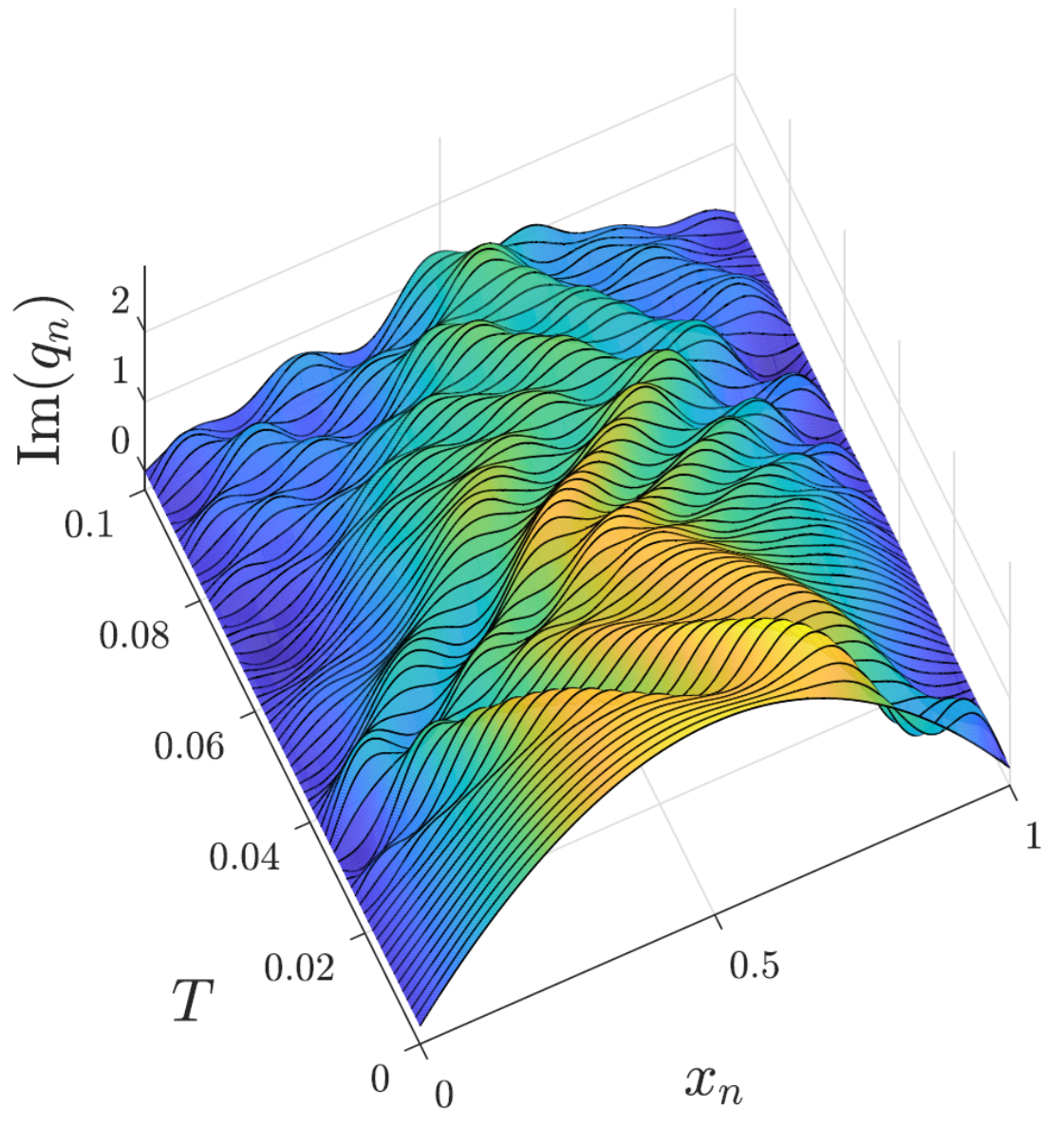}
  			\caption{}
  			\label{}
		\end{subfigure}
		\begin{subfigure}[t]{.45\textwidth}
			\centering
  			\includegraphics[width=1\linewidth]{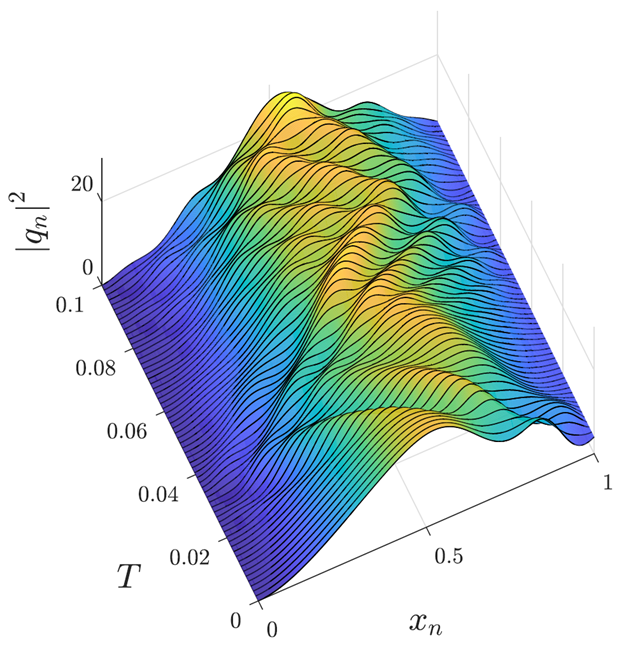}
  			\caption{}
  			\label{}
		\end{subfigure}\hfill
		\begin{subfigure}[t]{.45\textwidth}
			\centering
  			\includegraphics[width=1.1\linewidth]{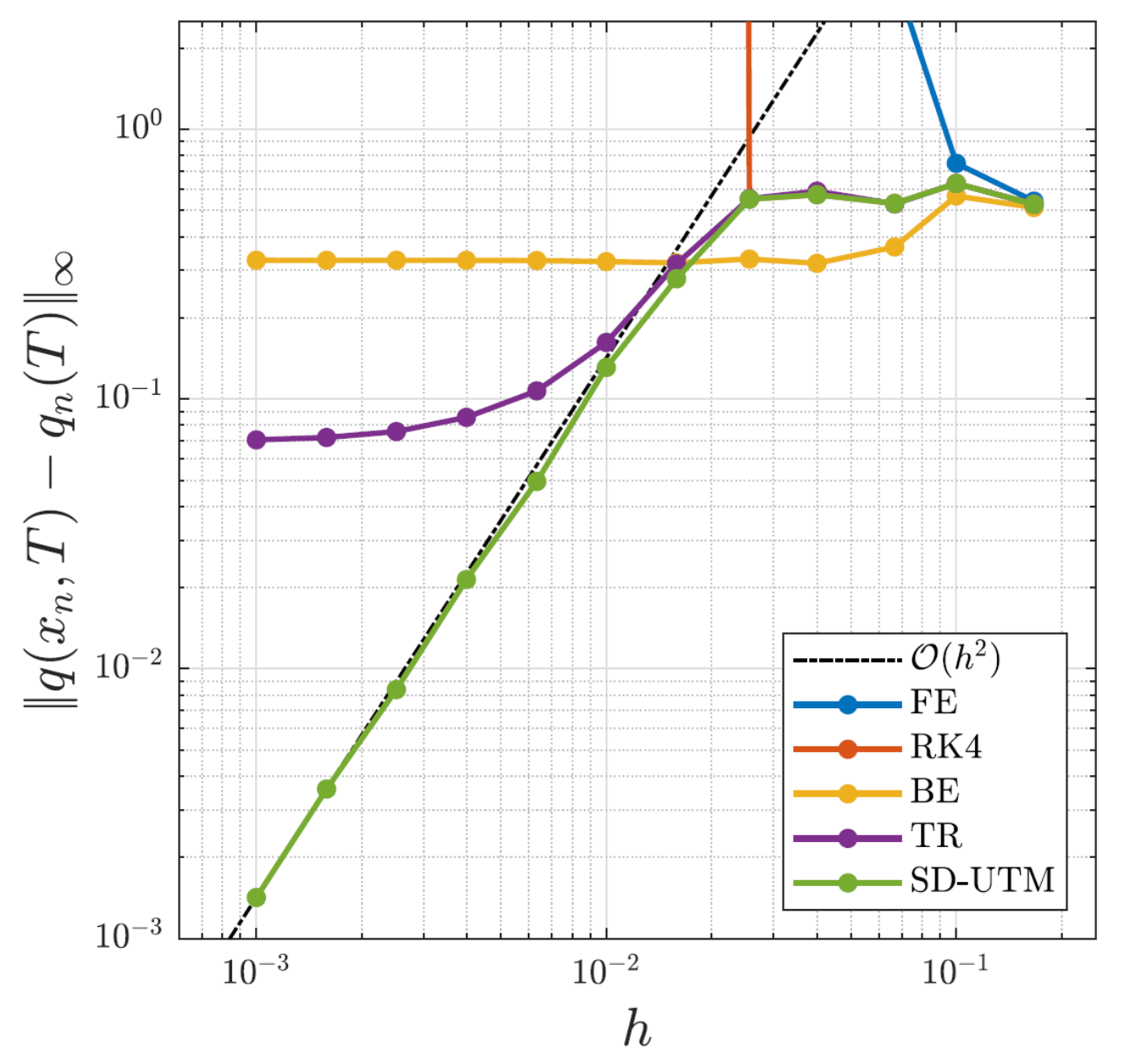}
  			\caption{}
  			\label{LS_errorplot_UTM1_FI}
		\end{subfigure}
		\caption{(a) - (c) Real and imaginary parts and modulus squared of the semidiscrete solution \eqref{class_soln_LS_centered_FI} evaluated at various $T$ for IBVP \eqref{LS_numerical1_FI} with $h = 0.01$. (d) Error plot of the semidiscrete solution \eqref{class_soln_LS_centered_FI} and finite-difference schemes relative to the exact solution as $h \rightarrow 0$ with $T = 0.1$ and $\Delta t  = 3.90625 \times 10^{-4}$.}	
		\label{LS_UTM1_FI}
	\end{figure}

	Using \eqref{class_soln_LS_centered_FI}, we examine the numerical solution to
\begin{equation}
\begin{dcases}
	q_t = \tfrac{i}{2} q_{xx},& 0 < x < 1,\, t > 0, \\ 
	q(x,0) = \phi(x) = 2(6+5i)x - 10(1+i) x^2 + \frac{1}{2}\sin\left(4 \pi x^3 \right),& 0 < x < 1,  \\
	q(0,t) = u^{(0)}(t) = 0,& t > 0, \\
	q(1,t) = v^{(0)}(t) = 2,& t > 0.
\end{dcases}
\label{LS_numerical1_FI}
\end{equation}	
	For this IBVP, the continuous solution is traditionally determined using Fourier series. The modified PDE is $p_t = (i/2)p_{xx} + (i h^2/24) p_{4x}$, which is dispersive. The higher-order dispersive term is $\mathcal{O}(h^2)$, and we can decrease any excess dispersion in numerical computations by decreasing $h$. Figure \ref{LS_UTM1_FI} shows the dispersive nature of the real and imaginary components of the SD-UTM solution, along with the square of the modulus. Once more, the SD-UTM outperforms the standard finite-difference methods, where BE's dissipative behavior practically dampens all oscillations. From the error plot in Figure \ref{LS_errorplot_UTM1_FI}, TR attempts to capture these oscillations, but not as accurately as the SD-UTM.


\subsection{Centered Discretization of $\bms{q_t = (i/2)q_{xx}}$ with Neumann boundary conditions} \label{LS_neumann_halfline}
	Lastly, we consider the second-order discretization \eqref{LS_centered_D}, but with Neumann boundary conditions at both ends of the interval:
	\begin{equation}\begin{dcases}
		q_t = \tfrac{i}{2} q_{xx},& 0 < x < L,\, t > 0, \\
		q(x,0) = \phi(x),& 0 < x < L,\\
		q_x(0,t) = u^{(1)}(t),& t > 0,\\
		q_x(L,t) = v^{(1)}(t),& t > 0.
		\label{LS_prob_N}
	\end{dcases}\end{equation}
	The local and dispersion relations, \eqref{LR_LS_centered} and \eqref{W_LS_centered} respectively, transfer over from the previous section. Now without Dirichlet data, $q_{0}(t)$ and $q_{N+1}(t)$ are unknown, and the global relation is 
	\begin{align}
		e^{WT} \hat{q}(k,T) - \hat{q}(k,0) - \frac{i}{2}\left[ \frac{f_{-1} - e^{ikh} f_{0} + e^{-i k L}\left( g_1 - e^{-ikh} g_0\right)}{h} \right] &= 0, \quad k \in \mathbb{C}.
		\label{GR_LS_centered_N}
	\end{align}
	Solving for $\hat{q}(k,T)$ and using the inverse transform, we obtain
	\begin{align}\begin{split}
		q_n(T) &= \frac{1}{2 \pi} \int_{-\pi/h}^{\pi/h} e^{iknh} e^{-WT}\hat{q}(k,0)\,dk + \frac{1}{2 \pi} \int_{-\pi/h}^{\pi/h}\frac{ i e^{iknh} e^{-WT}}{2}\left[ \frac{f_{-1} - e^{ikh} f_{0} }{h} \right]\,dk\\
		&\quad\, + \frac{1}{2 \pi} \int_{-\pi/h}^{\pi/h}\frac{ i e^{ik(nh-L)} e^{-WT}}{2}\left[ \frac{g_1 - e^{-ikh} g_0 }{h} \right]\,dk.
		\label{soln1_LS_centered_N}
	\end{split}\end{align}
	The global relations \eqref{GR_LS_centered_N} with $k$ and $k \rightarrow -k$ remove two of the four unknowns. To not introduce new unknowns, we apply the first-order backward discretization to $q_x(0,t)$ and the first-order forward discretization to $q_x(L,t)$. Upon taking time transforms, we solve the system
	\[\begin{dcases*}
		e^{WT} \hat{q}(k,T) - \hat{q}(k,0) - \frac{i}{2}\left[ \frac{ f_{-1} - e^{i k h} f_0  + e^{-ikL} \left(g_1 - e^{-i k h} g_0  \right)}{h} \right] = 0, \\
		e^{WT} \hat{q}(-k,T) - \hat{q}(-k,0) - \frac{i}{2}\left[ \frac{ f_{-1} - e^{-i k h} f_0  + e^{ikL} \left(g_1 - e^{i k h} g_0 \right)}{h} \right] = 0,\\
		\frac{f_0 - f_{-1}}{h} = U^{(1)}, \\
		\frac{g_1 - g_{0}}{h} = V^{(1)},
	\end{dcases*}\]
	 to remove all four unknowns from ``solution'' \eqref{soln1_LS_centered_N}, where
	 $$U^{(1)}(W,T) = \int_0^T e^{Wt} u^{(1)}(t)\, dt \quad \quad \text{ and }\quad\quad V^{(1)}(W,T) = \int_0^T e^{Wt} v^{(1)}(t)\, dt.$$
	Doing so and substituting into \eqref{soln1_LS_centered_N} after deforming to $\partial \tilde{V}^{\pm}$, depicted in Figure \ref{LS_centered_paths_V}, gives the first-order accurate solution
	\begin{align}
	\begin{split}
		\hspace*{-25pt}q_n(T) &= \frac{1}{2 \pi} \int_{-\pi/h}^{\pi/h} e^{iknh} e^{-WT}\hat{q}(k,0)\,dk \\
		&\quad\, - \frac{1}{2 \pi} \int_{\partial \tilde{V}^+} e^{iknh} e^{-WT}\left[ \frac{ e^{2 i k (L+h)} \hat{q}(k,0) + e^{i k h} \hat{q}(-k,0) }{e^{2 i k (L+h)} - 1} - \frac{ i \left(1+e^{i k h}\right) \left(U^{(1)} - e^{i k (L+h)} V^{(1)}\right) }{2\left(e^{2 i k (L+h)} - 1\right)} \right]\,dk\\
		&\quad\, - \frac{1}{2 \pi} \int_{\partial \tilde{V}^-} e^{iknh} e^{-WT}\left[ \frac{ \hat{q}(k,0) + e^{i k h} \hat{q}(-k,0) }{e^{2 i k (L+h)} - 1} - \frac{i \left(1+e^{i k h}\right) \left(U^{(1)} - e^{i k (L+h)} V^{(1)}\right) }{2\left(e^{2 i k (L+h)} - 1\right)} \right]\,dk,
		\label{soln_LS_centered_N}
	\end{split}
	\end{align}
	after applying similar techniques as before to remove the integral terms depending on $\hat{q}(\pm k,T)$. In the continuum limit, \eqref{soln_LS_centered_N} converges to the continuous UTM solution \cite{bernard_fokas}: 
	\begin{align}
	\begin{split}
		q(x,T) &= \frac{1}{2 \pi} \int_{-\pi/h}^{\pi/h} e^{ikx} e^{-\tilde{W}T}\hat{q}(k,0)\,dk \\
		&\quad\, - \frac{1}{2 \pi} \int_{\partial \tilde{\Omega}^+} e^{i k x} e^{-\tilde{W}T} \left[ \frac{e^{2i k L} \hat{q}(k,0) + \hat{q}(-k,0)  - i\left( F_1 - e^{i k L} G_1\right)}{e^{2 i k L} - 1} \right] \,dk \\
		&\quad\, - \frac{1}{2 \pi} \int_{\partial \tilde{\Omega}^-} e^{ik x} e^{-\tilde{W}T} \left[ \frac{\hat{q}(-k,0) + \hat{q}(k,0) - i\left(F_1- e^{i k L} G_1\right)}{e^{2 i k L} - 1}\right]\,dk,
		\label{soln_LS_cont_FI_Neumann}
	\end{split}
	\end{align}
	where $\lim_{h \rightarrow 0} U^{(1)} = F_1$ and $\lim_{h \rightarrow 0} V^{(1)} = G_1$.
	
	\subsubsection{\textbf{Series Representation}} 
	Proceeding as in the previous sections, the series representation for \eqref{soln_LS_centered_N} is
	\begin{align}
	\begin{split}
		\hspace{-25pt}q_n(T) &= \frac{L }{L +h} \sum_{\ell = 1}^{N+1} e^{-W_{\ell}T} \cos\left(\frac{\pi \ell \left(n + \tfrac{1}{2}\right)h}{L + h}\right) \left[ b_{\ell}  - \frac{i}{L} \cos\left[\frac{\pi \ell h}{2\left(L + h\right)}\right] H(W_{\ell},T)\right] + \frac{L b_0 - i H(W_{0},T)}{2(L +h)}.
		\label{class_soln_LS_centered_N}
	\end{split}
	\end{align}
	where
	$$b_{\ell} = \frac{2 h}{L} \sum_{m=0}^{N+1} \cos \left(\frac{\pi \ell \left(m+\tfrac{1}{2}\right) h}{L + h}\right)\phi_m , \quad \quad H(W_{\ell},T) = U^{(1)}(W_{\ell},T) + (-1)^{\ell+1} V^{(1)}(W_{\ell},T).$$

\begin{figure}[h!]	
		\raggedleft
		\begin{subfigure}[t]{.45\textwidth}
			\centering
  			\includegraphics[width=1\linewidth]{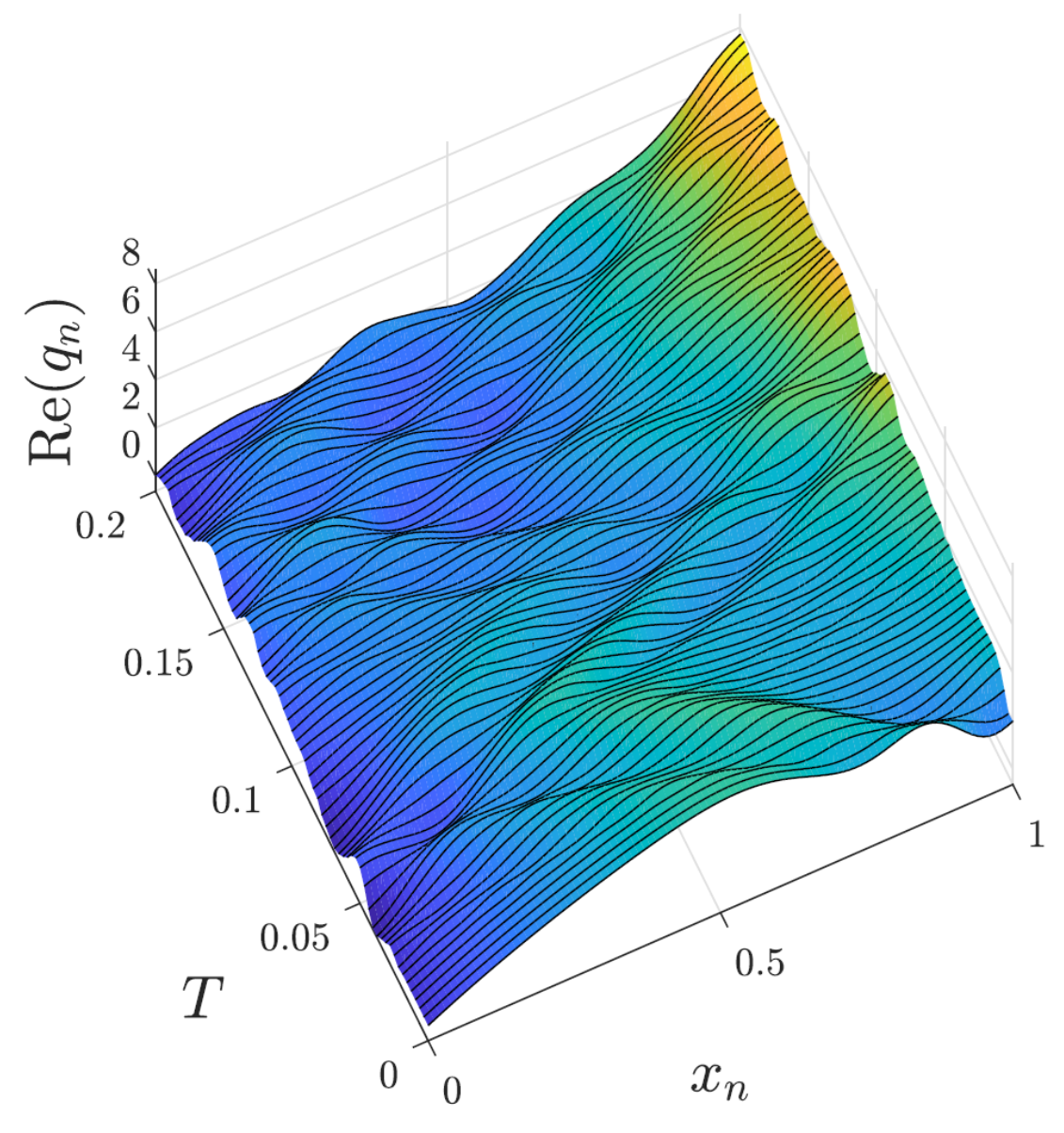}
  			\caption{}
  			\label{}
		\end{subfigure}\hfill 
		\begin{subfigure}[t]{.45\textwidth}
			\centering
  			\includegraphics[width=1\linewidth]{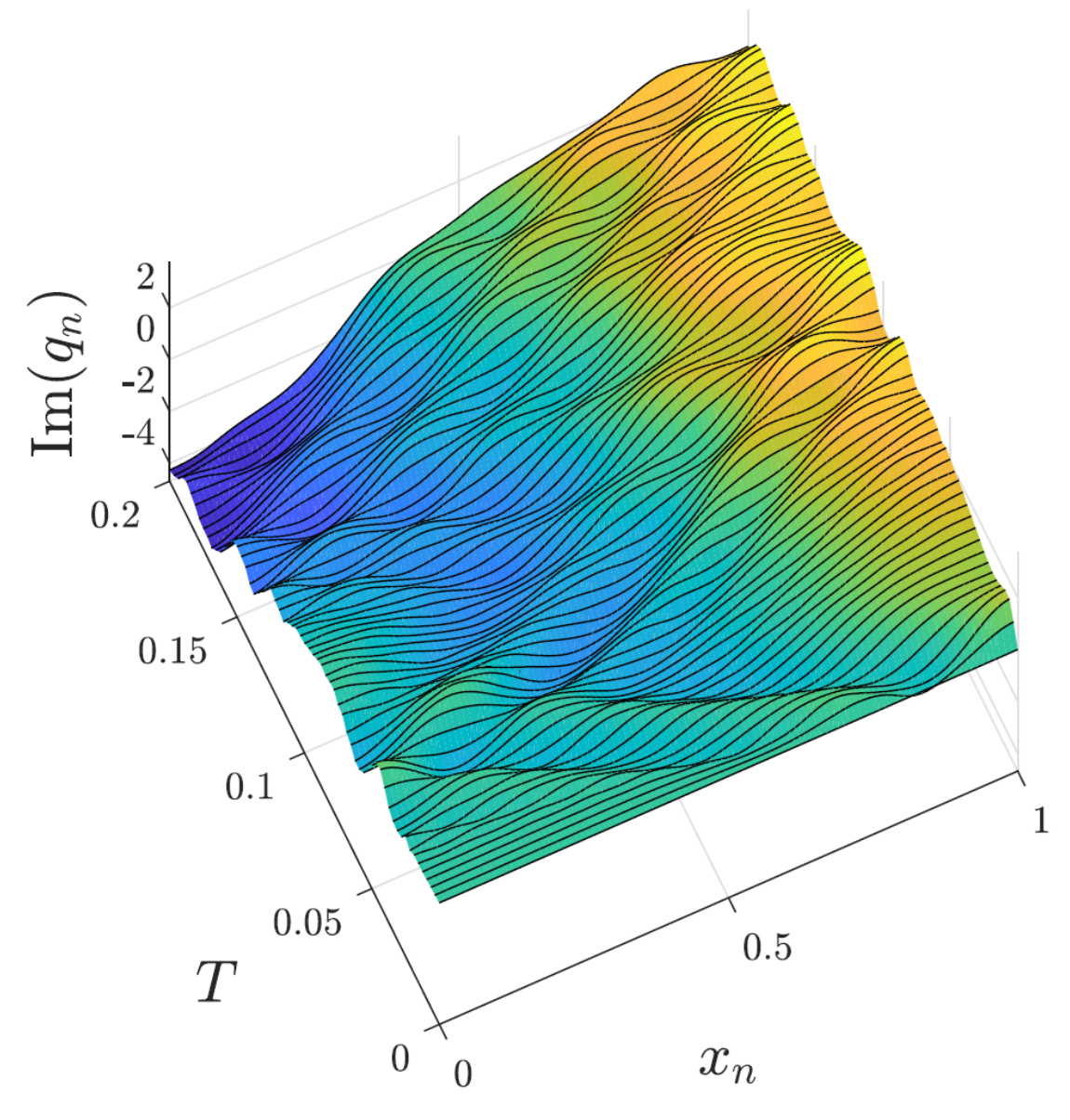}
  			\caption{}
  			\label{}
		\end{subfigure}
		\begin{subfigure}[t]{.45\textwidth}
			\centering
  			\includegraphics[width=1\linewidth]{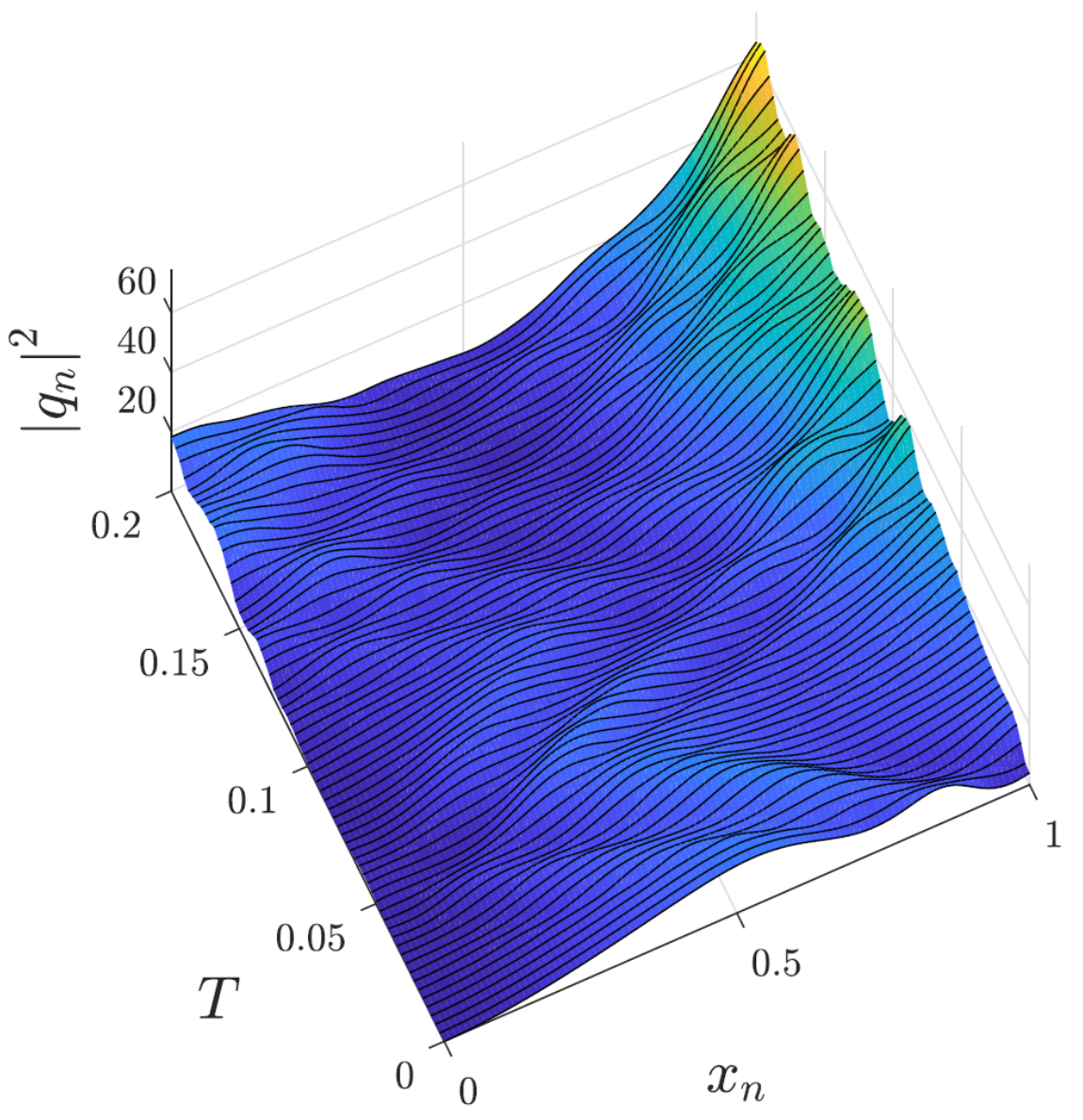}
  			\caption{}
  			\label{}
		\end{subfigure}\hfill
		\begin{subfigure}[t]{.45\textwidth}
			\centering
  			\includegraphics[width=1.1\linewidth]{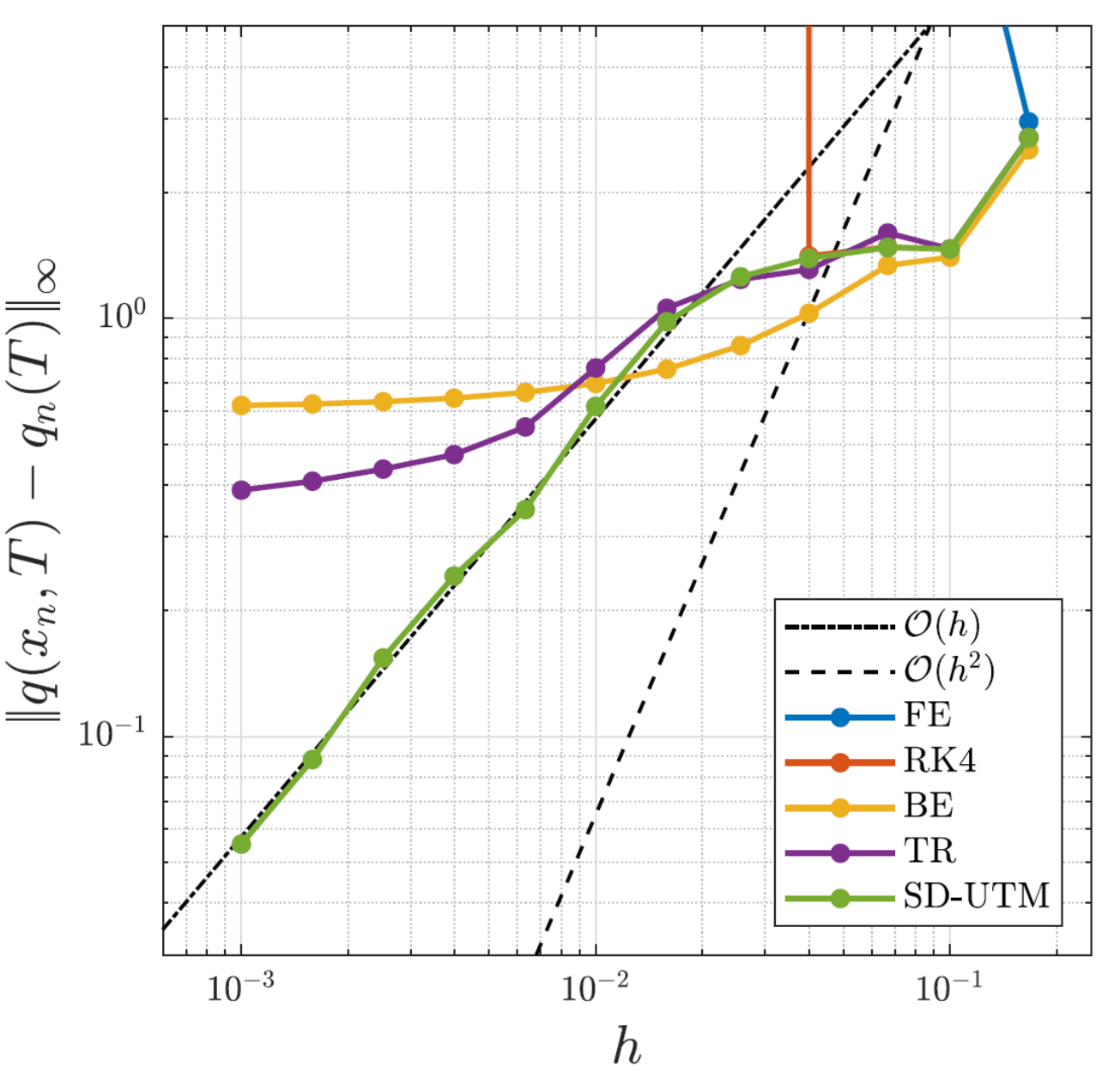}
  			\caption{}
  			\label{LS_errorplot_UTM1_FI_N}
		\end{subfigure}
		\caption{(a) - (c) Real and imaginary parts and modulus squared of the semidiscrete solution \eqref{class_soln_LS_centered_N} evaluated at various $T$ for IBVP \eqref{LS_numerical1_FI} with $h = 0.01$. (d) Error plot of the semidiscrete solution \eqref{class_soln_LS_centered_N} and finite-difference schemes relative to the exact solution as $h \rightarrow 0$ with $T = 0.2$ and $\Delta t  = 1.5625 \times 10^{-3}$.}	
		\label{LS_UTM1_FI_N}
	\end{figure}

	For a numerical example, consider the IBVP \eqref{LS_numerical1_FI}, but with Neumann boundary conditions:
	\begin{equation}
	\begin{dcases}
	q_t = \tfrac{i}{2} q_{xx},& 0 < x < 1,\, t > 0, \\
	q(x,0) = \phi(x) = 12 x - 10 x^2 + \frac{1}{2}\sin\left(4 \pi x^3 \right),& 0 < x < 1,  \\
	q_x(0,t) = u^{(1)}(t) = 12,& t > 0, \\
	q_x(1,t) = v^{(1)}(t) = 6\pi - 8,& t > 0.
	\end{dcases}
	\label{LS_numerical1_FI_N}
	\end{equation}	
	The dispersive nature of the LS equation is captured by the SD-UTM solution \eqref{class_soln_LS_centered_N} in the three $(x_n,t)$-plots of Figure \ref{LS_UTM1_FI_N}, while its first-order accuracy is presented in Figure \ref{LS_errorplot_UTM1_FI_N}.


\section{Computational Comparisons}\label{perform_comp}

	For each of the previous numerical examples, we compared the performance of the SD-UTM solutions to standard finite-difference schemes as $h \rightarrow 0$ for a fixed time step. In this section, we compare all numerical methods in terms of wall-clock time needed to achieve a target accuracy. We use the IBVP \eqref{LS_numerical1_FI} as an example, which poses the challenge of accurately capturing dispersive behavior as $T$ increases. 
	
\begin{figure}[tb]	
		\raggedleft
		\hspace*{-15pt}\begin{subfigure}[t]{.45\textwidth}
			\centering
  			\includegraphics[width=1.19\linewidth]{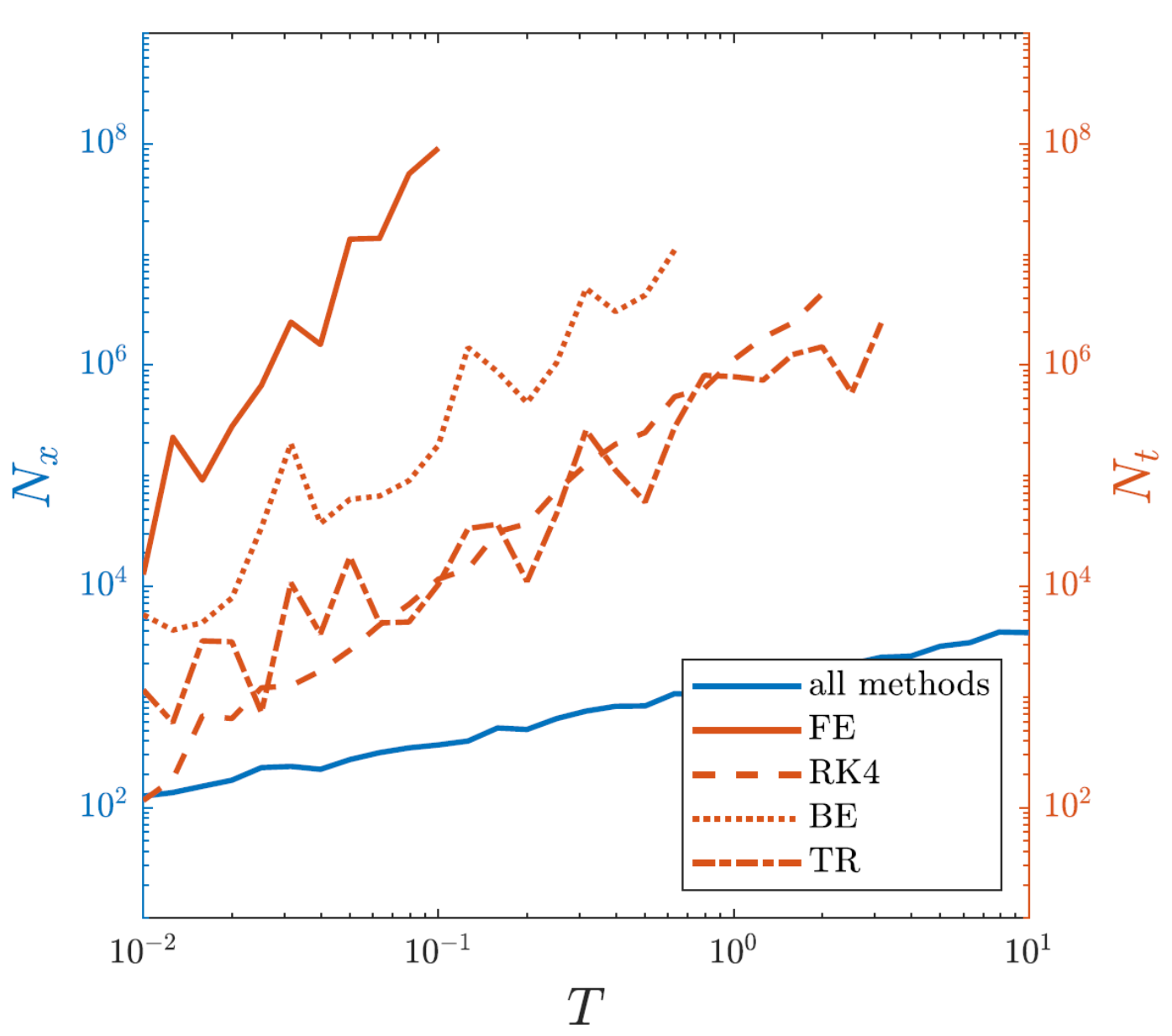}
  			\caption{}
  			\label{LS_Nx_Nt}
		\end{subfigure}\hfill 
		\begin{subfigure}[t]{.45\textwidth}
			\centering
  			\includegraphics[width=1.1\linewidth]{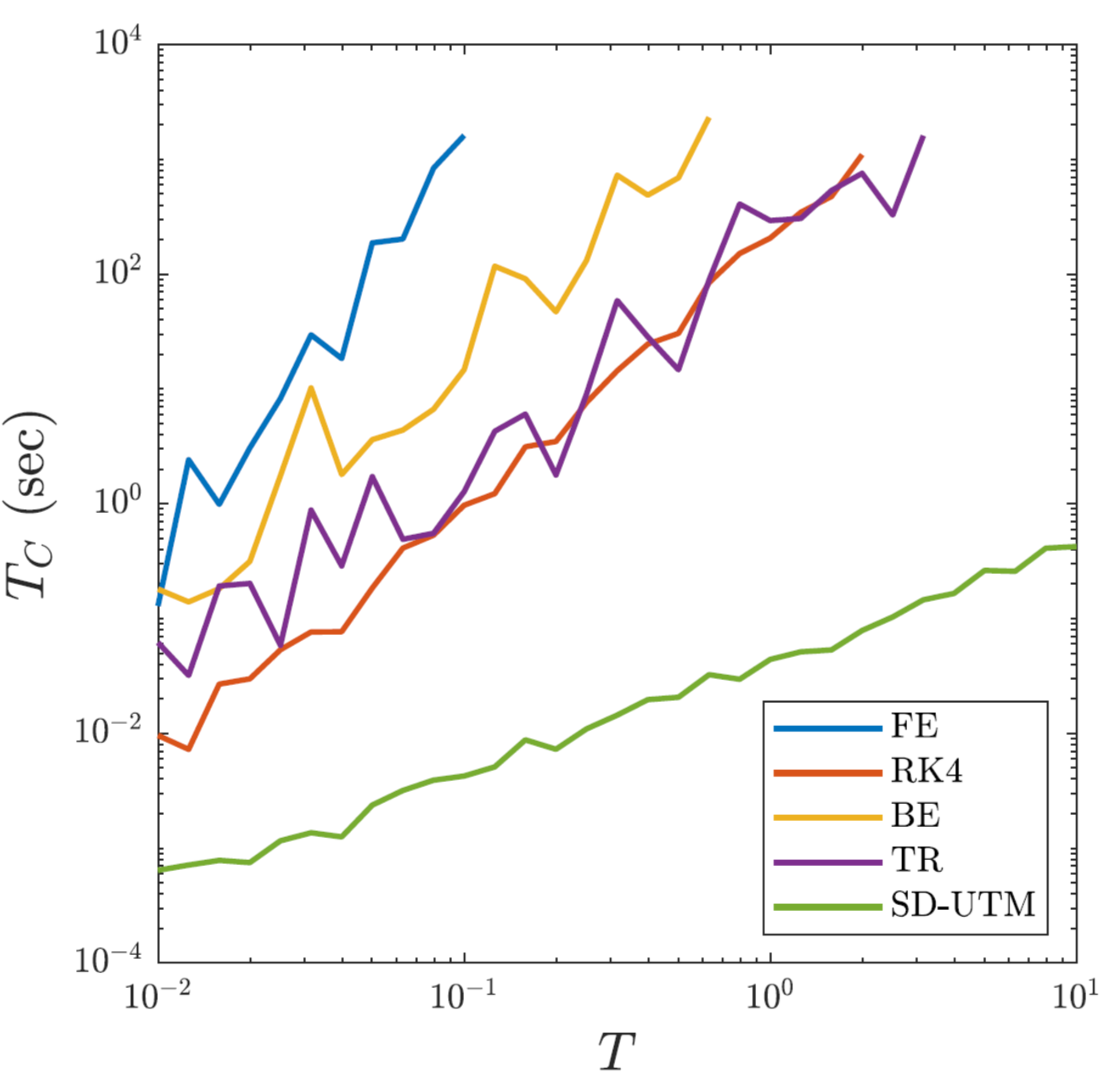}
  			\caption{}
  			\label{LS_Tcomp}
		\end{subfigure}	
		\caption{(a) The blue curve denotes the number of spatial grid points $N_x$ required for the SD-UTM solution to approximately reach the accuracy $E = 10^{-2}$ when solving the IBVP \eqref{LS_numerical1_FI}. Using this information, every finite-difference method uses the same spatial grid for each $T$ to determine how many time-steps $N_t$ are required to reach a similar accuracy as the SD-UTM solution. (b) The wall-clock computation time (averaged over 10 runs to rule out any effects due to background processes) required for each method to solve the IBVP \eqref{LS_numerical1_FI} with the selected $N_x$ and $N_t$ that approximately give an accuracy $E$.}
		\label{LS_comps}
	\end{figure}

	For all methods, we impose a mild target accuracy $\lVert q(x_n,T) - q_n(T) \rVert_{\infty} \approx E$ with $E = 10^{-2}$ for $10^{-2} \leq T \leq 10^1$. First, we determine the number of spatial grid points $N_x$ needed for the SD-UTM solution to reach the target accuracy. Using this spatial grid, we determine the number of time steps $N_t$ needed for each finite-difference solution to reach a similar accuracy. We use the SD-UTM series representation \eqref{class_soln_LS_centered_FI} with $f_0(W_{\ell},T) = 0$ and $g_0(W_{\ell},T) = 2\big(1 - e^{-W_{\ell}T}\big)/W_{\ell}$, while the finite difference solutions are set up in a standard method-of-lines approach with sparse tridiagonal linear systems that are efficiently solved. For every $T$, Figure \ref{LS_comps} presents $N_x$ and $N_t$ together (left panel) and the methods' wall-clock computation times $T_C$ (right panel) when solving the IBVP to the target accuracy $E$. The simulations were conducted in \MATLAB R2021b on an Intel Core i7-8705G processor with 12GB of RAM. 
	
	Starting from $T = 10^{-2}$, we stop computing the finite-difference solutions if $T_C > 10^3$ seconds for the most recent $T$. This threshold roughly translates to 91 million time steps for FE, 4 million time steps for RK4, 11 million time steps for BE, and 2 million time steps for TR, which we deem impractical, and terminate the computations to avoid excessive wall-clock times in trying to reach the final $T = 10$. All methods achieved an accuracy of $E \pm 10^{-4}$ for all successful $T$. Although the SD-UTM solution itself does not rely on a time-stepping procedure, we must refine the spatial mesh as $T$ increases to remove higher-order dispersive effects (see Section \ref{LS_D_section} for the modified PDE). Note that at the finest mesh of 3810 grid points, the computation time is still less than 1 second.
	


\section{Small-Time Increments}\label{small_time_sec}

	Our overall goal is to numerically solve IBVPs for quasilinear PDEs \eqref{ibvp_eq}, where the most nonlocal stencil is applied to the linear problem. Our proposed approach is a split-step technique: iteratively combine the solutions of the separated $M^{\text{th}}$-order linear IBVP $q_t = c \, q_{Mx}$ and the nonlinear IBVP $q_t = F\left(q,q_x,\ldots,q_{(M-1)x} \right)$ \cite{macnamara, randy_splitstep1, randy_splitstep2}. This requires the repeated computations of both linear and nonlinear solutions for time steps $\tau = T - t_0 \ll 1$, where $t_0$ is an arbitrary initial time. 
	
	For half-line problems \cite{SDUTM_HL}, we showed how to derive small-time approximate solutions with predetermined accuracy for SD-UTM integral representations. A similar approach is applied for finite-interval SD-UTM integral representations. For the forward-discretized advection equation \eqref{advec1_forward}, we show how to briefly derive these approximate solutions. 
	
	To include an arbitrary initial time $t_0$, we generalize the original IBVP \eqref{advec1_prob} to
	\begin{equation}\begin{dcases}
		q_t = c \, q_{x},& 0 < x < L,\, t > t_0, \\
		q(x,t_0) = \phi(x),& 0 < x < L,\\
		q(L,t) = v^{(0)}(t),& t > t_0,
		\label{advec1_prob_splitstep}
	\end{dcases}\end{equation}
	where $\phi(x)$ is the output from the previous split-step. Following the procedure in Section \ref{advec_forward_finiteinterval}, the solution to \eqref{advec1_prob_splitstep} is
	\begin{align}
		q_n(\tau) &= \frac{1}{2\pi} \int_{-\pi/h}^{\pi/h} e^{iknh} e^{-W \tau} \hat{q}(k,t_0)\,dk + \frac{c}{2\pi} \int_{-\pi/h}^{\pi/h} e^{ik(nh - L+h)} e^{-WT} g_{0}\,dk,
		\label{soln_advec_forward_splitstep}
	\end{align}
	where
	$$g_j\left(W,\tau\right) = \int_{0}^\tau e^{W(t + t_0)} q_{N+1+j}(t + t_0) \, dt, \quad k \in \mathbb{C},$$
	after making the substitution $\tau = T - t_0$. We expand \eqref{soln_advec_forward_splitstep} around $\tau=0$ to obtain a convenient approximation for a split-step method. We do not expand the initial-condition integral term, since the only time dependence is through $e^{-W \tau}$. For this term, we simply follow the steps in Section \ref{advec_forward_finiteinterval_series} to rewrite the integral as a series. The boundary-condition term of \eqref{soln_advec_forward_splitstep} has a more intricate dependence on time. Rewriting this integral into a series and then expanding the time-dependent terms leads to representations unique to the IBVP that generally cannot be addressed (see Remark \ref{ss_remark}). First, we expand the integrand and then rewrite each of the resulting expansion terms as finite series. Doing so up to third-order terms, we have
	\begin{align}
		q_n(\tau) &= e^{-c\tau/h} \sum_{m = 0}^{N-n} \left(\frac{c\tau}{h}\right)^{m} \frac{\phi_{n+m}}{m !}  \,+\, K_1(n) \tau \,+\, K_2(n) \tau^2 \,+\, K_3(n) \tau^3 \,+\, \mathcal{O}(\tau^4) , 
		\label{soln_splitstep1}
	\end{align}
	where $L = (N+1)h$ and
	\begin{align*}
		K_1(n) &= \frac{c}{2\pi} \int_{-\pi/h}^{\pi/h} e^{i k(n - N)h} v^{(0)}(t_0) \, dk = \frac{c v^{(0)}(t_0)}{h} \delta_{N n}, \\
		K_2(n) &= \frac{c}{4\pi} \int_{-\pi/h}^{\pi/h} e^{i k(n - N)h}\left[ \dot{v}^{(0)}(t_0) - W v^{(0)}(t_0)\right] \, dk,\\
		K_3(n) &= \frac{c}{12\pi} \int_{-\pi/h}^{\pi/h} e^{i k(n - N)h}\left[ \ddot{v}^{(0)}(t_0) - W \dot{v}^{(0)}(t_0) + W^2 v^{(0)}(t_0)\right] \, dk.
	\end{align*}
	The only dependence on $k$ within the integrands of $K_{\ell}(n)$ for $\ell = 1,2, \ldots$, are through $e^{i k(n - N)h}$ and powers of $W(k,\tau)$, and we want to rewrite
	$$I_m(n) = \int_{-\pi/h}^{\pi/h} e^{i k (n - N)h} W^m \, dk, \quad m = 0, 1, \dots, $$
	as a series in terms of
	$$K_{\ell}(n) = \frac{c}{2 \pi \cdot {\ell}!} \sum_{j = 0}^{{\ell}-1} I_{{\ell} - 1 - j}(n) \left. \frac{d^j v^{(0)}(t)}{dt^j}\right|_{t = t_0}.$$
	With $z = e^{ikh}$, the dispersion relation for this IBVP gives
	$$I_m(n) = \frac{1}{i h} \left(\frac{c}{h}\right)^m \oint_{|z| = 1} \frac{(1-z)^m}{z^{N-n+1}} \, dz = \frac{2\pi (-1)^{N-n}}{h (N - n)!} \left( \frac{c}{h}\right)^m \frac{m!}{(m-N+n)!}.$$
	Note that $I_m(n)$ is nonzero only for $m + n \geq N$. Substituting into $K_{\ell}(n)$ yields the small-time approximation \eqref{soln_splitstep1}, free of integral computations, where
	$$K_{\ell}(n) = \frac{(-1)^{N-n}}{(N-n)! \,{\ell}!} \sum_{j = 0}^{{\ell}-1} \frac{(-1)^j ({\ell} - 1 -j)!}{({\ell} -1 -j -N + n)!} \left(\frac{c}{h}\right)^{{\ell}-j} \left. \frac{d^j v^{(0)}(t)}{dt^j}\right|_{t = t_0}.$$	
	A similar process can be repeated for other IBVPs.
	
	\begin{remark}\label{ss_remark}
	Although the process of rewriting an SD-UTM integral representation into a series is straightforward, the resulting formulas are specific to the IBVPs and their dispersion relations. To obtain these small-time approximation solutions in general, we must start from the integral representations -- not the series.

	\end{remark}


\section{Concluding Remarks}

	Like the continuous UTM, the SD-UTM is applied algorithmically using the following steps:
\begin{enumerate}
	\item rewrite the semidiscretized equation into divergence form to obtain the local relation and the dispersion relation $W(k)$, \label{sd_utm_1}
	\item sum over spatial indices and integrate over the temporal domain to obtain the global relation, \label{sd_utm_2}
	\item invert to obtain a representation of the solution that depends on unknown boundary data, \label{sd_utm_3}
	\item determine the symmetries $\nu_j(k)$ of $W(k)$, \label{sd_utm_4}
	\item determine where the global relations with $k \rightarrow \nu_j(k)$ are valid in $\mathbb{C}$, \label{sd_utm_5}
	\item if necessary, deform integral paths of the boundary terms appropriately, \label{sd_utm_6}
	\item \textit{if necessary, determine additional boundary conditions from the PDE}, \label{sd_utm_7}
	\item \textit{appropriately discretize boundary conditions}, \label{sd_utm_8}
	\item solve for unknowns using global relations with $k \rightarrow \nu_j(k)$ and time transforms of discretized boundary conditions, and \label{sd_utm_9}
	\item check integral terms involving $\hat{q}(\nu_j,T)$ vanish, resulting in the solution representation depending only on known quantities. \label{sd_utm_10}
\end{enumerate}
	To rewrite the SD-UTM solution formulas as series, more suitable for numerical computations, substitute the definitions for the Fourier and time transforms and apply Cauchy's Theorem, where additional deformations may be necessary. Note that this rewrite may not always be possible. Third-order problems, like the linear Korteweg-de Vries equations $q_t = \pm q_{xxx}$, on the half-line and finite interval will be presented in a forthcoming paper \cite{SDUTM_kdv}.
	
	For a given discretization of a PDE, the global relation and its regions of validity under the symmetries $\nu_j(k)$ impose which stencils can be selected for derivative boundary conditions, as we saw with the higher-order discretizations and the Neumann problems. Similar to how the continuous UTM determines which types of boundary conditions result in a well-posed problem \cite{fokas_book}, ``natural'' discretizations reduce the variety of stencils down to those that are compatible with the IBVP \cite{SDUTM_HL}. We re-iterate that a natural discretization for a PDE is (i) of the same order as the spatial order of the PDE, (ii) not purely one sided (except for first-order problems), and (iii) the one that optimally aligns with the available boundary conditions.


\section{Acknowledgments}
	
	This work was supported by the Graduate Opportunities \& Minority Achievement Program Fellowship from the University of Washington and the Ford Foundation Predoctoral Fellowship (JC). Any opinions, findings, and conclusions or recommendations expressed in this material are those of the authors and do not necessarily reflect the views of the funding sources.


\vspace{15pt}
\nocite{*}
{\small
\bibliographystyle{abbrv}
\bibliography{main}
}

\end{document}